\documentclass[letterpaper, 12pt, oneside]{book}
    
    \usepackage[left=108pt, right=74pt, top=108pt, bottom=81pt, dvips, pdftex]{geometry}
    
    \usepackage{fancyhdr}
    \renewcommand{\bibname}{References}

    \usepackage{amsmath,amsthm, amsfonts,amssymb}
    \usepackage{setspace}
    \usepackage[linktocpage,bookmarksopen,bookmarksnumbered,%
                pdftitle={Patience Sorting and Its Generalizations},%
                pdfauthor={Isaiah Lankham},%
                pdfsubject={UC Davis Ph.D. Doctoral Thesis},%
                pdfkeywords={Patience Sorting, Barred Permutation Patterns}]{hyperref}
    \usepackage{pstricks,pst-plot}
    \usepackage[vcentermath]{youngtab}
    
    \newtheorem{Theorem}{Theorem}[section]
    \newtheorem{Proposition}[Theorem]{Proposition}
    \newtheorem{Lemma}[Theorem]{Lemma}
    \newtheorem{Corollary}[Theorem]{Corollary}
    \theoremstyle{definition}
        \newtheorem{Definition}[Theorem]{Definition}
        \newtheorem{Example}[Theorem]{Example}

        \newtheorem{Algorithm}[Theorem]{Algorithm}
        \newtheorem{CardGame}[Theorem]{Card Game}
        \newtheorem{Strategy}[Theorem]{Strategy}
        \newtheorem{Question}[Theorem]{Question}
    \theoremstyle{remark}
        \newtheorem{Remark}[Theorem]{Remark}

    \newenvironment{ReuseTheorem}[2]{\par\vspace{12pt}\noindent{\bf #1~\ref{#2}$'$.}\it}{\par\vspace{12pt}}
    \newenvironment{RestateTheorem}[2]{\par\vspace{12pt}\noindent{\bf #1~\ref{#2}.}\it}{\par\vspace{12pt}}

    \numberwithin{equation}{section}
    
    \newcommand{\done}{d_{1}}
    \newcommand{\dtwo}{d_{2}}
    \newcommand{\dell}{d_{\ell}}
    \newcommand{\dk}{d_{k}}
    
    \newcommand{\eone}{e_{1}}
    \newcommand{\etwo}{e_{2}}
    \newcommand{\eell}{e_{\ell}}
    
    \newcommand{\n}{\vspace{12pt}}
    
    \newcommand{\newchapter}[3]
    {
        \chapter[#2]{#3}
        \chaptermark{#1}
        \thispagestyle{myheadings}
    }

\begin{document}

    \pagenumbering{roman}
    \pagestyle{plain}
       
    %
    %

    \singlespacing

    ~\vspace{-0.75in} 
    \begin{center}
        
        \begin{huge}
            Patience Sorting and Its Generalizations
        \end{huge}\\\n
        By\\\n
        {\sc Isaiah Paul Lankham}\\
        B.S. (California State University, Chico) 2002\\
        M.A. (University of California, Davis) 2004\\\n
        DISSERTATION\\\n
        Submitted in partial satisfaction of the requirements for the degree of\\\n
        DOCTOR OF PHILOSOPHY\\\n
        in\\\n
        MATHEMATICS\\\n
        in the\\\n
        OFFICE OF GRADUATE STUDIES\\\n
        of the\\\n
        UNIVERSITY OF CALIFORNIA\\\n
        DAVIS\\\n\n
        
        Approved:\\\n\n
        
%
%

        Craig A. Tracy (Chair)\\
        \rule{4in}{1pt}\\\n\n
        
        Eric M. Rains\\
        \rule{4in}{1pt}\\\n\n
        
        Jes\'{u}s A. De Loera\\
        \rule{4in}{1pt}\\

        \vfill
        
        Committee in Charge\\
        2007

    \end{center}

    \newpage

%
%
%

    %
    %

    \begin{quote}
        What you see and hear depends a great deal on where you are
        standing; it also depends on what sort of person you are.
        
        ~~~~---~C.~S.~Lewis, \emph{The Magician's Nephew}
    \end{quote}
    
    \newpage

    \doublespacing

    %
    %
    
    {\singlespacing
    \tableofcontents
    }
    
    \newpage

    %
    %

    ~\vspace{-1in} 
    \begin{flushright}
        \singlespacing
        Isaiah Paul Lankham\\
        June 2007\\
        Mathematics
    \end{flushright}
    
    \vspace{-0.125in}

    \centerline{\large Patience Sorting and Its Generalizations}
    
    \centerline{\textbf{\underline{Abstract}}}
            
        Despite having been introduced in the 1960s, the card game and
        combinatorial algorithm \emph{Patience Sorting} is only now
        beginning to receive significant attention.  This is due in
        part to recent results like the Baik-Deift-Johansson Theorem,
        which suggest connections with Probabilistic Combinatorics and
        Random Matrix Theory.

        Patience Sorting (a.k.a.~\emph{Floyd's Game}) can be viewed as
        an idealized model for the immensely popular single-person
        card game Klondike Solitaire.  Klondike is interesting from a
        mathematical perspective as it has long resisted the analysis
        of optimality for its strategies.  While there is a well-known
        optimal greedy strategy for Floyd's Game, we provide a
        detailed analysis of this strategy as well as a suitable
        adaption for studying more Klondike-like generalizations of
        Patience Sorting.

        At the same time, Patience Sorting can also be viewed as an
        iterated, non-recursive form of the Schensted Insertion
        Algorithm.  We study the combinatorial objects that result
        from this viewpoint and then extend Patience Sorting to a full
        bijection between the symmetric group and certain pairs of
        these combinatorial objects.  This \emph{Extended Patience
        Sorting Algorithm} is similar to the Robinson-Schensted-Knuth
        (or RSK) Correspondence, which is itself built from repeated
        application of the Schensted Insertion Algorithm.

        This analysis of Patience Sorting and its generalizations
        naturally encounters the language of barred pattern avoidance.
        We also introduce a geometric form for the Extended Patience
        Sorting Algorithm that is naturally dual to X.~G.~Viennot's
        celebrated Geometric RSK Algorithm.  Unlike Geometric RSK,
        though, the lattice paths coming from Patience Sorting are
        allowed to intersect.  We thus also give a characterization
        for the intersections of these lattice paths.

    \newpage

    %
    %

    \chapter*{\vspace{-1.5in}Acknowledgments and Thanks}

        First of all, I would like to take this opportunity to thank
        Craig Tracy for his eternally patient style of advising and
        his relentless ability to provide encouragement.  Craig has
        been incredibly generous with his time and ideas, and, without
        the direction and advice I have received from him, this
        Dissertation would not have been possible.
        
        I would next like to thank the other two members of my
        Dissertation Committee, Jes\'{u}s De Loera and Eric Rains.
        Their input and encouragement have also been an integral part
        of the Dissertation writing process, and I am especially
        grateful for the many mathematical conversations that I have
        been fortunate to share with them both.
        
        Other members of the UC Davis mathematics community that
        deserve extra special thanks include Celia Davis, Philip
        Sternberg, and Brian Wissman.  For as long as I've known her,
        Celia has been an indispensable source of professional
        assistance and emotional support.  Philip has been the best
        officemate that one could hope to share three different
        offices with over the course of a four and a half year period.
        Brian has been someone with whom I have enjoyed, and hope to
        continue enjoying, many wonderful conversations regarding
        teaching and pedagogical style.
        
        While there are many other people for whom thanks is due, I
        would like to close by extending extra special thanks to my
        wife, Kelly.  She has been kind enough to keep me company
        during many late nights of writing and revisions, and, without
        her support, none of my work as a graduate student would have
        been at all possible.\\

        \noindent This work was supported in part by the U.S. National
        Science Foundation under Grants DMS-0135345, DMS-0304414, and
        DMS-0553379.
    
    \newpage

    %
    %

    \pagestyle{fancy}
    \pagenumbering{arabic}
       
    %
    %

    \newchapter{Introduction}{Introduction}{Introduction}
    \label{sec:Intro}
    
        \section{Informal Overview and Motivation}
        \label{sec:Intro:Motivation}

                Given a positive integer $n \in \mathbb{Z}_{+}$, we
                use $\mathfrak{S}_{n}$ to denote the \emph{symmetric
                group} on the set $[n] = \{ 1, 2, \ldots, n \}$.  In
                other words, $\mathfrak{S}_{n}$ is the set of all
                bijective functions on $[n]$.  Each element $\sigma
                \in \mathfrak{S}_{n}$ is called a \emph{permutation},
                and $\sigma_{i} = \sigma(i)$ denotes the
                $i^{\text{th}}$ function value for each $i \in [n]$.
                However, even though $\sigma$ is defined as a
                \emph{function} on the set $[n]$, it is often
                convenient to instead regard $\sigma$ as a
                \emph{rearrangement} of the sequence of numbers $1, 2,
                \ldots, n$.  This, in particular, motives the
                so-called \emph{two-line notation}
                
                \vspace{-12pt}
                {\singlespacing
                \[
                    \sigma
                    =
                    \begin{pmatrix}
                        1 & 2 & \cdots & n \\
                        \sigma_{1} & \sigma_{2} & \cdots & \sigma_{n}
                    \end{pmatrix}
                \]
                }
                \noindent and its associated \emph{one-line notation}
                $\sigma = \sigma_{1}\sigma_{2}\cdots\sigma_{n}$.
                
                In keeping with this emphasis on the values
                $\sigma_{i}$, a \emph{subsequence}
                (a.k.a.~\emph{subpermutation}) of a permutation
                $\sigma = \sigma_{1}\sigma_{2}\cdots\sigma_{n}$ is any
                sequence of the form $\pi =
                \sigma_{i_{1}}\sigma_{i_{2}}\cdots\sigma_{i_{k}}$,
                where $k \in [n]$ and $i_{1} < i_{2} < \cdots <
                i_{k}$.  We denote the \emph{length} of the
                subsequence $\pi$ by $|\pi| = k$.  It is also common
                to call $\pi$ a \emph{partial permutation on} $[n]$
                since it is the restriction of the bijective function
                $\sigma$ to the subset $\{i_{1}, i_{2}, \ldots,
                i_{k}\}$ of $[n]$.  Note that the components of
                $\sigma_{i_{1}}\sigma_{i_{2}}\cdots\sigma_{i_{k}}$ are
                not required to be contiguous in $\sigma$.  As such,
                subsequences are sometimes called \emph{scattered
                subsequences} in order to distinguish them from
                so-called (contiguous) \emph{substrings}.

            \subsection{Longest Increasing Subsequences and Row Bumping}
            \label{sec:Intro:Motivation:LISs}

                Given two permutations (thought of as rearrangement of
                the numbers $1, 2, \ldots, n$), it is natural to ask
                whether one permutation is more ``out of order'' than
                the other when compared with the strictly increasing
                arrangement $12\cdots n$.  For example, one would
                probably consider $\sigma^{(1)} = 53241$ to be more
                ``out of order'' than something like $\sigma^{(2)} =
                21354$.  While there are various metrics for
                justifying such intuitive notions of ``disorder'', one
                of the most well-studied involves the examination of a
                permutation's subsequences that are themselves in
                strictly increasing order.
                
                An \emph{increasing subsequence} of a permutation is
                any subsequence that increases when read from left to
                right.  For example, $\mathbf{234}$ is an increasing
                subsequence of $\mathbf{2}1\mathbf{3}5\mathbf{4}$.
                One can also see that $\mathbf{235}$ is an increasing
                subsequence of $\mathbf{2}1\mathbf{35}4$.  This
                illustrates the nonuniqueness of \emph{longest
                increasing subsequences}, and such subsequences can
                even be disjoint as in $456123$.  The \emph{length} of
                every longest increasing subsequence is nonetheless a
                well-defined property for a given permutation $\sigma
                \in \mathfrak{S}_{n}$, and we denote this statistic by
                $\ell_{n}(\sigma)$.  For example, with notation as
                above, $\ell_{5}(\sigma^{(1)}) = \ell_{5}(53241) = 2$
                and $\ell_{5}(\sigma^{(2)}) = \ell_{5}(21354) = 3$,
                which provides one possible heuristic justification
                for regarding $\sigma^{(1)}$ as more ``out of order''
                than $\sigma^{(2)}$.
                
                Given a permutation $\sigma \in \mathfrak{S}_{n}$,
                there are various methods for calculating the length
                of the longest increasing subsequence
                $\ell_{n}(\sigma)$.  The most obvious algorithm
                involves directly examining every subsequence of
                $\sigma$, but such an approach is far from being
                computationally efficient as there are $2^{n}$ total
                subsequences to examine.  Since a given increasing
                subsequence is essentially built up from shorter
                increasing subsequences, a fairly routine application
                of so-called dynamic programming methodologies allows
                us to calculate $\ell_{n}(\sigma)$ using, in the worst
                possible case, $O(n^{2})$ operations.
                
                \begin{Algorithm}[Calculating $\ell_{n}(\sigma)$ via Dynamic Programming]
                \label{alg:DynamicProgrammingLISs}
                {\tt
                    ~\\
                    Input:~a permutation $\sigma =
                    \sigma_{1}\sigma_{2}\cdots\sigma_{n} \in
                    \mathfrak{S}_{n}$\\
                    Output:~the sequence of positive
                    integers $L_{1}(\sigma), \ldots, L_{n}(\sigma)$
                    
                    \begin{enumerate}
                        
                        \item First set $L_{1}(\sigma) = 1$.
                        
                        \item Then, for each $i = 2, \ldots, n$,
                        determine the value $L_{i}(\sigma)$ as
                        follows:
                        
                        \begin{enumerate}
                            
                            \item If ${\displaystyle \sigma_{i} > 
                            \min_{1 \, \leq \, j \, < \,
                            i}\{ \sigma_{j} \}}$, then set
                            \[
                                L_{i}(\sigma) = 1 + \max_{1 \, \leq \,
                                j \, < \, i}\{ L_{j}(\sigma) \ | \
                                \sigma_{i} > \sigma_{j} \}.
                            \]
                            
                            \item Otherwise, set $L_{i}(\sigma) = 1$.
                            
                        \end{enumerate}
                    
                    \end{enumerate}
    
                }
                \end{Algorithm}
                
                \noindent Each value $L_{i}(\sigma)$ computed in
                Algorithm~\ref{alg:DynamicProgrammingLISs} is the
                length of the longest increasing subsequence (when
                reading from left to right) that is terminated by the
                entry $\sigma_{i}$ in $\sigma$.  Given this data, it
                is then clear that
                \[
                    \ell_{n}(\sigma)
                    =
                    \max_{1 \, \leq \, i \, \leq \, n}\{ L_{i}(\sigma) \}.
                \]
                We illustrate
                Algorithm~\ref{alg:DynamicProgrammingLISs} in the
                following example.
                
                \begin{Example}
                \label{eg:DynamicProgrammingLISsAlgExample}
                    Given $\sigma =
                    \sigma_{1}\sigma_{2}\cdots\sigma_{9} = 364827159
                    \in \mathfrak{S}_{9}$, we use
                    Algorithm~\ref{alg:DynamicProgrammingLISs} to
                    compute the sequence $L_{1}(\sigma),
                    L_{2}(\sigma), \ldots, L_{9}(\sigma)$ by
                    
                    \begin{itemize}
                        \item first setting $L_{1}(\sigma) = 1$, and 
                        then,
                    
                        \item for each $i = 2, \ldots, 9$, computing the
                        values $L_{i}(\sigma)$ as follows:

                            \begin{itemize}
                            
                                \item Since $6 = \sigma_{2} > \min\{
                                \sigma_{1} \} = 3$, set $L_{2}(\sigma)
                                = 1 + L_{1}(\sigma) = 2$.
                            
                                \item Since $4 = \sigma_{3} > \min\{
                                \sigma_{1}, \sigma_{2} \} = 3$, set
                                $L_{3}(\sigma) = 1 + L_{1}(\sigma) =
                                2$.
                            
                                \item Since $8 = \sigma_{4} > \min\{
                                \sigma_{1}, \sigma_{2}, \sigma_{3} \}
                                = 3$, set $L_{4}(\sigma) = 1 +
                                L_{3}(\sigma) = 3$.
                            
                                \item Since $2 = \sigma_{5} < \min\{
                                \sigma_{1}, \ldots, \sigma_{4} \}
                                = 2$, set $L_{5}(\sigma) = 1$.
                            
                                \item Since $7 = \sigma_{6} > \min\{
                                \sigma_{1}, \ldots, \sigma_{5} \}
                                = 2$, set $L_{6}(\sigma) = 1 +
                                L_{3}(\sigma) = 3$.
                            
                                \item Since $1 = \sigma_{7} < \min\{
                                \sigma_{1}, \ldots, \sigma_{6} \}
                                = 2$, set $L_{7}(\sigma) = 1$.
                            
                                \item Since $5 = \sigma_{8} > \min\{
                                \sigma_{1}, \ldots, \sigma_{7} \}
                                = 2$, set $L_{8}(\sigma) = 1 +
                                L_{3}(\sigma) = 3$.
                            
                                \item Since $9 = \sigma_{9} > \min\{
                                \sigma_{1}, \ldots, \sigma_{8} \}
                                = 2$, set $L_{9}(\sigma) = 1 +
                                L_{8}(\sigma) = 4$.
                                
                            \end{itemize}
                            
                    \end{itemize}
                    
                    \noindent It follows that $\ell_{9}(\sigma) =
                    \ell_{9}(364827159) = \max\{ L_{i}(\sigma) \ | \ i
                    = 1, 2, \ldots, 9 \} = 4$, which can be checked by
                    direct inspection.  E.g., $3679$ and $3459$ are
                    two longest increasing subsequences in
                    $364827159$.
                    
                \end{Example}
                
                While each term in the sequence $L_{1}(\sigma),
                L_{2}(\sigma), \ldots, L_{n}(\sigma)$ has significance
                in describing various combinatorial properties of the
                permutation $\sigma \in \mathfrak{S}_{n}$, there is no
                need to explicitly calculate every value if one is
                only interested in finding the length
                $\ell_{n}(\sigma)$ of the longest increasing
                subsequence in $\sigma$.  To see this, suppose that
                the value $v \in \mathbb{Z}_{+}$ occurs in the
                sequence $L_{1}(\sigma), L_{2}(\sigma), \ldots,
                L_{n}(\sigma)$ at positions $i_{1}, i_{2}, \ldots,
                i_{k} \in [n]$, where $i_{1} < i_{2} < \cdots <
                i_{k}$.  Then, from the definitions of
                $L_{i_{1}}(\sigma), L_{i_{2}}(\sigma), \ldots,
                L_{i_{k}}(\sigma)$, we must have that $\sigma_{i_{1}}
                > \sigma_{i_{2}} > \cdots > \sigma_{i_{k}}$.  (In
                other words,
                $\sigma_{i_{1}}\sigma_{i_{2}}\cdots\sigma_{i_{k}}$ is
                a decreasing subsequence of $\sigma$, which we call
                the $v^{\rm th}$ \emph{left-to-right minima
                subsequences} of $\sigma$.  Such subsequences will
                play an important role in the analysis of Patience
                Sorting throughout this Dissertation.  See
                Section~\ref{sec:PSasAlgorithm:PilesAndNEshadows})
                Moreover, given a particular element $\sigma_{i_{j}}$
                in the subsequence
                $\sigma_{i_{1}}\sigma_{i_{2}}\cdots\sigma_{i_{k}}$,
                $L_{m}(\sigma) > v$ for some $m > i_{j}$ if and only
                if $\sigma_{m} > \sigma_{i_{j}}$.  Consequently, it
                suffices to solely keep track of the element
                $\sigma_{i_{j}}$ in the subsequence
                $\sigma_{i_{1}}\sigma_{i_{2}}\cdots\sigma_{i_{k}}$ in
                order to determine the value of $L_{i_{j +
                1}}(\sigma)$.  This observation allows us to
                significantly reduce the number of steps required for
                computing $\ell_{n}(\sigma)$.  The resulting
                \emph{Single Row Bumping Algorithm} was first
                implicitly introduced by Craige Schensted
                \cite{refSchensted1961} in 1961 and then made explicit
                by Knuth \cite{refKnuth1970} in 1970.
                
                \begin{Algorithm}[Single Row Bumping]
                \label{alg:SingleRowBumping}
                {\tt
                    ~\\
                    Input:~a permutation $\sigma =
                    \sigma_{1}\sigma_{2}\cdots\sigma_{n} \in
                    \mathfrak{S}_{n}$\\
                    Output:~the partial permutation $w$
                    
                    \begin{enumerate}
                        
                        \item First set $w = \sigma_{1}$.
                        
                        \item Then, for each $i = 2, \ldots, n$,
                        insert $\sigma_{i}$ into $w =
                        w_{1}w_{2}\cdots w_{k}$ using the following 
                        rule:
                        
                        \begin{enumerate}
                            
                            \item If $\sigma_{i} > w_{k}$, then append
                            $\sigma_{i}$ to $w$ to obtain $w =
                            w_{1}w_{2}\cdots w_{k}\sigma_{i}$.
                            
                            \item Otherwise, use $\sigma_{i}$ to
                            "bump" the left-most element $w_{j}$ of
                            $w$\\ that is larger than $\sigma_{i}$.  In
                            other words, set
                            \[
                                w = w_{1}w_{2}\cdots
                                w_{j-1}\sigma_{i}w_{j+1}\cdots w_{k} \
                                \ \mathtt{where} \ \ j = \min_{1 \, 
                                \leq \,
                                m \, \leq \, k} \{ m \ | \ \sigma_{i} <
                                \sigma_{m} \}.
                            \]
                            
                        \end{enumerate}
                        
                    \end{enumerate}
                }
                \end{Algorithm}
                
                \noindent In particular, one can show that
                $L_{i}(\sigma) \geq r$ in
                Algorithm~\ref{alg:DynamicProgrammingLISs} if and only
                if $\sigma_{i}$ was inserted into $w$ at position $r$
                during some iteration of
                Algorithm~\ref{alg:SingleRowBumping}.  It follows that
                $\ell_{n}(\sigma)$ is equal to the length $|w|$ of
                $w$.  (The interplay between the sequence
                $L_{1}(\sigma), \ldots, L_{n}(\sigma)$ and the
                formation of $w$ can most easily be seen by comparing
                Example~\ref{eg:DynamicProgrammingLISsAlgExample} with
                Example~\ref{eg:SingleRowBumpingExample} below.)
                
                In terms of improved computational efficiency over
                Algorithm~\ref{alg:DynamicProgrammingLISs}, note that
                the elements in $w$ must necessarily increase when
                read from left to right.  Consequently, Step 2(b) of
                Algorithm~\ref{alg:SingleRowBumping} can be
                accomplished with a Binary Search (see
                \cite{refCLRS2002}), under which Single Row Bumping
                requires, in the worst possible case, $O(n\log(n))$
                operations in order to calculate $\ell_{n}(\sigma)$
                for a given permutation $\sigma \in \mathfrak{S}_{n}$.
                This can actually be even further reduced to
                $O(n\log(\log(n)))$ operations if $w$ is formed as a
                special type of associative array known as a \emph{van
                Emde Boas tree} \cite{refBS2000}; see
                \cite{refvEBKZ1977} for the appropriate definitions.
                (Hunt and Szymanski \cite{refHS1977} also
                independently gave an $O(n\log(\log(n)))$ algorithm in
                1977.  Their algorithm, however, computes
                $\ell_{n}(\sigma)$ as a special case of the length of
                the longest common subsequence in two permutations,
                where one of the two permutations is taken to be $12
                \cdots n$.)
                
                In the following example, we illustrate Single Row
                Bumping (Algorithm~\ref{alg:SingleRowBumping}) using a
                row of boxes for reasons that will become clear in
                Section~\ref{sec:Intro:Motivation:RSK} below.  We also
                explicitly indicate each ``bumped'' value, using the
                null symbol ``$\emptyset$'' to denote the empty
                partial permutation.  (E.g., ``$\young(36) \leftarrow
                \mathbf{4} = \young(34) \leadsto 6$'' means that
                $\mathbf{4}$ has been inserted into $\young(36)$ and
                has ``bumped'' $6$, whereas ``$\young(34) \leftarrow
                \mathbf{8} = \young(348) \leadsto \emptyset$'' means
                that $\mathbf{8}$ has been inserted into \young(34)
                and nothing has been ``bumped''.)

                \begin{Example}
                \label{eg:SingleRowBumpingExample}
                    Given $\sigma = 364827159 \in \mathfrak{S}_{9}$,
                    we use Single Row Bumping
                    (Algorithm~\ref{alg:SingleRowBumping}) to form the
                    partial permutation $w$ as follows:
                    
                    \begin{itemize}
                        \item Start with $w = \emptyset$, and insert
                        $\mathbf{3}$ to obtain $w = \emptyset
                        \leftarrow \mathbf{3} = \young(3)$.
                    
                        \item Append $\mathbf{6}$ so that $w =
                        \young(3) \leftarrow \mathbf{6} = \young(36) 
                        \leadsto \emptyset$.
                    
                        \item Use $\mathbf{4}$ to bump $6$: $w =
                        \young(36) \leftarrow \mathbf{4} = 
                        \young(34) \leadsto 6$.
                    
                        \item Append $\mathbf{8}$ so that $w =
                        \young(34) \leftarrow \mathbf{8} = \young(348) 
                        \leadsto \emptyset$.
                        
                        \item Use $\mathbf{2}$ to bump $3$: $w =
                        \young(348) \leftarrow \mathbf{2} = 
                        \young(248) \leadsto 3$.
                    
                        \item Use $\mathbf{7}$ to bump $8$: $w =
                        \young(248) \leftarrow \mathbf{7} = 
                        \young(247) \leadsto 8$.
                    
                        \item Use $\mathbf{1}$ to bump $2$: $w =
                        \young(247) \leftarrow \mathbf{1} = 
                        \young(147) \leadsto 2$.
                    
                        \item Use $\mathbf{5}$ to bump $7$: $w =
                        \young(147) \leftarrow \mathbf{5} = 
                        \young(145) \leadsto 7$.
                    
                        \item Append $\mathbf{9}$ so that $w =
                        \young(145) \leftarrow \mathbf{9} =
                        \young(1459) \leadsto \emptyset$.
                    \end{itemize}
                    
                    \noindent Even though $w$ is not a longest
                    increasing subsequence of $\sigma$, one can
                    nonetheless check that $\ell_{9}(\sigma) = 4 =
                    |w|$.
                    (Cf.~Example~\ref{eg:DynamicProgrammingLISsAlgExample}.)
                    
                \end{Example}
                
                Algorithm~\ref{alg:SingleRowBumping} provides an
                efficient method for computing the length of the
                longest increasing subsequence statistic, but it is
                also combinatorially wasteful.  In particular, it is
                reasonable to anticipate that even more combinatorial
                information would be obtained by placing additional
                structure upon the ``bumped'' values.  The most
                classical and well-studied generalization involves
                recursively reapplying
                Algorithm~\ref{alg:SingleRowBumping} in order to
                create additional partial permutations from the
                ``bumped'' values.  The resulting construction is
                called the \emph{Schensted Insertion Algorithm}.
                Historically, this was the original framework within
                which Schensted \cite{refSchensted1961} invented Row
                Bumping while studying the length of the longest
                increasing subsequence statistic.  We review Schensted
                Insertion and some of its more fundamental properties
                (including the well-studied and widely generalized
                \emph{RSK Correspondence} based upon it) in
                Section~\ref{sec:Intro:Motivation:RSK} below.
                
                The remainder of this Dissertation then describes
                various parallels and differences between Schensted
                Insertion and another natural extension of Single Row
                Bumping called \emph{Patience Sorting}.  We first
                describe Patience Sorting in
                Section~\ref{sec:Intro:Motivation:PS}.  Further
                background material on permutation patterns is then
                given in Section~\ref{sec:Intro:Motivation:Patterns}.
                We then provide a summary of the main results of this
                Dissertation in Section~\ref{sec:Intro:Summary}.
            
            \subsection[Schensted Insertion and the RSK Correspondence]{Schensted Insertion and the RSK Correspondence}
            \label{sec:Intro:Motivation:RSK}
                As discussed in
                Section~\ref{sec:Intro:Motivation:LISs} above, Single
                Row Bumping (Algorithm~\ref{alg:SingleRowBumping}) can
                be viewed as combinatorially wasteful since nothing is
                done with the values as they are ``bumped''.  The most
                classical extension repeatedly employs Single Row
                Bumping in order to construct a collection of partial
                permutations $P = P(\sigma) = (w_{1}, w_{2}, \ldots,
                w_{r})$.  This results in the following algorithm, in
                which we use the same ``$\leadsto$'' notation for
                ``bumping'' as in
                Example~\ref{eg:SingleRowBumpingExample}.
                
                \begin{Algorithm}[Schensted Insertion]
                \label{alg:SchenstedInsertion}
                {\tt
                    ~\\
                    Input:~a permutation $\sigma =
                    \sigma_{1}\sigma_{2}\cdots\sigma_{n} \in
                    \mathfrak{S}_{n}$\\
                    Output:~the sequence of
                    partial permutations $w_{1}, w_{2}, \ldots, w_{r}$
                    
                    \begin{enumerate}
                        
                        \item First set $w_{1} = \sigma_{1}$.
                        
                        \item Then, for each $i = 2, \ldots, n$,
                        insert $\sigma_{i}$ into the partial
                        permutations $w_{1}, w_{2}, \ldots, w_{m}$ as
                        follows:
                        
                        \begin{enumerate}
                            
                            \item Insert $\sigma_{i}$ into $w_{1}$
                            using Single Row Bumping
                            (Algorithm~\ref{alg:SingleRowBumping}).\\
                            If a value is "bumped", then denote it by
                            $\sigma_{1}^{*}$.  Otherwise, set\\
                            $\sigma_{1}^{*} = \emptyset$.  We denote
                            this redefinition of $w_{1}$ as
                            \[
                                w_{1} \leftarrow \sigma_{i}
                                =
                                (w_{1} \leftarrow \sigma_{i}) \leadsto \sigma_{1}^{*}.
                            \]
                            
                            \item For $j = 2, \ldots, m$, redefine
                            each $w_{j}$ as
                            \[
                                w_{j} \leftarrow \sigma_{j-1}^{*}
                                =
                                (w_{j} \leftarrow \sigma_{j-1}^{*}) 
                                \leadsto \sigma_{j}^{*}
                            \]
                            using the convention that $w_{j} = w_{j}
                            \leftarrow \emptyset = (w_{j} \leftarrow
                            \emptyset) \leadsto \emptyset = w_{j}
                            \leadsto \emptyset$.
                            
                        \end{enumerate}
                        
                    \end{enumerate}
                }
                \end{Algorithm}
                
                \noindent In other words, one forms $w = w_{1}$ as
                usual under Single Row Bumping
                (Algorithm~\ref{alg:SingleRowBumping}) while
                simultaneously forming a new partial permutation
                $w_{2}$ (also using Single Row Bumping) from each
                value as it is ``bumped'' from $w_{1}$.  The values
                ``bumped'' from $w_{2}$ are furthermore employed to
                form a new partial permutation $w_{3}$ (again using
                Single Row Bumping), and so forth until $w_{1},
                \ldots, w_{r}$ have been formed.
                
                We illustrate Schensted Insertion in the following
                example, where we form the object $P = (w_{1}, w_{2},
                \ldots, w_{r})$ as a collection of left- and
                top-justified boxes.  This extends the notation used
                in Example~\ref{eg:SingleRowBumpingExample} so that
                $P$ becomes a so-called \emph{standard Young tableau}
                written in \emph{English notation}.
                
                \begin{Example}
                \label{eg:SchenstedInsertionExample}
                    Given $\sigma = 364827159 \in \mathfrak{S}_{9}$,
                    we form $P(\sigma) = (w_{1}, w_{2}, w_{3}, w_{4})$
                    under Schensted Insertion
                    (Algorithm~\ref{alg:SchenstedInsertion}) as
                    follows:
                    
                    {\singlespacing
                    \begin{itemize}
                        \item Start with $P = \emptyset$ and insert 
                        $\mathbf{3}$ into $P$ to obtain
                        \[
                            P
                            =
                            \emptyset \leftarrow \mathbf{3}
                            =
                            \young(3).
                        \]
                    
                        \item Append $\mathbf{6}$ to obtain
                        \[
                            P
                            =
                            \young(3) \leftarrow \mathbf{6}
                            =
                            \young(36) \leadsto \emptyset.
                        \]
                    
                        \item Use $\mathbf{4}$ to bump $6$:
                        \[
                            P
                            =
                            \young(36) \leftarrow \mathbf{4}
                            =
                            \begin{array}{ll}
                                \young(34) & \leadsto 6
                                \\
                                \ \emptyset & \leftarrow 6
                            \end{array}
                            =
                            \ \young(34,6).
                        \]
                    
                        \item Append $\mathbf{8}$ to obtain
                        \[
                            P
                            =
                            \young(34,6) \raisebox{6pt}{$\ \leftarrow\mathbf{8}$}
                            =
                            \young(348,6).
                        \]
                        
                        \item Use $\mathbf{2}$ to bump $3$:
                        \[
                            P
                            =
                            \young(348,6) \raisebox{6pt}{$~\leftarrow\mathbf{2}$}
                            = 
                            \young(248,6)
                            \begin{array}{l}
                                \leadsto 3
                                \\
                                \leftarrow 3
                            \end{array}
                            =
                            \begin{array}{ll}
                                \young(248,3) & \raisebox{-8pt}{$\leadsto 6$}
                                \\
                                \ \emptyset & \leftarrow 6
                            \end{array}
                            =
                            \ \young(248,3,6).
                        \]
                    
                        \item Use $\mathbf{7}$ to bump $8$:
                        \[
                            P
                            =
                            \young(248,3,6) \raisebox{12pt}{$~\leftarrow\mathbf{7}$}
                            = 
                            \young(247,3,6)
                            \begin{array}{l}
                                \leadsto 8
                                \\
                                \leftarrow 8
                                \\
                                ~
                            \end{array}
                            =
                            \ \young(247,38,6).
                        \]
                    
                        \item Use $\mathbf{1}$ to bump $2$:
                        \[
                            \begin{array}{rclcrcl}
                                P
                                &
                                =
                                &
                                \young(247,38,6)
                                \raisebox{12pt}{$~\leftarrow\mathbf{1}$}
                                &
                                =
                                &
                                \young(147,38,6)
                                \begin{array}{l}
                                    \leadsto 2
                                    \\
                                    \leftarrow 2
                                    \\
                                    ~
                                \end{array}
                                &
                                =
                                &
                                \young(147,28,6)
                                \begin{array}{l}
                                    ~
                                    \\
                                    \leadsto 3
                                    \\
                                    \leftarrow 3
                                \end{array}
                                \\
                                \\
                                &
                                &
                                &
                                =
                                &
                                \begin{array}{ll}
                                    \young(147,28,3) & \raisebox{-14pt}{$\leadsto 6$}
                                    \\
                                    \ \emptyset & \leftarrow 6
                                \end{array}
                                &
                                =
                                &
                                \young(147,28,3,6).
                            \end{array}
                        \]
                    
                        \item Use $\mathbf{5}$ to bump $7$:
                        \[
                            P
                            =
                            \young(147,28,3,6) \raisebox{20pt}{$~\leftarrow\mathbf{5}$}
                            =
                            \young(145,28,3,6)
                            \begin{array}{l}
                                \leadsto 7
                                \\
                                \leftarrow 7
                                \\
                                ~
                                \\
                                ~
                            \end{array}
                            =
                            \young(145,27,3,6)
                            \begin{array}{l}
                                ~
                                \\
                                \leadsto 8
                                \\
                                \leftarrow 8
                                \\
                                ~
                            \end{array}
                            =
                            \young(145,27,38,6).
                        \]
                    
                        \item Append $\mathbf{9}$ to obtain
                        \[
                            P
                            =
                            \young(145,27,38,6) \raisebox{20pt}{$~\leftarrow\mathbf{9}$}
                            =
                            \young(1459,27,38,6).
                        \]
                    \end{itemize}
                    }
                    
                    \noindent It follows that $w_{1} = 1459$, $w_{2} =
                    27$, $w_{3} = 38$, and $w_{4} = 6$.  As in
                    Example~\ref{eg:SingleRowBumpingExample}, the
                    length of the longest increasing subsequence of
                    $\sigma$ is given by $\ell_{9}(\sigma) = 4 =
                    |w_{1}|$.  This is the exact result due to
                    Schensted \cite{refSchensted1961} in 1961, with a
                    sweeping generalization involving more than just
                    $w_{1}$ also given by Greene \cite{refGreene1974}
                    in 1974.
                    
                \end{Example}
                
                The standard Young tableau $P(\sigma) = (w_{1}, w_{2},
                \ldots, w_{r})$ contains a wealth of combinatorial
                information about the permutation $\sigma \in
                \mathfrak{S}_{n}$.  To get a sense of this, we first
                define the \emph{shape} $\lambda$ of $P$ to be the
                sequence of lengths $\lambda = (|w_{1}|, |w_{2}|,
                \cdots, |w_{r}|)$.  One can prove (see
                \cite{refFulton1997}) that the partial permutations
                constituting $P$ much satisfy the conditions
                \[
                    |w_{1}| \geq |w_{2}| \geq \cdots \geq |w_{r}|
                    \mbox{\quad and \quad}
                    |w_{1}| + |w_{2}| + \cdots + |w_{r}| = n
                \]
                so that $\lambda$ is an example of a so-called
                \emph{partition} of $n$.  (This is denoted $\lambda
                \vdash n$.)  It follows that $P$ can always be
                represented by a top- and left-justified shape as in
                Example~\ref{eg:SchenstedInsertionExample} above.  One
                can furthermore prove (again, see
                \cite{refFulton1997}) that the entries in $P$ must
                increase when read down each column and when read
                (from left to right) across each row.
                
                In general, one defines a standard Young tableau to be
                any ``filling'' of a partition shape $\lambda \vdash
                n$ using each element from the set $[n]$ exactly once
                and such that both this row and column condition are
                satisfied.
                
                When presented with a combinatorial algorithm like
                Schensted Insertion, it is natural to explore the
                invertibility of the steps involved.  Such
                considerations lead not only to a better understanding
                of the combinatorial objects output from the algorithm
                but also to potentially useful inverse constructions.
                With this in mind, the following is a natural first
                question to ask:
                
                \begin{Question}
                \label{huh:RSKQuestion1}
                    Given a standard Young tableau $P$, is there some
                    ``canonical'' permutation $\sigma$ such that
                    Schensted Insertion applied to $\sigma$ yields $P
                    = P(\sigma)$?
                \end{Question}
                
                To answer this questions, we first exploit the
                column-filling condition for standard Young tableaux.
                In particular, since elements must increase when read
                down each column, it follows that Schensted Insertion
                applied to a decreasing sequence will yield a tableau
                having exactly one column:
                
                {\singlespacing
                \begin{displaymath}
                    P(d_{1} > d_{2} > \cdots > d_{k}) =
                    \begin{array}{l}
                        \young(\done,\dtwo)  \\
                        \ \vdots             \\
                        \young(\dk) 
                    \end{array}.
                \end{displaymath}
                }
                
                \noindent Now, suppose that we have two decreasing
                sequences
                \[
                    d_{1} > d_{2} > \cdots > d_{k}
                    \mbox{\quad and \quad}
                    e_{1} > e_{2} > \cdots > e_{\ell}
                \]
                with $k \geq \ell$ and each $d_{i} < e_{i}$ for $i =
                1, \ldots, \ell$.  Then it is easy to see that
                
                {\singlespacing
                \begin{equation}
                \label{eqn:TableauForTwoDecreasingSequences}
                    P(d_{1} d_{2} \cdots d_{k} e_{1} e_{2} \cdots 
                    e_{\ell}) =
                    \begin{array}{l}
                        \young(\done\eone,\dtwo\etwo)  \\
                        \ \vdots \ \ \ \vdots          \\
                        \young(\dell\eell)             \\
                        \ \vdots                       \\
                        \young(\dk) 
                    \end{array}.
                \end{equation}
                }
                
                \noindent This motivates the definition of the
                \emph{column word} $\mbox{col}(P)$ of a standard Young
                tableau $P$, which is formed by reading up each of the
                columns of $P$ from left to right.
                
                \begin{Example}
                \label{eg:TableauColumnWordExample}
                    We have that
                    \[
                        \mbox{col}\left( \ \young(1459,27,38,6) \ \right)
                        =
                        632187459.
                    \]
                    
                    \noindent One can also check that
                    \[
                        P\left(\mbox{col}\left( \ \young(1459,27,38,6) \ \right)\right)
                        =
                        P(632187459)
                        =
                        \young(1459,27,38,6).
                    \]
                    This illustrates the following important property
                    of column words (which is proven by induction and
                    repeated use of
                    Equation~\eqref{eqn:TableauForTwoDecreasingSequences}):
                    
                \end{Example}
                
                \begin{Lemma}
                \label{lem:ColumnWordIdentity}
                    
                    Given a permutation $\sigma \in \mathfrak{S}_{n}$,
                    $P(\mathrm{col}(P(\sigma))) = P(\sigma)$.
                
                \end{Lemma}
                
                For a more algebraic view of this result, denote by
                $\mathfrak{T}_{n}$ the set of all standard Young
                tableaux with some partition shape $\lambda \vdash n$.
                Schensted Insertion and the above column word
                operation can then be viewed as maps $P:
                \mathfrak{S}_{n} \to \mathfrak{T}_{n}$ and
                $\mbox{col}: \mathfrak{T}_{n} \to \mathfrak{S}_{n}$,
                respectively.  With this notation,
                Lemma~\ref{lem:ColumnWordIdentity} becomes
                
                \begin{ReuseTheorem}{Lemma}{lem:ColumnWordIdentity}
                    The composition $P\circ\mathrm{col}$ is the
                    identity map on the set $\mathfrak{T}_{n}$.
                \end{ReuseTheorem}
                
                \noindent In particular, even though
                $\mbox{col}(P(\mathfrak{S}_{n})) =
                \mbox{col}(\mathfrak{T}_{n})$ is a proper subset of
                the symmetric group $\mathfrak{S}_{n}$, we nonetheless
                have that $P(\mbox{col}(P(\mathfrak{S}_{n}))) =
                P(\mathfrak{S}_{n}) = \mathfrak{T}_{n}$.  As such, it
                makes sense to define the following non-trivial
                equivalence relation on $\mathfrak{S}_{n}$, with each
                element of $\mbox{col}(\mathfrak{T}_{n})$ being the
                most natural choice of representative for the distinct
                equivalence class to which it belongs.
                
                \begin{Definition}
                \label{defn:KnuthEquivalence}
                    Two permutation $\sigma, \tau \in
                    \mathfrak{S}_{n}$ are \emph{Knuth equivalent},
                    written $\sigma \stackrel{K}{\sim} \tau$, if they
                    yield the same standard Young tableau $P(\sigma) =
                    P(\tau)$ under Schensted Insertion
                    (Algorithm~\ref{alg:SchenstedInsertion}).
                \end{Definition}
                
                \begin{Example}
                \label{eg:KnuthEquivalenceExample}
                    One can check (using the so-called Hook Length
                    Formula; see \cite{refFulton1997}) that there are
                    216 permutations $\sigma \in \mathfrak{S}_{9}$
                    satisfying $\sigma \stackrel{K}{\sim} 632187459$.
                    E.g., as illustrated in
                    Example~\ref{eg:TableauColumnWordExample},
                    $364827159 \stackrel{K}{\sim} 632187459$.
                \end{Example}
                
                In order to characterize the equivalence classes
                formed under $\stackrel{K}{\sim}$, it turns out (see
                \cite{refFulton1997}) that the following two
                examples are generic enough to characterize Knuth
                equivalence up to transitivity and order-isomorphism:
                \[
                    P(213) =  P(231) = \young(13,2) \mbox{\quad and 
                    \quad } P(312) = P(132) = \young(12,3).
                \]
                \noindent In particular, these motivate the so-called
                \emph{Knuth relations}:
                
                \begin{Definition}
                \label{defn:KnuthRelations}
                    Given two permutations $\sigma, \tau \in
                    \mathfrak{S}_{n}$, we define the \emph{Knuth
                    relations} $\stackrel{K1}{\sim}$ and
                    $\stackrel{K2}{\sim}$ on $\mathfrak{S}_{n}$ as
                    
                    \begin{itemize}
                        \item[(K1)] $\sigma \stackrel{K1}{\sim} \tau$
                        if $\sigma$ can be obtained from $\tau$ either
                        by
                            \begin{itemize}
                                \item[(K1-1)] changing a substring
                                order-isomorphic to $213$ in $\sigma$
                                into a substring order-isomorphic to
                                $231$ in $\tau$
                                
                                \item[(K1-2)] or by changing a
                                substring order-isomorphic to $231$ in
                                $\sigma$ into a substring
                                order-isomorphic to $213$ in $\tau$.
                            \end{itemize}
                            
                        \item[(K2)] $\sigma \stackrel{K2}{\sim} \tau$
                        if $\sigma$ can be obtained from $\tau$ either
                        by
                            \begin{itemize}
                                \item[(K2-1)] changing a substring
                                order-isomorphic to $312$ in $\sigma$
                                into a substring order-isomorphic to
                                $132$ in $\tau$
                                
                                \item[(K2-2)] or by changing a
                                substring order-isomorphic to $132$ in
                                $\sigma$ into a substring
                                order-isomorphic to $312$ in $\tau$.
                            \end{itemize}
                    \end{itemize}
                
                \end{Definition}
                
                \noindent One can show (again, see
                \cite{refFulton1997}) that Knuth equivalence
                $\stackrel{K}{\sim}$ on the symmetric group
                $\mathfrak{S}_{n}$ is the equivalence relation
                generated by $\stackrel{K1}{\sim}$ and
                $\stackrel{K2}{\sim}$.
                
                \begin{Example}
                \label{eg:KnuthRelationsExample}
                    From Example~\ref{eg:KnuthEquivalenceExample}, we
                    see that $364827159 \stackrel{K}{\sim} 632187459$.
                    This equivalence under Schensted Insertion can be
                    obtained by the following sequence of Knuth
                    relations, where we have underlined the
                    appropriate order-isomorphic substring being
                    permuted in order to move from one step to the
                    next.  (E.g., the substring $364$ in
                    $\underline{364}827159$ is order-isomorphic to
                    $132$ because $3$, $6$, and $4$ occur in the same
                    order with the same relative magnitudes as do $1$,
                    $3$, and $2$ in $132$.)
                    \begin{displaymath}
                        \begin{array}{ccccccccc}
                            \underline{364}827159 & 
                            \stackrel{K2}{\sim}   &
                            6\underline{342}87159 &
                            \stackrel{K1}{\sim}   &
                            63248\underline{715}9 &
                            \stackrel{K2}{\sim}   &
                            632\underline{481}759 &
                            \stackrel{K1}{\sim}   &
                            63\underline{241}8759
                            \\
                                                  &
                            \stackrel{K1}{\sim}   & 
                            6321\underline{487}59 &
                            \stackrel{K2}{\sim}   & 
                            63218\underline{475}9 &
                            \stackrel{K2}{\sim}   &
                            632187459.            &
                                                  &
                                                  
                        \end{array}
                    \end{displaymath}
                \end{Example}
                
               \noindent Given how many intermediate permutations are
               needed in order to realizing the Knuth equivalence of
               $364827159$ and $632187459$ via the two Knuth
               relations, one might expect Knuth equivalence classes
               to be fairly large in general.  This suggests the next
               question, which is a natural follow-up to
               Question~\ref{huh:RSKQuestion1}.

                \begin{Question}
                \label{huh:RSKQuestion2}
                    Given a standard Young tableau $P$, how much extra
                    ``bookkeeping'' is necessary in order to uniquely
                    specify a permutation $\sigma$ such that $P =
                    P(\sigma)$?
                \end{Question}
                
                One can show (see \cite{refFulton1997} for a detailed
                account) that, due to the row- and column-filling
                conditions on $P$, it suffices to keep track of the
                order in which new boxes are added to the shape of
                $P$.  In particular, each entry in $P$ has a unique
                ``bumping path'' by which it reached its final
                position, and so this path can be inverted in order to
                ``unbump'' (a.k.a.~\emph{reverse row bump}) the entry.
                This motivates a bijective extension of Schensted
                Insertion called the \emph{RSK Correspondence}, which
                originated from the study of representations of the
                symmetric group by Robinson \cite{refRobinson1938} in
                1938.  The modern ``bumping'' version that we present,
                though, resulted from work by Schensted
                \cite{refSchensted1961} in 1961 that extended
                Robinson's algorithm to combinatorial words.  Further
                generalization was also done by Knuth
                \cite{refKnuth1970} for so-called
                $\mathbb{N}$-matrices (or, equivalently, so-called
                two-line arrays) in 1970.  (See Fulton
                \cite{refFulton1997} for the appropriate definitions
                and for a detailed account of the differences between
                these algorithms; Knuth also includes a description in
                \cite{refKnuth1998}.)
                
                \begin{Algorithm}[RSK Correspondence]
                \label{alg:RSKAlgorithm}
                {\tt
                    ~\\
                    Input:~a permutation $\sigma =
                    \sigma_{1}\cdots\sigma_{n} \in \mathfrak{S}_{n}$\\
                    Output:~the pair of standard Young tableaux $P =
                    P(\sigma)$ and $Q = Q(\sigma)$
                    
                    \begin{itemize}
                        
                        \item Use Schensted Insertion
                        (Algorithm~\ref{alg:SchenstedInsertion}) to
                        build $P$.
                        
                        \item For each $i = 1, \ldots, n$, when the
                        $i^{\rm th}$ box is added to the shape of 
                        $P$,\\
                        add the box $\young(i)$ to $Q$ so that $P$ and
                        $Q$ maintain the same shape.
                        
                    \end{itemize}
                }
                \end{Algorithm}
                
                \noindent We call $P(\sigma)$ the \emph{insertion
                tableau} and $Q(\sigma)$ the \emph{recording tableau}
                corresponding to $\sigma \in \mathfrak{S}_{n}$ under
                the RSK Correspondence.  Since the construction of the
                recording tableau suffices to invert Schensted
                Insertion, this correspondence yields a bijection
                between the symmetric group $\mathfrak{S}_{n}$ and the
                set of all ordered pairs of standard Young tableaux
                such that the two tableaux have the same partition
                shape.  We denote this bijection at the element level
                by $\sigma \stackrel{RSK}{\longleftrightarrow}
                (P(\sigma), Q(\sigma))$.
                
                \begin{Example}
                \label{eg:RSKExample}
                    Given the permutation $\sigma = 364827159 \in
                    \mathfrak{S}_{9}$, we form the insertion tableau
                    $P = P(\sigma)$
                    (cf.~Example~\ref{eg:SchenstedInsertionExample}
                    above) and the recording tableau $Q = Q(\sigma)$
                    under the RSK Correspondence
                    (Algorithm~\ref{alg:RSKAlgorithm}) as follows:
                    
                    {\singlespacing
                    \begin{itemize}
                        \item Start with $P = \emptyset$ and insert
                        $\mathbf{3}$ into $P$ to obtain $P =
                        \young(3)$ and $Q = \young(1)$.
                    
                        \item Append $\mathbf{6}$ to $P$ so that
                        \[
                            P = \young(36)
                            \mbox{\quad and\quad}
                            Q = \young(12).
                        \]
                    
                        \item Use $\mathbf{4}$ to bump $6$ in $P$ so 
                        that
                        \[
                            P = \young(34,6)
                            \mbox{\quad and\quad}
                            Q = \young(12,3).
                        \]
                    
                        \item Append $\mathbf{8}$ to the top row of $P$ so that
                        \[
                            P = \young(348,6)
                            \mbox{\quad and\quad}
                            Q = \young(124,3).
                        \]
                        
                        \item Use $\mathbf{2}$ to bump $3$ in $P$ so 
                        that
                        \[
                            P = \young(248,3,6)
                            \mbox{\quad and\quad}
                            Q = \young(124,3,5).
                        \]
                        
                        \item Use $\mathbf{7}$ to bump $8$ in $P$ so 
                        that
                        \[
                            P = \young(247,38,6)
                            \mbox{\quad and\quad}
                            Q = \young(124,36,5).
                        \]
                    
                        \item Use $\mathbf{1}$ to bump $2$ in $P$ so 
                        that
                        \[
                            P = \young(147,28,3,6)
                            \mbox{\quad and\quad}
                            Q = \young(124,36,5,7).
                        \]
                    
                        \item Use $\mathbf{5}$ to bump $7$ in $P$ so 
                        that
                        \[
                            P = \young(145,27,38,6)
                            \mbox{\quad and\quad}
                            Q = \young(124,36,58,7).
                        \]
                        
                        \item Append $\mathbf{9}$ to the top row of $P$ so that
                        \[
                            P = \young(1459,27,38,6)
                            \mbox{\quad and\quad}
                            Q = \young(1249,36,58,7).
                        \]
                    \end{itemize}
                    }
                    
                    \noindent It follows that $364827159
                    \stackrel{RSK}{\longleftrightarrow} \left( \
                    \young(1459,27,38,6)\raisebox{-22pt}{,} \quad
                    \young(1249,36,58,7) \ \right)$.\bigskip
                    
                    One can also check that the inverse of
                    $\sigma$ (when thought of a function) satisfies
                    \[  
                        \sigma^{-1}
                        =
                        751382649
                        \stackrel{RSK}{\longleftrightarrow}
                        \left( \
                        \young(1249,36,58,7)\raisebox{-22pt}{,}
                        \quad \young(1459,27,38,6) \
                        \right).
                    \]
                    This illustrates the following amazing fact about
                    the RSK Correspondence:
                    
                \end{Example}
                
                \begin{Theorem}[Sch\"{u}tzenberger Symmetry for the RSK Correspondence]
                \label{thm:RSKSchuetzenbergerSymmetryProperty}
                    Given a permutation $\sigma \in \mathfrak{S}_{n}$,
                    $\sigma \stackrel{RSK}{\longleftrightarrow}
                    (P(\sigma), Q(\sigma))$ if and only if
                    $\sigma^{-1} \stackrel{RSK}{\longleftrightarrow}
                    (Q(\sigma), P(\sigma))$.
                \end{Theorem}
                
                Sch\"{u}tzenberger Symmetry, which was first proven by
                a direct combinatorial argument in
                \cite{refSchutzenberger1963}, is only one of many
                remarkable properties of the RSK Correspondence.  A
                particularly good account of the numerous consequences
                of the RSK Correspondence in such fields as
                Representation Theory can be found in Sagan
                \cite{refSagan2000}.  One can also read about the RSK
                Correspondence on words and on $\mathbb{N}$-matrices
                in Fulton \cite{refFulton1997}.
                
                Another remarkable fact about the RSK Correspondence
                is that it can be realized without explicit
                involvement of Single Row Bumping
                (Algorithm~\ref{alg:SingleRowBumping}).  One of the
                more striking alternatives involves the so-called
                \emph{shadow diagram} of a permutation, as introduced
                by Viennot \cite{refViennot1977} in the context of
                further explicating Sch\"{u}tzenberger Symmetry.  As
                we review in Section~\ref{sec:ExtendingPS:GeometricPS}
                below,
                Theorem~\ref{thm:RSKSchuetzenbergerSymmetryProperty}
                follows trivially from Viennot's use of shadow
                diagrams.
                
                At the same time, it is also interesting to look at
                when the full RSK Correspondence is essentially
                unnecessary.  This motivates the following question:
                
                \begin{Question}
                \label{huh:RSKQuestion3}
                    Given a standard Young tableau $P$, under what
                    conditions is no extra ``bookkeeping'' necessary
                    in order to uniquely recover a permutation
                    $\sigma$ such that $P = P(\sigma)$?
                \end{Question}
                
                As a direct consequence of the Sch\"{u}tzenberger
                Symmetry property for the RSK Correspondence, there is
                a bijection between the set $\mathfrak{I}_{n}$ of
                involutions in $\mathfrak{S}_{n}$ and the set
                $\mathfrak{T}_{n}$ of all standard Young tableaux.
                (An \emph{involution} is any permutation that is equal
                to its own inverse.)  However, this bijection doesn't
                yield any information about the size of the
                equivalence classes containing each involution.  Thus,
                in order to answer Question~\ref{huh:RSKQuestion3},
                one must instead consider permutations that ``avoid''
                the order-isomorphic substrings (a.k.a.~\emph{block
                patterns} as defined in
                Section~\ref{sec:Intro:Motivation:Patterns} below)
                used to define the Knuth relations.  A Knuth
                equivalence class cannot be a singleton set unless the
                only permutation it contains cannot be transformed
                into a Knuth equivalent permutation via a Knuth
                relation.  Given how restrictive a condition this is,
                a tableau $P$ satisfies
                Question~\ref{huh:RSKQuestion3} if and only if $P$ is
                a single row or column.
                
                Chapters~\ref{sec:PSasAlgorithm} and
                \ref{sec:ExtendingPS} of this Dissertation
                are focused upon addressing
                Questions~\ref{huh:RSKQuestion1},
                \ref{huh:RSKQuestion2}, and \ref{huh:RSKQuestion3}
                when adapted to the combinatorial algorithm Patience
                Sorting.  As we will see, the notion of ``pattern
                avoidance'', which is made explicit in
                Section~\ref{sec:Intro:Motivation:Patterns}, provides
                the unifying language for characterizing our
                responses.
            
            \subsection[Patience Sorting as ``Non-recursive'' Schensted Insertion]{Patience Sorting as\\ ``Non-recursive'' Schensted Insertion}
            \label{sec:Intro:Motivation:PS}
                The term \emph{Patience Sorting} was introduced in
                1962 by Mallows \cite{refMallows1962, refMallows1963}
                as the name of a two-part card sorting algorithm
                invented by A.~S.~C.~Ross.  The first part of this
                algorithm, which Mallows referred to as a ``patience
                sorting procedure'', involves partitioning a shuffled
                deck of cards into a collection of sorted subsequences
                called \emph{piles}.  Unless otherwise noted, we take
                our ``(shuffled) deck of cards'' to be a permutation.
                
                \begin{Algorithm}[Mallows' Patience Sorting Procedure]
                \label{alg:MallowsPSprocedure}
                {\tt
                    ~\\
                    Input:~a shuffled deck of cards $\sigma = c_{1}
                    c_{2} \cdots c_{n} \in \mathfrak{S}_{n}$\\
                    Output:~the sequence of partial permutations
                    $r_{1}, r_{2}, \ldots, r_{m}$

                    \begin{enumerate}
                        \item First form a new pile $r_{1}$ using the
                        card $c_{1}$.

                        \item Then, for each $i = 2, \ldots, n$,
                        consider the cards $d_{1}, d_{2}, \ldots,
                        d_{k}$ atop the piles $r_{1}, r_{2}, \ldots,
                        r_{k}$ that have already been formed.

                            \begin{enumerate}
                                \item If ${\displaystyle c_{i} >
                                \max_{1 \, \leq \, j \, \leq \,
                                k}\{d_{j}\}}$, then form a new pile
                                $r_{k + 1}$ using $c_{i}$.

                                \item Otherwise, find the left-most
                                card $d_{j}$ that is larger than
                                $c_{i}$\\ and place the card $c_{i}$
                                atop pile $r_{j}$.
                                In other words,
                                set
                                \[
                                    d_{j} = c_{i} \
                                    \ \mathtt{where} \ \ j = \min_{1 
                                    \, \leq \,
                                    m \, \leq \, k} \{ m \ | \ c_{i} <
                                    d_{m} \}.
                                \]
                            \end{enumerate}
                    \end{enumerate}
                }
                \end{Algorithm}

                \noindent We call the collection of subsequences $R =
                R(\sigma) = \{ r_{1}, r_{2}, \ldots, r_{m} \}$ the
                \emph{pile configuration} associated to the deck of
                cards $\sigma \in \mathfrak{S}_{n}$ and illustrate
                their formation in Example~\ref{eg:NormalPSexample}
                below.  In keeping with the language of ``piles'', we
                will often write the constituent subsequences $r_{1},
                r_{2}, \ldots, r_{m}$ vertically, in order, and
                bottom-justified with respect to the largest value in
                the pile.  This is illustrated in the following
                example.
                
                \begin{Example}
                \label{eg:PileConfigurationExample}
                    The pile configuration $R(64518723) = \{641, 52,
                    873\}$ (formed in Example~\ref{eg:NormalPSexample}
                    below) has piles $r_{1} = 641$, $r_{2} = 52$, and
                    $r_{3} = 873$.  We represent this visually as

                    {\singlespacing
                    \begin{center}
                        $R(64518723) \, = $
                        \begin{tabular}{l l l}
                            1 &   & 3 \\
                            4 & 2 & 7 \\
                            6 & 5 & 8 \end{tabular}.\\
                    \end{center}
                    }
                    
                \end{Example}
                
                \bigskip
                
                \begin{Example}
                \label{eg:NormalPSexample}
                    Given the deck of cards $\sigma = 64518723 \in
                    \mathfrak{S}_{8}$, we form $R(\sigma)$ under
                    Patience Sorting (Algorithm
                    \ref{alg:MallowsPSprocedure}) as follows:
                    
                    {\singlespacing
                    \begin{tabular}{l p{80pt} l r}
                        \begin{minipage}[c]{110pt}
                            First, use \textbf{6} to form a new pile:
                        \end{minipage} 
                        &
                        \begin{tabular}{l l l}
                                        &    &  \\
                                        &    &  \\
                             \textbf{6} &    &  
                        \end{tabular}
                        &
                        \begin{minipage}[c]{110pt}
                            Then place \textbf{4} atop this new pile:
                        \end{minipage}
                        &
                        \begin{tabular}{l l l}
                                        &    &  \\
                             \textbf{4} &    &  \\
                                     6  &    &  
                        \end{tabular}
                    \end{tabular}\\ \\

                    \begin{tabular}{l p{80pt} l r}
                        \begin{minipage}[c]{110pt}
                            Use \textbf{5} to form a new pile:
                        \end{minipage}
                        &
                        \begin{tabular}{l l l}
                               &            &  \\
                             4 &            &  \\
                             6 & \textbf{5} &  
                        \end{tabular}
                        &
                        \begin{minipage}[c]{110pt}
                            Then place \textbf{1} atop the left-most pile: 
                        \end{minipage}
                        &
                        \begin{tabular}{l l l}
                            \textbf{1} &    &  \\
                                    4  &    &  \\
                                    6  &  5 &  
                        \end{tabular}
                    \end{tabular}\\ \\
                        
                    \begin{tabular}{l p{80pt} l r}
                        \begin{minipage}[c]{110pt}
                            Use \textbf{8} to form a new pile: 
                        \end{minipage}
                        &
                        \begin{tabular}{l l l}
                            1 &   &            \\
                            4 &   &            \\
                            6 & 5 & \textbf{8}
                        \end{tabular}
                        &
                        \begin{minipage}[c]{110pt}
                            Then place \textbf{7} atop this new pile: 
                        \end{minipage}
                        &
                        \begin{tabular}{l l l}
                            1 &   &            \\
                            4 &   & \textbf{7} \\
                            6 & 5 &         8
                        \end{tabular}
                    \end{tabular}\\ \\

                    \begin{tabular}{l p{80pt} l r}
                        \begin{minipage}[c]{110pt}
                            Place \textbf{2} atop the middle pile: 
                        \end{minipage}
                        &
                        \begin{tabular}{l l l}
                            1 &            &   \\
                            4 & \textbf{2} & 7 \\
                            6 &         5  & 8
                        \end{tabular}
                        &
                        \begin{minipage}[c]{110pt}
                            Finally, place \textbf{3} atop the right-most pile:
                        \end{minipage}
                        &
                        \begin{tabular}{l l l}
                            1 &   & \textbf{3} \\
                            4 & 2 &         7  \\
                            6 & 5 &         8
                        \end{tabular}
                    \end{tabular}\\ \\
                    }
                    
                    \noindent Now, in order to affect sorting, cards
                    can be removed one at a time from these piles in
                    the order $1, 2, \ldots, 8$.  Note that, by
                    construction, the cards in each pile decrease when
                    read from bottom to top.  Consequently, the
                    appropriate card will always be atop a pile at
                    each step of this removal process:
                    
                    {\singlespacing
                    \begin{tabular}{l p{70pt} l r}
                        \begin{minipage}[c]{120pt}
                            After removing $\mathbf{1}$: 
                        \end{minipage}
                        &
                        \begin{tabular}{l l l}
                              &   & 3 \\
                            4 & 2 & 7 \\
                            6 & 5 & 8
                        \end{tabular}
                        &
                        \begin{minipage}[c]{120pt}
                            After removing $\mathbf{2}$:
                        \end{minipage}
                        &
                        \begin{tabular}{l l l}
                              &   & 3 \\
                            4 &   & 7 \\
                            6 & 5 & 8
                        \end{tabular}
                    \end{tabular}\\ \\

                    \begin{tabular}{l p{70pt} l r}
                        \begin{minipage}[c]{120pt}
                            After removing $\mathbf{3}$: 
                        \end{minipage}
                        &
                        \begin{tabular}{l l l}
                              &   &   \\
                            4 &   & 7 \\
                            6 & 5 & 8
                        \end{tabular}
                        &
                        \begin{minipage}[c]{120pt}
                            After removing $\mathbf{4}$:
                        \end{minipage}
                        &
                        \begin{tabular}{l l l}
                              &   &   \\
                              &   & 7 \\
                            6 & 5 & 8
                        \end{tabular}
                    \end{tabular}\\ \\
                    
                    \begin{tabular}{l p{70pt} l r}
                        \begin{minipage}[c]{120pt}
                            After removing $\mathbf{5}$: 
                        \end{minipage}
                        &
                        \begin{tabular}{l l l}
                              &                &   \\
                              &                & 7 \\
                            6 & $\phantom{5}$  & 8
                        \end{tabular}
                        &
                        \begin{minipage}[c]{120pt}
                            And so on.
                        \end{minipage}
                    \end{tabular}\\
                    }

                \end{Example}

                When applying Algorithm~\ref{alg:MallowsPSprocedure}
                to a permutation $\sigma = c_{1} c_{2} \cdots c_{n}
                \in \mathfrak{S}_{n}$, each card $c_{i}$ is either
                larger than the top card of every pile or is placed
                atop the left-most top card $d_{j}$ larger than it.
                As such, the cards $d_{1}, d_{2}, \ldots, d_{k}$ atop
                the piles will always be in increasing order (from
                left to right) at each step of the algorithm, and it
                is in this sense that
                Algorithm~\ref{alg:MallowsPSprocedure} resembles
                Schensted Insertion
                (Algorithm~\ref{alg:SchenstedInsertion}).  The
                distinction is that cards remain in place and have
                other cards placed on top of them when ``bumped''
                rather than endure insertion into successive rows of a
                standard Young tableau through recursive application
                of Single Row Bumping
                (Algorithm~\ref{alg:SingleRowBumping}).
                
                Note, in particular, that the cards $d_{1}, d_{2},
                \ldots, d_{k}$ correspond exactly to the top row of
                the standard Young tableau formed at each stage of
                Schensted Insertion.  Consequently, the number of
                piles formed under Patience Sorting is exactly equal
                to the length of the longest increasing subsequence in
                the deck of cards.
                
                Given the algorithmic simplicity of Patience Sorting,
                the pile configuration $R(\sigma)$ is commonly offered
                as a tangible realization for the length
                $\ell_{n}(\sigma)$ of the longest increasing
                subsequence in a permutation $\sigma \in
                \mathfrak{S}_{n}$.  This is particularly true in the
                context of Probabilistic Combinatorics since the
                asymptotic number of piles formed must follow the
                highly celebrated Baik-Deift-Johansson Theorem from
                \cite{refBDJ1999}.  In other words, even though it is
                easy to construct $R(\sigma)$, there is no simple way
                to completely describe the number of piles formed
                $\ell_{n}(\sigma)$ without borrowing from Random
                Matrix Theory.  (A good historical survey of attempts
                previous to \cite{refBDJ1999} can be found in both
                \cite{refAD1999} and \cite{refStanley2006}.)
                
                According to the Baik-Deift-Johansson Theorem, the
                distribution for the number of piles formed under
                Patience Sorting converges asymptotically to the
                Tracy-Widom $F_{2}$ distribution (up to an appropriate
                rescaling).  Remarkably, though, $F_{2}$ originated in
                work by Tracy and Widom \cite{refTW1994} as the
                asymptotic distribution for the largest eigenvalue of
                a random Hermitian matrix (again, up to rescaling).
                Because of this deep connection between Patience
                Sorting and Probabilistic Combinatorics, it has been
                suggested (see, e.g., \cite{refKuykendall1999},
                \cite{refMackenzie1998} and \cite{refPeterson1999};
                cf.~\cite{refConrey2003}) that studying
                generalizations of Patience Sorting might be the key
                to tackling certain open problems that can be viewed
                from the standpoint of Random Matrix Theory --- the
                most notable being the Riemann Hypothesis.
                
                Another natural direction of study involves
                characterizing the objects output from both
                Algorithm~\ref{alg:MallowsPSprocedure} and an
                appropriate bijective extension.
                Chapter~\ref{sec:PSasAlgorithm} is largely devoted to
                the combinatorics that arises from various
                characterizations for pile configurations.  Then, in
                Chapter~\ref{sec:ExtendingPS}, we study the
                combinatorics that arises from a full, non-recursive
                analog of the RSK Correspondence.  In particular, we
                mimic the RSK recording tableau construction so that
                ``recording piles'' $S(\sigma)$ are assembled along
                with the usual pile configuration $R(\sigma)$ under
                Patience Sorting (which by analogy to RSK we will
                similarly now call ``insertion piles'').  We refer to
                the resulting (ordered) pair of pile configurations as
                a \emph{stable pair} and denote $\sigma
                \stackrel{XPS}{\longleftrightarrow} (R(\sigma),
                S(\sigma))$.

                \begin{Algorithm}[Extended Patience Sorting]
                \label{alg:ExtendedPSalgorithm}
                {\tt
                    ~\\
                    Input:~a shuffled deck of cards $\sigma = c_{1}
                    c_{2} \cdots c_{n} \in \mathfrak{S}_{n}$\\
                    Output:~the ordered pair $(R(\sigma), S(\sigma))$
                    of pile configurations, where we\\ denote $R(\sigma)
                    = \{r_{1}, r_{2}, \ldots, r_{m}\}$ and $S(\sigma)
                    = \{s_{1}, s_{2}, \ldots, s_{m}\}$

                    \begin{enumerate}
                        \item First form a new pile $r_{1}$ using the
                        card $c_{1}$ and set $s_{1} = 1$.

                        \item Then, for each $i = 2, \ldots, n$,
                        consider the cards $d_{1}, d_{2}, \ldots,
                        d_{k}$ atop the piles $r_{1}, r_{2}, \ldots,
                        r_{k}$ that have already been formed.

                            \begin{enumerate}
                                \item If ${\displaystyle c_{i} >
                                \max_{1 \, \leq \, j \, \leq \,
                                k}\{d_{j}\}}$, then form a new pile
                                $r_{k + 1}$ using $c_{i}$ and set\\
                                $s_{k + 1} = i$.

                                \item Otherwise, find the left-most
                                card $d_{j}$ that is larger than
                                $c_{i}$\\ and place the card $c_{i}$
                                atop pile $r_{j}$  while 
                                simultaneously\\ placing $i$ at the 
                                bottom of pile $s_{j}$.  In other words, set
                                \[
                                    d_{j} = c_{i} \
                                    \ \mathtt{where} \ \ j = \min_{1 
                                    \, \leq \,
                                    m \, \leq \, k} \{ m \ | \ c_{i} <
                                    d_{m} \}.
                                \]
                                and insert $i$ at the bottom of pile
                                $s_{j}$.
                            \end{enumerate}
                    \end{enumerate}
                }
                \end{Algorithm}

                By construction, the pile configurations in the
                resulting stable pair must have the same notion of
                ``shape'', which we define as follows.

                \begin{Definition}
                    Given a pile configuration $R=\{r_{1}, r_{2},
                    \ldots, r_{m}\}$ that has been formed from $n$
                    cards, we call the $m$-tuple $\mathrm{sh}(R)$ the
                    \emph{shape} of $R$, where
                    \[
                        \mathrm{sh}(R)
                        =
                        (|r_{1}|, |r_{2}|, \ldots, |r_{m}|).
                    \]
                \end{Definition}
                
                \noindent Note, in particular, that $sh(R)$ satisfies
                \[
                    |r_{1}|, |r_{2}|, \ldots, |r_{m}| \in [n]
                    \mbox{ \ and \ }
                    |r_{1}| + |r_{2}| + \cdots + |r_{m}| = n
                \]
                and thus is an example of a so-called
                \emph{composition} of $n$.  By partial analogy to the
                notation for a partition of $n$ (and since partitions
                are a special case of compositions), this is denoted
                by $\mathrm{sh}(R) \models n$.
                
                Having established this shape convention, we now
                illustrate Extended Patience Sorting
                (Algorithm~\ref{alg:ExtendedPSalgorithm}) in the
                following example.

                \begin{Example}
                \label{eg:ExtendedPSexample}
                    Given $\sigma = 64518723 \in \mathfrak{S}_{8}$, we
                    form the following stable pair $(R(\sigma),
                    S(\sigma))$ under
                    Algorithm~\ref{alg:ExtendedPSalgorithm}, with
                    shape $\mathrm{sh}(R(\sigma)) =
                    \mathrm{sh}(S(\sigma)) = (3, 2, 3) \models 8$.
                    
                    {\singlespacing
                    \begin{tabular}{l p{56pt} p{72pt} l p{56pt} l}
                        \begin{minipage}[c]{48pt}$\phantom{foo}$\end{minipage}
                        &
                        \begin{minipage}[c]{48pt}insertion piles\end{minipage}
                        &
                        \begin{minipage}[c]{48pt}recording piles\end{minipage}
                        &
                        \begin{minipage}[c]{48pt}$\phantom{foo}$\end{minipage}
                        &
                        \begin{minipage}[c]{48pt}insertion piles\end{minipage}
                        &
                        \begin{minipage}[c]{48pt}recording piles\end{minipage}
                     \end{tabular}\\

                    \begin{tabular}{l p{56pt} p{72pt} l p{56pt} l}
                        \begin{minipage}[c]{48pt}
                            After inserting \textbf{6}:
                        \end{minipage}
                        &
                        \begin{tabular}{l l l}
                                              & & \\
                                              & & \\
                               \textbf{6} & &
                        \end{tabular}
                        &
                        \begin{tabular}{l l l}
                                             & & \\
                                             & & \\
                              \textbf{1} & &
                        \end{tabular}
                        &
                        \begin{minipage}[c]{48pt}
                            After inserting \textbf{4}:
                        \end{minipage}
                        &
                        \begin{tabular}{l l l}
                                              & & \\
                               \textbf{4} & & \\
                                          6 & &
                        \end{tabular}
                        &
                        \begin{tabular}{l l l}
                                           & & \\
                                       1 & & \\
                            \textbf{2} & &
                        \end{tabular}
                    \end{tabular}\\ \\

                    \begin{tabular}{l p{56pt} p{72pt} l p{56pt} l}
                        \begin{minipage}[c]{48pt}
                            After inserting \textbf{5}:
                        \end{minipage}
                        &
                        \begin{tabular}{l l l}
                               & & \\
                            4 & & \\
                            6 & \textbf{5} &
                        \end{tabular}
                        &
                        \begin{tabular}{l l l}
                              & & \\
                           1 & & \\
                           2 & \textbf{3} &
                        \end{tabular}
                        &
                        \begin{minipage}[c]{48pt}
                            After inserting \textbf{1}:
                        \end{minipage}
                        &
                        \begin{tabular}{l l l}
                            \textbf{1} & & \\
                                       4 & & \\
                                       6 & 5 &
                        \end{tabular}
                        &
                        \begin{tabular}{l l l}
                                       1 & & \\
                                       2 & & \\
                            \textbf{4} & 3 &
                        \end{tabular}
                    \end{tabular}\\ \\

                    \begin{tabular}{l p{56pt} p{72pt} l p{56pt} l}
                        \begin{minipage}[c]{48pt}
                            After inserting \textbf{8}:
                        \end{minipage}
                        &
                        \begin{tabular}{l l l}
                            1 & & \\
                            4 & & \\
                            6 & 5 & \textbf{8}
                        \end{tabular}
                        &
                        \begin{tabular}{l l l}
                           1 & & \\
                           2 & & \\
                           4 & 3 & \textbf{5}
                        \end{tabular}
                        &
                        \begin{minipage}[c]{48pt}
                            After inserting \textbf{7}:
                        \end{minipage}
                        &
                        \begin{tabular}{l l l}
                            1 & & \\
                            4 & & \textbf{7} \\
                            6 & 5 & 8
                        \end{tabular}
                        &
                        \begin{tabular}{l l l}
                            1 & & \\
                            2 & & 5 \\
                            4 & 3 & \textbf{6}
                        \end{tabular}
                    \end{tabular}\\ \\

                    \begin{tabular}{l p{56pt} p{72pt} l p{56pt} l}
                        \begin{minipage}[c]{48pt}
                            After inserting \textbf{2}:
                        \end{minipage}
                        &
                        \begin{tabular}{l l l}
                            1 & & \\
                            4 & \textbf{2} & 7 \\
                            6 & 5 & 8
                        \end{tabular}
                        &
                        \begin{tabular}{l l l}
                           1 & & \\
                           2 & 3 & 5 \\
                           4 & \textbf{7} & 6
                        \end{tabular}
                        &
                        \begin{minipage}[c]{48pt}
                            After inserting \textbf{3}:
                        \end{minipage}
                        &
                        \begin{tabular}{l l l}
                            1 & & \textbf{3} \\
                            4 & 2 & 7 \\
                            6 & 5 & 8
                        \end{tabular}
                        &
                        \begin{tabular}{l l l}
                           1 & & 5 \\
                           2 & 3 & 6 \\
                           4 & 7 & \textbf{8}
                        \end{tabular}
                    \end{tabular}\\
                    }
                \end{Example}

                \vspace{-5mm}
                
                Given a permutation $\sigma \in \mathfrak{S}_{n}$, the
                recording piles $S(\sigma)$ indirectly label the order
                in which cards are added to the insertion piles
                $R(\sigma)$ under Patience Sorting
                (Algorithm~\ref{alg:MallowsPSprocedure}).  With this
                information, $\sigma$ can be reconstructed by removing
                cards from $R(\sigma)$ in the opposite order that they
                were added to $R(\sigma)$.  (This removal process,
                though, should not be confused with the process of
                sorting that was illustrated in
                Example~\ref{eg:NormalPSexample}.)  Consequently,
                reversing Extended Patience Sorting
                (Algorithm~\ref{alg:ExtendedPSalgorithm}) is
                significantly less involved than reversing the RSK
                Correspondence (Algorithm~\ref{alg:RSKAlgorithm})
                through recursive ``reverse row bumping''.  The
                trade-off is that entries in the pile configurations
                resulting from the former are not independent (see
                Section~\ref{sec:ExtendingPS:StablePairs}) whereas the
                tableaux generated under the RSK Correspondence are
                completely independent (up to shape).

                To see that $S(\sigma) = \{s_{1}, s_{2}, \ldots,
                s_{m}\}$ records the order in which cards are added to
                the insertion piles, instead add cards atop new piles
                $s'_{j}$ in Algorithm~\ref{alg:ExtendedPSalgorithm}
                rather than to the bottoms of the piles $s_{j}$.  This
                yields modified recording piles $S'(\sigma)$ from
                which each original recording pile $s_{j} \in
                S(\sigma)$ can be recovered by simply reflecting the
                corresponding pile $s'_{j}$ vertically.  Moreover, it
                is easy to see that the entries in $S'(\sigma)$
                directly label the order in which cards are added to
                $R(\sigma)$ under Patience Sorting.

                \begin{Example}
                \label{eg:PileReflectionExample}
                    As in Example~\ref{eg:ExtendedPSexample} above,
                    let $\sigma = 6 4 5 1 8 7 2 3 \in
                    \mathfrak{S}_{8}$.  Then $R(\sigma)$ is formed as
                    before and

                    {\singlespacing
                    \begin{center}
                        \begin{tabular}{r c c c l}
                            $S'(\sigma) \ \ =$ &
                            \begin{tabular}{l l l}
                                4 & & 8 \\
                                2 & 7 & 6 \\
                                1 & 3 & 5
                            \end{tabular}
                            & 
                            $\stackrel{\mathit{reflect}}{\leadsto}$ &
                            \begin{tabular}{l l l}
                               1 & & 5 \\
                               2 & 3 & 6 \\
                               4 & 7 & 8
                            \end{tabular}
                            & $= \ \ S(\sigma)$
                        \end{tabular}.
                    \end{center}
                    }

                \end{Example}
                
                As Example~\ref{eg:PileReflectionExample} illustrates,
                Extended Patience Sorting
                (Algorithm~\ref{alg:ExtendedPSalgorithm}) is only one
                possible way to bijectively extend Patience Sorting
                (Algorithm~\ref{alg:MallowsPSprocedure}).  What
                suggests Extended Patience Sorting as the right
                extension is a symmetric property (see
                Section~\ref{sec:ExtendingPS:Involutions}) directly
                analogous to Sch\"{u}tzenberger Symmetry for the RSK
                Correspondence
                (Theorem~\ref{thm:RSKSchuetzenbergerSymmetryProperty}):
                $\sigma \in \mathfrak{S}_{n}$ yields $(R(\sigma),
                S(\sigma))$ under Extended Patience Sorting if and
                only if $\sigma^{-1}$ yields $(S(\sigma), R(\sigma))$.
                Moreover, this symmetry result follows trivially from
                the geometric realization of Extended Patience Sorting
                given in
                Section~\ref{sec:ExtendingPS:GeometricPS:SWshadows},
                which is also a direct analogy to the geometric
                realization for the RSK Correspondence reviewed in
                Section~\ref{sec:ExtendingPS:GeometricPS:NEshadows}.
                
                While such interesting combinatorics result from
                questions suggested by the resemblance between
                Patience Sorting and Schensted Insertion, this is not
                the only possible direction of study.  In particular,
                after applying Algorithm~\ref{alg:MallowsPSprocedure}
                to a deck of cards, it is easy to recollect each card
                in ascending order from amongst the current top cards
                of the piles (and thus complete A.~S.~C.~Ross'
                original card sorting algorithm as in
                Example~\ref{eg:NormalPSexample}).  While this is not
                necessarily the fastest sorting algorithm that one
                might apply to a deck of cards, the \emph{patience} in
                \emph{Patience Sorting} is not intended to describe a
                prerequisite for its use.  Instead, it refers to how
                pile formation in
                Algorithm~\ref{alg:MallowsPSprocedure} resembles the
                placement of cards into piles when playing the popular
                single-person card game \emph{Klondike Solitaire}, and
                Klondike Solitaire is often called \emph{Patience} in
                the UK. This is more than a coincidence, though, as
                Algorithm~\ref{alg:MallowsPSprocedure} also happens to
                be an optimal strategy (in the sense of forming as few
                piles as possible; see
                Section~\ref{sec:PSasGame:Strategies}) when playing an
                idealized model of Klondike Solitaire known as
                \emph{Floyd's Game}:

                \begin{CardGame}[Floyd's Game]
                \label{game:FloydsGame}
                    Given a deck of cards $c_{1}c_{2}\cdots c_{n} \in
                    \mathfrak{S}_{n}$,

                    \begin{itemize}
                        \item place the first card $c_{1}$ from the
                        deck into a pile by itself.

                        \item Then, for each card $c_{i}$ ($i = 2,
                        \ldots, n$), either

                            \begin{itemize}

                                \item put $c_{i}$ into a new pile by
                                itself

                                \item or play $c_{i}$ on top of any
                                pile whose current top card is larger
                                than $c_{i}$.

                            \end{itemize}

                        \item The object of the game is to end with as
                        few piles as possible.

                    \end{itemize}

                \end{CardGame}

                \noindent In other words, cards are played one at a
                time according to the order that they appear in the
                deck, and piles are created in much the same way that
                they are formed under Patience Sorting.  According to
                \cite{refAD1999}, Floyd's Game was developed
                independently of Mallows' work during the 1960s as an
                idealized model for Klondike Solitaire in unpublished
                correspondence between computer scientists Bob Floyd
                and Donald Knuth.

                Note that, unlike Klondike Solitaire, there is a known
                strategy (Algorithm~\ref{alg:MallowsPSprocedure}) for
                Floyd's Game under which one will always win.  In
                fact, Klondike Solitaire --- though so popular that it
                has come pre-installed on the vast majority of
                personal computers shipped since 1989 --- is still
                poorly understood mathematically.  (Recent progress,
                however, has been made in developing an optimal
                strategy for a version called \emph{thoughtful
                solitaire} \cite{refYDRR2005}.)  As such, Persi
                Diaconis (\cite{refAD1999} and private communication
                with the author) has suggested that a deeper
                understanding of Patience Sorting and its
                generalization would undoubtedly help in developing a
                better mathematical model for analyzing Klondike
                Solitaire.
                
                Chapter~\ref{sec:PSasGame} is largely dedicated to the
                careful study of strategies for Floyd's Game as well
                as a generalization of Patience Sorting called
                \emph{Two-color Patience Sorting}.

            \subsection{The Language of Pattern Avoidance}
            \label{sec:Intro:Motivation:Patterns}

                Given two positive integers $m, n \in \mathbb{Z}_{+}$,
                with $1 < m < n$, we begin this section with the
                following definition.

                \begin{Definition}
                \label{defn:ClassicalPattern}
                    Let $\pi = \pi_{1}\pi_{2}\cdots\pi_{m} \in
                    \mathfrak{S}_{m}$.  Then we say that $\pi$ is a
                    \emph{(classical permutation) pattern contained}
                    in $\sigma = \sigma_{1}\sigma_{2}\cdots\sigma_{n}
                    \in \mathfrak{S}_{n}$ if $\sigma$ contains a
                    subsequence
                    $\sigma_{i_{1}}\sigma_{i_{2}}\ldots\sigma_{i_{m}}$
                    that is order-isomorphic to $\pi$.  I.e., for each
                    $j, k \in [m]$,
                    \[
                        \sigma_{i_{j}} < \sigma_{i_{k}}
                        \mbox{ \ if and only if \ }
                        \pi_{j} < \pi_{k}.
                    \]
                \end{Definition}
                
                Despite its ostensive complexity, the importance of
                Definition~\ref{defn:ClassicalPattern} cannot be
                overstated.  In particular, the containment (and,
                conversely, avoidance) of patterns provides a
                remarkably flexible language for characterizing
                collections of permutations that share a common
                combinatorial property.  Such a point of view is
                sufficiently general that it provides the foundation
                for an entire field of study commonly called
                \emph{Pattern Avoidance}.  There are also many natural
                ways to extend the definition of a classical pattern,
                as we will discuss after the following examples.
                
                \begin{Example}
                \label{eg:ClassicalPatternExamples}
                    ~
                    \begin{enumerate}
                        \item Given $k \in \mathbb{Z}_{+}$, a pattern
                        of the form $\textrm{\emph{\i}}_k = 1 2 \cdots
                        k \in \mathfrak{S}_{k}$ is called a
                        (classical) \emph{monotone increasing}
                        pattern, and $\sigma \in \mathfrak{S}_{n}$
                        contains $\textrm{\emph{\i}}_k$ if and only if
                        the length of the longest increasing
                        subsequence in $\sigma$ satisfies
                        $\ell_{n}(\sigma) \geq k$.
                        
                        E.g., $\sigma = 364827159 \in
                        \mathfrak{S}_{9}$ has several increasing
                        subsequences of length four, and so $\sigma$
                        contains quite a few occurrences of the
                        patterns $12$, $123$, and $1234$.  Moreover,
                        as we saw in
                        Example~\ref{eg:SingleRowBumpingExample},
                        $\ell_{9}(364827159) = 4$.
                    
                        \item The permutation $\sigma = 364827159 \in
                        \mathfrak{S}_{9}$ also contains occurrences of
                        such patterns as $231 \in \mathfrak{S}_{3}$
                        (e.g., via the subsequence $271$), $2431 \in
                        \mathfrak{S}_{4}$ (e.g., via the subsequence
                        $3871$), $23541 \in \mathfrak{S}_{5}$ (e.g.,
                        via the subsequence $34871$), and $235416 \in
                        \mathfrak{S}_{6}$ (e.g., via the subsequence
                        $368719$).
                    \end{enumerate}
                    
                \end{Example}

                Given a permutation $\sigma \in \mathfrak{S}_{n}$
                containing a pattern $\pi \in \mathfrak{S}_{m}$, note
                that the components of the order-isomorphic
                subsequence
                $\sigma_{i_{1}}\sigma_{i_{2}}\cdots\sigma_{i_{m}}$ are
                not required to occur contiguously within $\sigma$.
                In other words, the difference between consecutive
                indices in the subsequence is allowed to be more than
                one.  It is sometimes convenient, though, to consider
                patterns formed using only consecutive subsequences.
                This gives rise to the definition of a \emph{block
                pattern} (a.k.a~\emph{consecutive pattern} or
                \emph{segment pattern}).

                \begin{Definition}
                \label{defn:BlockPattern}
                    Let $\pi \in \mathfrak{S}_{m}$ be a permutation.
                    Then we say that $\pi$ is a \emph{block
                    (permutation) pattern contained} in $\sigma =
                    \sigma_{1}\sigma_{2}\cdots\sigma_{n} \in
                    \mathfrak{S}_{n}$ if $\sigma$ contains a
                    subsequence $\sigma_{i}\sigma_{i +
                    1}\ldots\sigma_{i + m - 1}$ that is
                    order-isomorphic to $\pi$.  I.e., $\sigma$
                    contains $\pi$ as a classical pattern but with the
                    subsequence $\sigma_{i}\sigma_{i +
                    1}\ldots\sigma_{i + m - 1}$ necessarily comprised
                    of consecutive entries in $\sigma$.
                \end{Definition}
                
                In effect, block patterns can be seen as the most
                restrictive special case of a classical pattern, and
                both notions of pattern have numerous combinatorial
                applications.  E.g., an \emph{inversion} in a
                permutation is an occurrence of a (classical) $21$
                pattern, while a \emph{descent} is an occurrence of a
                (block) $21$ pattern.  At the same time, though, there
                is no reason to insist on either extreme.  The
                following definition, which originated from the
                classification of so-called Mahonian Statistics by
                Babson and Steingr\'{\i}msson \cite{refBabS2000},
                naturally generalizes both
                Definition~\ref{defn:ClassicalPattern} and
                Definition~\ref{defn:BlockPattern}.

                \begin{Definition}
                \label{defn:GeneralizedPattern}
                    Let $\pi \in \mathfrak{S}_{m}$ be a permutation,
                    and let $D$ be a distinguished subset $D =
                    \{d_{1}, d_{2}, \ldots, d_{k} \} \subset [m - 1]$
                    of zero or more of its subscripts.  Then we say
                    that $\pi$ is a \emph{generalized (permutation)
                    pattern contained} in $\sigma =
                    \sigma_{1}\sigma_{2}\cdots\sigma_{n} \in
                    \mathfrak{S}_{n}$ if $\sigma$ contains a
                    subsequence
                    $\sigma_{i_{1}}\sigma_{i_{2}}\ldots\sigma_{i_{m}}$
                    that is order-isomorphic to $\pi$ such that
                    \[
                        i_{j + 1} - i_{j} = 1
                        \mbox{ \ if \ }
                        j \in [m - 1] \setminus D.
                    \]
                    I.e., $\sigma$ contains $\pi$ as a classical
                    pattern but with $\sigma_{i_{j}}$ and
                    $\sigma_{i_{j + 1}}$ necessarily consecutive
                    entries in $\sigma$ unless $j \in D$ is a
                    distinguished subscript of $\pi$.
                \end{Definition}
                
                When denoting a generalized pattern $\pi \in
                \mathfrak{S}_{m}$ with distinguished subscript set $D
                = \{d_{1}, d_{2}, \ldots, d_{k} \}$, we will often
                insert a dash between the entries $\pi_{d_{j}}$ and
                $\pi_{d_{j + 1}}$ in $\pi$ for each $j \in [k]$.  We
                illustrate this convention in the following examples.
                
                \begin{Example}
                \label{eg:GeneralizedPatternExamples}
                    ~
                    \begin{enumerate}
                        \item Let $\pi \in \mathfrak{S}_{m}$ be a
                        generalized pattern.  Then $\pi$ is a
                        classical pattern if the subscript set $D = [m
                        - 1]$.  I.e.,
                        \[
                            \pi = 
                            \pi_{1}   \mathrm{-}
                            \pi_{2}   \mathrm{-}
                            \cdots    \mathrm{-}
                            \pi_{m-1} \mathrm{-}
                            \pi_{m}.
                        \]
                        Similarly, $\pi$ is a block pattern if the
                        subscript set $D = \emptyset$, in which case
                        it is written with no dashes: $\pi =
                        \pi_{1}\pi_{2}\cdots\pi_{m}$.
                    
                        \item The permutation $\sigma = 364827159 \in
                        \mathfrak{S}_{9}$ from
                        Example~\ref{eg:ClassicalPatternExamples}
                        contains many examples of generalized
                        patterns.  The subsequence $271$ is an
                        occurrence of (the block pattern) $231$, while
                        $371$ is not.  However, both of the
                        subsequences $271$ and $371$ are occurrences
                        of $2\mathrm{-}31$, while $381$ is not.  In a
                        similar way, both $361$ and $481$ are
                        occurrences of $23\mathrm{-}1$, while $381$ is
                        again not.
                        
                        Finally, we note that each of the subsequences
                        $362$, $361$, $342$, $341$, $382$, $381$,
                        $371$, $682$, $681$, $685$, $671$, $675$,
                        $482$, $481$, $471$, and $271$ form an
                        occurrence of the generalized pattern
                        $2\mathrm{-}3\mathrm{-}1$ (a.k.a.~the
                        classical pattern $231$) in $364827159$.
                    \end{enumerate}
                    
                \end{Example}

                An important further generalization of
                Definition~\ref{defn:GeneralizedPattern} requires that
                the context in which the occurrence of a generalized
                pattern occurs also be taken into account.  The
                resulting concept of \emph{barred pattern} first arose
                within the study of so-called stack-sortability of
                permutations by West \cite{refWest1990} (though West's
                barred patterns were based upon the definition of a
                classical pattern and not upon the definition of a
                generalized pattern as below).  As we will illustrated
                in Section~\ref{sec:Intro:Summary}, these barred
                patterns arise as a natural combinatorial tool in the
                study of Patience Sorting.

                \begin{Definition}
                \label{defn:BarredPattern}
                    Let $\pi = \pi_{1}\pi_{2}\cdots\pi_{m} \in
                    \mathfrak{S}_{m}$ be a generalized pattern with an
                    additional distinguished subset $B = \{b_{1},
                    b_{2}, \ldots, b_{\ell} \} \subsetneq [m]$ of zero or
                    more of its subscripts.  Then we say that $\pi$ is
                    a \emph{barred (generalized permutation) pattern
                    contained} in $\sigma =
                    \sigma_{1}\sigma_{2}\cdots\sigma_{n} \in
                    \mathfrak{S}_{n}$ if $\sigma$ contains a
                    subsequence
                    $\sigma_{i_{1}}\sigma_{i_{2}}\ldots\sigma_{i_{m}}$
                    such that
                    \begin{itemize}
                        \item the (index restricted) subsequence
                        $\sigma_{i_{1}}\sigma_{i_{2}}\ldots\sigma_{i_{m}}
                        |_{[m] \setminus S} =
                        \sigma_{i_{b_{1}}}\sigma_{i_{b_{2}}}\ldots\sigma_{i_{b_{\ell}}}$
                        is order-isomorphic to the subsequence
                        $\pi|_{[m] \setminus S} =
                        \pi_{b_{1}}\pi_{b_{2}}\cdots\pi_{b_{\ell}}$ of
                        $\pi$
                    
                        \item \textbf{and} the (unrestricted)
                        subsequence
                        $\sigma_{i_{1}}\sigma_{i_{2}}\ldots\sigma_{i_{m}}$
                        is \textbf{not} order-isomorphic to $\pi$.
                    \end{itemize}
                    I.e., $\sigma$ contains the subsequence of $\pi$
                    indexed by $[m] \setminus B$ as a generalized
                    pattern unless this occurrence of $\pi|_{[m]
                    \setminus B}$ is part of an occurrence of $\pi$ as
                    a generalized pattern.
                \end{Definition}
                
                When denoting a barred (generalized) pattern $\pi \in
                \mathfrak{S}_{m}$ with distinguished subscript set $B
                = \{b_{1}, b_{2}, \ldots, b_{\ell} \}$, we will often
                place overbars atop the entries $\pi_{b_{j}}$ in $\pi$
                for each $j \in [\ell]$.  This convention is illustrated
                in the following examples.

                \begin{Example}
                \label{eg:BarredPatternExamples}
                    ~
                    \begin{enumerate}
                        \item Let $\pi \in \mathfrak{S}_{m}$ be a
                        barred pattern.  Then $\pi$ is an ordinary
                        generalized pattern (as in
                        Definition~\ref{eg:GeneralizedPatternExamples})
                        if the subscript set $B = \emptyset$.
                    
                        \item A permutation contains an occurrence of
                        the barred pattern
                        $3\mathrm{-}\overline{1}\mathrm{-}42$ (where
                        $\pi = 3142 \in \mathfrak{S}_{4}$, $D = \{1, 2
                        \}$, and $B = \{ 2 \}$) if it contains an
                        occurrence of the generalized pattern
                        $2\mathrm{-}31$ that is not part of an
                        occurrence of the generalized pattern
                        $3\mathrm{-}1\mathrm{-}42$.
                        
                        E.g., the permutation $\sigma = 364827159 \in
                        \mathfrak{S}_{9}$ from
                        Example~\ref{eg:ClassicalPatternExamples}
                        contains the generalized pattern
                        $2\mathrm{-}31$ (and hence
                        $3\mathrm{-}\overline{1}\mathrm{-}42$) via
                        each of the subsequences $382$, $371$, $682$,
                        $671$, $482$, $471$, and $271$.  This is
                        because none of these subsequences occur
                        within a larger occurrence of the generalized
                        pattern $3\mathrm{-}1\mathrm{-}42$.  In fact,
                        one can check that no subsequence of
                        $364827159$ constitutes an occurrence of the
                        pattern $3\mathrm{-}1\mathrm{-}42$.
                    \end{enumerate}
                \end{Example}

                When a permutation $\sigma \in \mathfrak{S}_{n}$ fails
                to contain a given pattern $\pi \in \mathfrak{S}_{m}$,
                then we say that $\sigma$ \emph{avoids} $\pi$.
                Questions of pattern avoidance have tended to dominate
                the study of patterns in general.  This is at least
                partially due to the initial choice of terminology
                \emph{forbidden subsequence}, which was only later
                refined to distinguish between \emph{pattern
                containment} and \emph{pattern avoidance}.
                Nonetheless, such bias continues to be reflected in
                the following (now standard) definitions.

                \begin{Definition}
                    Given any collection $\pi^{(1)}, \pi^{(2)},
                    \ldots, \pi^{(k)}$ of patterns, we denote by
                    \[
                        S_{n}(\pi^{(1)}, \pi^{(2)},\ldots, \pi^{(k)})
                        =
                        \bigcap_{i=1}^{k}S_{n}(\pi^{(i)})
                        = \bigcap_{i=1}^{k}
                        \{\sigma \in \mathfrak{S}_{n} \ | \
                          \sigma \textrm{ avoids } \pi^{(i)}
                        \}
                    \]
                    the \emph{avoidance set} consisting of all
                    permutations $\sigma \in \mathfrak{S}_{n}$ such
                    that $\sigma$ simultaneously avoids each of the
                    patterns $\pi^{(1)}, \pi^{(2)}, \ldots,
                    \pi^{(k)}$.  Furthermore, the set
                    \[
                        Av(\pi^{(1)}, \pi^{(2)}, \ldots, \pi^{(k)})
                        =
                        \bigcup_{n \geq 1} S_{n}(\pi^{(1)}, \pi^{(2)}, \ldots, \pi^{(k)})
                    \]
                    is called the \emph{(pattern) avoidance class}
                    with \emph{basis} $\{\pi^{(1)}, \pi^{(2)}, \ldots,
                    \pi^{(k)}\}$.  We will also refer to the sequence
                    of cardinalities
                    \[
                        |S_{1}(\pi^{(1)}, \pi^{(2)}, \ldots, \pi^{(k)})|,
                        |S_{2}(\pi^{(1)}, \pi^{(2)}, \ldots, \pi^{(k)})|,
                        \ldots,
                        |S_{n}(\pi^{(1)}, \pi^{(2)}, \ldots, \pi^{(k)})|,
                        \ldots
                    \]
                    as the \emph{avoidance sequence} for the basis
                    $\{\pi^{(1)}, \pi^{(2)}, \ldots, \pi^{(k)}\}$.
                \end{Definition}
                
                In effect, avoidance sets provide a language for
                characterizing permutations that share common
                combinatorial properties, and the associated avoidance
                sequences are commonly used to motivate the
                formulation of (often otherwise unmotivated)
                bijections with equinumerous combinatorial sets.
                Section~\ref{sec:PSasAlgorithm:Enumerating3142}
                contains a discussion along these lines for the
                avoidance set
                $S_{n}(3\mathrm{-}\overline{1}\mathrm{-}42)$, which
                turns out to be especially important in the study of
                Patience Sorting.  (See also
                Section~\ref{sec:Intro:Summary} below.)
                
                More information about permutation patterns in general
                can be found in \cite{refBona2004}.

        \section{Summary of Main Results}
        \label{sec:Intro:Summary}

            In this section, we summarize the main results in this
            Dissertation concerning Patience Sorting as an algorithm.
            Many of these results first appeared in
            \cite{refBLFPSAC05, refBLPP05, refBLFPSAC06}.
            (Chapter~\ref{sec:PSasGame} is largely concerned with
            comparing Patience Sorting to other strategies for Floyd's
            Game, so we do not summarize it here.)
            
            In Section~\ref{sec:PSasAlgorithm:PSEquivalence}, we
            define a ``column word'' operation on pile configurations
            (Definition~\ref{defn:ReversePatienceWords}) called the
            \emph{reverse patience word} ($RPW$).  This is by direct
            analogy with the column word of a standard Young tableau
            (see Example~\ref{eg:TableauColumnWordExample}).  In
            Section~\ref{sec:PSasAlgorithm:Enumerating3142}, we then
            provide the following characterization of reverse patience
            words in terms of pattern avoidance.

            \begin{RestateTheorem}{Theorem}{thm:EnumeratingSn3142}
                The set of permutations
                $S_n(3\mathrm{-}\overline{1}\mathrm{-}42)$ avoiding
                the barred pattern
                $3\mathrm{-}\overline{1}\mathrm{-}42$ is exactly the
                set $RPW(R(\mathfrak{S}_{n}))$ of reverse patience
                words obtainable from the symmetric group
                $\mathfrak{S}_{n}$.  Moreover, the number of elements
                in $S_n(3\mathrm{-}\overline{1}\mathrm{-}42)$ is given
                by the $n^{\mathrm{th}}$ Bell number
                \[
                    B_{n} =
                    \frac{1}{e}\sum_{k=0}^{\infty}\frac{k^{n}}{k!}.
                \]
            \end{RestateTheorem}
            
            This theorem is actually a special case of a more general
            construction.  In
            Section~\ref{sec:PSasAlgorithm:PSEquivalence}, we define
            the equivalence relation $\stackrel{PS}{\sim}$
            (Definition~\ref{defn:PSequivalence}) on the symmetric
            group $\mathfrak{S}_{n}$ by analogy to Knuth equivalence
            (Definition~\ref{defn:KnuthEquivalence}).  Specifically,
            two permutations $\sigma, \tau \in \mathfrak{S}_{n}$ are
            called \emph{patience sorting equivalent} (denoted $\sigma
            \stackrel{PS}{\sim} \tau$) if they result in the same pile
            configurations $R(\sigma) = R(\tau)$ under Patience
            Sorting (Algorithm~\ref{alg:MallowsPSprocedure}).  We then
            provide the following characterization of this equivalence
            relation:

            \begin{RestateTheorem}{Theorem}{thm:PSequivalence}
               Let $\sigma, \tau \in \mathfrak{S}_{n}$.  Then $\sigma
               \stackrel{PS}{\sim} \tau$ if and only if $\sigma$ and
               $\tau$ can be transformed into the same permutation by
               changing one or more occurrences of the pattern
               $2\mathrm{-}31$ into occurrences of the pattern
               $2\mathrm{-}13$ such that none of these $2\mathrm{-}31$
               patterns are contained within an occurrence of a
               $3\mathrm{-}1\mathrm{-}42$ pattern.

               In other words, $\stackrel{PS}{\sim}$ is the
               equivalence relation generated by changing
               $3\mathrm{-}\overline{1}\mathrm{-}42$ patterns into
               $3\mathrm{-}\overline{1}\mathrm{-}24$ patterns.
            \end{RestateTheorem}
            
            In Section~\ref{sec:PSasAlgorithm:Invertibility}, we then
            use this result to additionally characterize those
            permutations within singleton equivalence classes under
            $\stackrel{PS}{\sim}$.  (We also enumerate the resulting
            avoidance set, which involves convolved Fibonacci numbers
            as defined in \cite{refOEIS}.)

            \begin{RestateTheorem}{Theorem}{thm:PermutationsWithUniquePileConfigurations}
                A pile configuration pile $R$ has a unique preimage
                $\sigma \in \mathfrak{S}_{n}$ under Patience Sorting
                if and only if $\sigma \in
                S_n(3\mathrm{-}\bar{1}\mathrm{-}42,3\mathrm{-}\bar{1}\mathrm{-}24)$.
            \end{RestateTheorem}
            
            In Chapter~\ref{sec:ExtendingPS}, we next turn our
            attention to Extended Patience Sorting
            (Algorithm~\ref{alg:ExtendedPSalgorithm}).  Specifically,
            we first characterize the resulting stable pairs in
            Section~\ref{sec:ExtendingPS:StablePairs} using further
            pattern avoidance (where $S'$ is the same notation as
            in Example~\ref{eg:PileReflectionExample}):

            \begin{RestateTheorem}{Theorem}{thm:ExtendedPSbijection}
                Extended Patience Sorting gives a bijection between
                the symmetric group $\mathfrak{S}_{n}$ and the set of
                all ordered pairs of pile configurations $(R, S)$ such
                that both $\mathrm{sh}(R) = \mathrm{sh}(S)$ and
                $(RPW(R),RPW(S'))$ avoids simultaneous occurrences of
                the pairs of patterns $(31\mathrm{-}2,13\mathrm{-}2)$,
                $(31\mathrm{-}2,32\mathrm{-}1)$ and
                $(32\mathrm{-}1,13\mathrm{-}2)$ at the same positions
                in $RPW(R)$ and $RPW(S')$.
            \end{RestateTheorem}
            
            We also give a geometric realization for Extended Patience
            Sorting that is, in the sense described in
            Section~\ref{sec:ExtendingPS:GeometricPS}, naturally dual
            to the Geometric RSK Correspondence reviewed in
            Section~\ref{sec:ExtendingPS:GeometricPS:NEshadows}.

            \begin{RestateTheorem}{Theorem}{thm:GeometricPS}
                The Geometric Patience Sorting process described in
                Section~\ref{sec:ExtendingPS:GeometricPS:SWshadows}
                yields the same pair of pile configurations as
                Extended Patience Sorting.
            \end{RestateTheorem}
            
            Unlike the Geometric RSK Correspondence, though, Geometric
            Patience Sorting can result in intersecting lattice paths.
            In
            Section~\ref{sec:ExtendingPS:GeometricPS:CharacterizingCrossings},
            we provide the following characterization of those
            permutations that do no result in intersecting lattice
            paths.

            \begin{RestateTheorem}{Theorem}{thm:NoncrossingPilesCondition}
                Geometric Patience Sorting applied to $\sigma \in
                \mathfrak{S}_{n}$ results in non-crossing lattice
                paths at each iteration of the algorithm if and only
                if every row in both $R(\sigma)$ and $S(\sigma)$ is
                monotone increasing from left to right.
            \end{RestateTheorem}
                        
    %
    %

    \newchapter{Patience Sorting as a Card Game: Floyd's Game}{Patience Sorting as a Card Game: Floyd's Game}{Patience Sorting as a Card Game: Floyd's Game}
    \label{sec:PSasGame}

            As discussed in Section~\ref{sec:Intro:Motivation:PS},
            Patience Sorting (Algorithm~\ref{alg:MallowsPSprocedure})
            can be simultaneously viewed as a card sorting algorithm,
            as a tangible realization of the length of the longest
            increasing subsequence in a permutation, and as an optimal
            strategy for Floyd's Game (Card
            Game~\ref{game:FloydsGame}).  From this last standpoint,
            there are two natural directions of study: comparing
            Patience Sorting to other strategies for Floyd's Game and
            appropriately modifying Patience Sorting so that it
            becomes a strategy for generalizations of Floyd's Game.
            Note that the most important property of any such strategy
            is the number of piles formed since this is the statistic
            that determines whether or not the strategy is optimal.
        
            In Section~\ref{sec:PSasGame:KlondikeModel}, we explicitly
            describe how Floyd's Game should be viewed as an idealized
            model for Klondike Solitaire.  Then, in
            Section~\ref{sec:PSasGame:Strategies}, we motive two
            particular optimal strategies.  These are the \emph{Greedy
            Strategy}, as defined in
            Section~\ref{sec:PSasGame:GreedyStrategies}, and the
            \emph{Look-ahead from Right Strategy}, which we introduce
            in Section~\ref{sec:PSasGame:LookAhead}.  Even though the
            Greedy Strategy in
            Section~\ref{sec:PSasGame:GreedyStrategies} is
            conceptually more general than
            Algorithm~\ref{alg:MallowsPSprocedure} (which is also
            called the ``Greedy Strategy'' in \cite{refAD1999} but
            which we will call the \emph{Simplified Greedy Strategy}),
            we will see in Section~\ref{sec:PSasGame:SimplifiedGreedy}
            that both strategies are fundamentally the same when
            applied to Floyd's Game.
            
            The distinction between Mallows' ``patience sorting
            procedure'' and the Greedy Strategy of
            Section~\ref{sec:PSasGame:GreedyStrategies} becomes
            apparent when both strategies are extended to
            generalizations of Floyd's Game.  In particular, for the
            game Two-color Patience Sorting introduced in
            Section~\ref{sec:PSasGame:TCPS}, the Greedy Strategy
            remains optimal, \emph{mutatis mutandis}, while the
            Simplified Greedy Strategy does not.
    
        \section[Floyd's Game as an Idealized Model for Klondike Solitaire]{Floyd's Game as an Idealized Model for\\ Klondike Solitaire}
        \label{sec:PSasGame:KlondikeModel}

            One of the most popular games in the world is what many
            people commonly refer to as either \emph{Solitaire} (in
            the U.S.) or \emph{Patience} (in the U.K.).  Properly
            called \emph{Klondike Solitaire} (and also sometimes
            called \emph{Demon Patience} or \emph{Fascination}), this
            is the game that many card players in the English-speaking
            world use to while away the hours if left alone with a
            deck of playing cards; and yet, no one knows the odds of
            winning or if there is even an optimal strategy for
            maximizing one's chance of winning.
            
            Klondike Solitaire is played with a standard deck of 52
            cards, with each card uniquely labeled by a combination of
            suit and rank.  There are four suits ($\clubsuit,
            \heartsuit, \spadesuit, \diamondsuit$) and thirteen ranks,
            which are labeled in increasing order as
            \begin{center}
                $A$ (for \emph{Ace}), 2, 3, \ldots, 10, $J$ (for
                \emph{Jack}), $Q$ (for \emph{Queen}), and $K$ (for
                \emph{King}).
            \end{center}
            The majority of actual gameplay involves placing cards
            into a so-called \emph{tableau} that consists of (at most)
            seven piles.  Each pile starts with exactly one face-up
            card in it, and additional cards are then placed atop
            these piles according to the rule that the ranks must
            decrease in order (from $K$ to $A$) and that cards must
            alternate between red suits ($\heartsuit, \diamondsuit$)
            and black suits ($\clubsuit, \spadesuit$).
            A player is also allowed to move
            multiple cards between piles, and it has been frequently
            suggested (see, e.g., \cite{refAD1999},
            \cite{refKuykendall1999}, \cite{refPeterson1999}) that
            this form of player choice is responsible for the
            difficulties encountered in analyzing Klondike Solitaire.
            
                    
            Floyd's Game can be viewed as a particularly simplistic
            idealized model for Klon\-dike Solitaire.  In particular,
            Floyd's Game abstracts the formation of piles with
            descending value (yet allows gaps in the sequence
            $n,\ldots, 1$), and all notions of card color and repeated
            rank values are eliminated.  Furthermore, the objective in
            Floyd's Game it to end with as few piles as possible,
            while Klondike Solitaire is concerned with forming four
            so-called \emph{foundation piles}.  One also expects to
            end up with many piles under Floyd's Game since cards
            cannot be moved once they are placed in a pile.
            
            When viewed as an abstraction of Klondike Solitaire, there
            are two natural directions for generalizing Floyd's Game
            so that it to be more like Klondike Solitaire.  Namely,
            one can introduce either repeated card ranks or card
            colors (or both).  In the former case, analysis of pile
            formation is relatively straightforward, though care must
            be taken to specify whether or not identical ranks can be
            played atop each other.  (Both ``repeated rank'' variants
            are briefly considered in \cite{refAD1999} and
            \cite{refFulman2000} under the names ``ties allowed'' and
            ``ties forbidden''.)  In the ``card colors'' case, though,
            analyzing pile formation becomes significantly more
            subtle.  We introduce a two-color generalization of
            Floyd's Game called \emph{Two-color Patience Sorting} in
            Section~\ref{sec:PSasGame:TCPS} below.

        \section[Strategies for Floyd's Game and their Optimality]{Strategies for Floyd's Game and their\\ Optimality}
        \label{sec:PSasGame:Strategies}

            There are essentially two broad classes of strategies with
            which one can play Floyd's Game: \emph{look-ahead} and
            \emph{non-look-ahead} strategies.  In the former class,
            one is allowed to take into account the structure of the
            entire deck of cards when deciding the formation of card
            piles.  When actually playing a game, though, such
            strategies are often undesirable in comparison to
            non-look-ahead strategies.  By taking a more restricted
            view of the deck, non-look-ahead strategies require that
            each card be played in order without detailed
            consideration of the cards remaining to be played.  In
            this sense, look-ahead strategies are \emph{global}
            algorithms applied to the deck of cards, while
            non-look-ahead strategies are \emph{local} algorithms.
        
            One might intuitively guess that look-ahead strategies are
            superior to non-look-ahead strategies since ``looking
            ahead'' eliminates any surprises about how the remainder
            of the deck of cards is ordered.  Thus, it may come as a
            surprise that look-ahead strategies can do no better than
            non-look-ahead strategies when playing Floyd's Game.  (We
            also prove an analogous result for Two-color Patience
            Sorting in Section~\ref{sec:PSasGame:TCPS}).
            Specifically, the prototypical non-look-ahead Greedy
            Strategy defined in
            Section~\ref{sec:PSasGame:GreedyStrategies} is easily
            shown to be \emph{optimal} using a so-called ``strategy
            stealing'' argument.  In other words, the Greedy Strategy
            always results in the fewest number of piles possible.
            
            Even though one cannot form fewer piles under a look-ahead
            strategy for Floyd's Game, it can still be useful to
            compare optimal look-ahead strategies to the Greedy
            Strategy.  In particular, the \emph{Look-ahead from Right}
            Strategy introduced in
            Section~\ref{sec:PSasGame:LookAhead} yields exactly the
            same piles as the Greedy Strategy, but it also brings
            additional insight into the structure of these piles since
            pile formation is more combinatorial.
        
        \section{The Greedy Strategy for Floyd's Game}
        \label{sec:PSasGame:GreedyStrategies}

            Under the Greedy Strategy for Floyd's Game, one plays as
            follows.
            
            \begin{Strategy}[Greedy Strategy for Floyd's Game]
            \label{alg:GreedyStrategyforPS}
                Given $\sigma~=~\sigma_{1}\sigma_{2}\cdots\sigma_{n}
                \in \mathfrak{S}_{n}$, construct the set of piles
                $\{p_{i}\}$ by
                
                \begin{enumerate}
                    \item first forming a new pile $p_{1}$ with top
                    card $\sigma_{1}$.
                    
                    \item Then, for $l = 2, \ldots, n$, suppose that
                    $\sigma_{1}, \sigma_{2}, \cdots, \sigma_{l-1}$
                    have been used to form the piles
                    \begin{center}
                        \begin{tabular}{l l l l}
                            $p_{1} = \left\{ \begin{array}{l}
                                        \sigma_{1 s_{1}} \\
                                        \vdots \\
                                        \sigma_{1 1}
                                    \end{array}
                                    \right.$,
                            &
                            $p_{2} = \left\{ \begin{array}{l}
                                        \sigma_{2 s_{2}} \\
                                        \vdots \\
                                        \sigma_{2 1}
                                    \end{array}
                                    \right.$,
                            &       
                            \ldots \ ,
                            &
                            $p_{k} = \left\{ \begin{array}{l}
                                        \sigma_{k s_{k}} \\
                                        \vdots \\
                                        \sigma_{k 1}
                                    \end{array}
                                    \right.$.
                        \end{tabular}
                    \end{center}
                    \begin{enumerate}
                        \item If $\sigma_{l} > \sigma_{j s_{j}}$ for
                        each $j = 1, \ldots, k$, form a new pile
                        $p_{k+1}$ with top card $\sigma_{l}$.
                        
                        \item Otherwise redefine pile $p_{m}$ to be
                            \begin{displaymath}
                                p_{m} =
                                    \left\{
                                        \begin{array}{l} \sigma_{l} \\
                                            \sigma_{m s_{m}} \\
                                            \vdots \\
                                            \sigma_{m 1}
                                        \end{array}
                                    \right.  \ \ \mathrm{where} \ \
                                    \sigma_{m s_{m}} = \min_{1 \, \leq
                                    \, j \, \leq \, k} \{ \sigma_{j
                                    s_{j}} \ | \ \sigma_{l} <
                                    \sigma_{j s_{j}} \}.
                            \end{displaymath}
                    \end{enumerate}
                \end{enumerate}
                
            \end{Strategy}
            
            \noindent In other words,
            Strategy~\ref{alg:GreedyStrategyforPS} initially creates
            a single pile using the first card from the deck.  Then,
            if possible, each remaining card is played atop the
            pre-existing pile having smallest top card that is larger
            than the given card.  Otherwise, if no such pile exists,
            one forms a new pile.
            
            The objective of Floyd's Game is to form as few piles as
            possible, so, intuitively, an optimal strategy is one that
            forms a new pile only when absolutely necessary.  The
            Greedy Strategy fulfills this intuition by ``eliminating''
            as few possible future plays.  We make this explicit in
            the proof of the following theorem.
            
            \begin{Theorem}
            \label{thm:GreedyStrategyOptimalForPSviaSS}
                The Greedy Strategy
                (Strategy~\ref{alg:GreedyStrategyforPS}) is an optimal
                strategy for Floyd's Game (Game~\ref{game:FloydsGame}) in
                the sense that it forms the fewest number of piles
                possible under any strategy.
            \end{Theorem}
            
            \begin{proof}
                We use an inductive \emph{strategy stealing} argument
                to show that the position in which each card is played
                under the Greedy Strategy cannot be improved upon so
                that fewer piles are formed: Suppose that, at a given
                moment in playing according to the Greedy Strategy,
                card $c$ will be played atop pile $p$; suppose
                further that, according to some optimal strategy $S$,
                card $c$ is played atop pile $q$.  We will show that
                it is optimal to play $c$ atop $p$ by ``stealing'' the
                latter strategy.
                
                Denote by $c_p$ the card currently atop pile $p$ and
                by $c_q$ the card currently atop pile $q$.  Since the
                Greedy Strategy places each card atop the pile having
                top card that is both larger than $c$ and smaller than
                all other top cards that are larger than $c$, we have
                that $c < p \leq q $.  Thus, if we were to play $c$ on
                top of pile $q$, then, from that point forward, any
                card playable atop the resulting pile could also be
                played atop pile $p$.  As such, we can construct a new
                optimal strategy $T$ that mimics $S$ verbatim but with
                the roles of piles $p$ and $q$ interchanged from the
                moment that $c$ is played.  Playing card $c$ atop pile
                $p$ is therefore optimal.
                
                Since the above argument applied equally well to each
                card in the deck, it follows that no strategy can form
                fewer piles than are formed under the Greedy Strategy.
            \end{proof}
            
            \begin{Remark}
                Even though we call
                Strategy~\ref{alg:GreedyStrategyforPS} the ``Greedy
                Strategy'' for Floyd's Game, it is subtly different
                from what Aldous and Diaconis call the ``Greedy
                Strategy'' in \cite{refAD1999}.  The difference lies
                in how much work is performed when choosing where to
                play each card atop a pre-existing pile.  However, as
                we will see in
                Section~\ref{sec:PSasGame:SimplifiedGreedy} below, the
                distinction is actually somewhat artificial.
            \end{Remark}

        \section{Simplifying the Greedy Strategy}
        \label{sec:PSasGame:SimplifiedGreedy}

            Under the Simplified Greedy Strategy for Floyd's Game, one plays as
            follows.
            
            \begin{Strategy}[Simplified Greedy Strategy for Floyd's Game]
            \label{alg:SimplifiedGreedyStrategyforPS}
                Given a permutation
                $\sigma~=~\sigma_{1}\sigma_{2}\cdots\sigma_{n} \in
                \mathfrak{S}_{n}$, construct the set of piles
                $\{p_{i}\}$ by
                
                \begin{enumerate}
                    \item first forming a new pile $p_{1}$ with top
                    card $\sigma_{1}$.
                    
                    \item Then, for $l = 2, \ldots, n$, suppose that
                    $\sigma_{1}, \sigma_{2}, \cdots, \sigma_{l-1}$
                    have been used to form the piles
                    \begin{center}
                        \begin{tabular}{l l l l}
                            $p_{1} = \left\{ \begin{array}{l}
                                        \sigma_{1 s_{1}} \\
                                        \vdots \\
                                        \sigma_{1 1}
                                    \end{array}
                                    \right.$,
                            &
                            $p_{2} = \left\{ \begin{array}{l}
                                        \sigma_{2 s_{2}} \\
                                        \vdots \\
                                        \sigma_{2 1}
                                    \end{array}
                                    \right.$,
                            &       
                            \ldots \ ,
                            &
                            $p_{k} = \left\{ \begin{array}{l}
                                        \sigma_{k s_{k}} \\
                                        \vdots \\
                                        \sigma_{k 1}
                                    \end{array}
                                    \right.$.
                        \end{tabular}
                    \end{center}
                    \begin{enumerate}
                        \item If $\sigma_{l} > \sigma_{j s_{j}}$ for
                        each $j = 1, \ldots, k$, form a new pile
                        $p_{k+1}$ with top card $\sigma_{l}$.
                        
                        \item Otherwise redefine pile $p_{m}$ to be
                            \begin{displaymath}
                                p_{m} =
                                    \left\{
                                        \begin{array}{l} \sigma_{l} \\
                                            \sigma_{m s_{m}} \\
                                            \vdots \\
                                            \sigma_{m 1}
                                        \end{array}
                                    \right.  \ \ \mathrm{where} \ \ m
                                    = \min_{1 \, \leq \, j \, \leq \,
                                    k} \{ j \ | \ \sigma_{l} <
                                    \sigma_{j s_{j}} \}.
                            \end{displaymath}
                    \end{enumerate}
                \end{enumerate}
                
            \end{Strategy}
            
            It is not difficult to see that the Simplified Greedy
            Strategy
            (Strategy~\ref{alg:SimplifiedGreedyStrategyforPS})
            produces the same piles as the Greedy Strategy
            (Strategy~\ref{alg:GreedyStrategyforPS}).  However, before
            providing a proof of this fact, we first give an
            independent proof of the optimality of the Simplified
            Greedy Strategy.  In particular, we use the following
            two-step approach due to Aldous and Diaconis in
            \cite{refAD1999}: First, we prove a lower bound on the
            number of piles formed under an strategy.  Then we show
            that Strategy~\ref{alg:SimplifiedGreedyStrategyforPS}
            achieves this bound.

            Given a positive integer $n \in \mathbb{Z}_{+}$, recall
            that an \emph{increasing subsequence} of a permutation
            $\sigma = \sigma_{1}\sigma_{2}\cdots\sigma_{n} \in
            \mathfrak{S}_{n}$ is any subsequence $s =
            \sigma_{i_{1}}\sigma_{i_{2}}\cdots\sigma_{i_{k}}$ (meaning
            that $i_{1} < i_{2} < \cdots < i_{k}$) for which
            $\sigma_{i_{1}} < \sigma_{i_{2}} < \cdots <
            \sigma_{i_{k}}$.  Note that, while a permutation may have
            many longest increasing subsequences, the length of the
            longest increasing subsequence of $\sigma$ (which we
            denote by $\ell_{n}(\sigma)$) is nonetheless well-defined.
            
            \begin{Lemma}
            \label{lem:PSPilesBound}
                The number of piles that result from applying any
                strategy for Floyd's Game to a permutation $\sigma =
                \sigma_{1}\sigma_{2}\cdots\sigma_{n} \in
                \mathfrak{S}_{n}$ is bounded from below by the length
                $\ell_{n}(\sigma)$ of the longest increasing
                subsequence of $\sigma$.
            \end{Lemma}
            
            \begin{proof}
                Let $s =
                \sigma_{i_{1}}\sigma_{i_{2}}\cdots\sigma_{i_{k}}$ be a
                longest increasing subsequence of $\sigma$.  Since the
                cards played atop each pile must decrease in value
                while the components of $s$ increase in value, it
                follows that each component of $s$ must be played into
                a different pile.  Consequently, there must be at
                least $k$ piles, regardless of one's strategy.
            \end{proof}

            \begin{Proposition}
            \label{prop:GreedyStrategyUpperBoundPSviaLIS}
                The number of piles that result from applying the
                Simplified Greedy Strategy
                (Strategy~\ref{alg:SimplifiedGreedyStrategyforPS}) to
                a permutation $\sigma =
                \sigma_{1}\sigma_{2}\cdots\sigma_{n} \in
                \mathfrak{S}_{n}$ is bounded from above by the length
                $\ell_{n}(\sigma)$ of the longest increasing
                subsequence of $\sigma$.
            \end{Proposition}
            
            \begin{proof}
                While applying the Simplified Greedy Strategy to
                $\sigma$, we impose the following bookkeeping device:
                
                \begin{itemize}
                    \item If a card is placed on top of the left most pile, 
                    then proceed as normal.
                
                    \item If a card is placed on top of any other pile,
                    then draw an arrow from it to the current top card of 
                    the pile immediately to the left.
                \end{itemize}
                
                \noindent When the Simplified Greedy Strategy
                terminates, denote by $\sigma_{i_{k}}$ (~$ = \sigma_{k
                s_{k}}$) the top card of the right most pile.  Then
                there must exist a sequence of arrows from right to
                left connecting the cards
                \[
                    \sigma_{i_{1}} < \sigma_{i_{2}} < \cdots < \sigma_{i_{k}},
                \]
                where $\sigma_{i_{j}}$ is contained in the
                $j^{\textrm{th}}$ pile and $k$ denotes the number of
                resulting piles.  Moreover, these cards must appear in
                increasing order from left to right in $\sigma$ by
                construction of the arrows.
            \end{proof}

            \begin{Corollary}
            \label{cor:GreedyStrategyOptimalForPSviaLIS}
                The Simplified Greedy Strategy
                (Strategy~\ref{alg:SimplifiedGreedyStrategyforPS}) is an
                optimal strategy for Floyd's Game
                (Game~\ref{game:FloydsGame}) in the sense that it
                forms the fewest number of piles possible under any
                strategy.
            \end{Corollary}
            
            \begin{proof}
                This follow from the combination of
                Proposition~\ref{prop:GreedyStrategyUpperBoundPSviaLIS}
                and Lemma~\ref{lem:PSPilesBound}.
            \end{proof}

            Finally, to see that the Greedy Strategy and the
            Simplified Greedy Strategy actually do produce the same
            piles, note first that the strategies differ only in their
            respective Step 2(b).  In particular, under the Greedy
            Strategy, one plays each card atop the pile whose top card
            is both larger than the current card and, at the same
            time, is the smallest top card among all top cards larger
            than it.  However, with a bit of thought, it is easy to
            see that the cards atop each pile naturally form an
            increasing sequence from left to right since new piles are
            always created to the right of all pre-existing piles.
            (Cf.~the formation of piles in
            Example~\ref{eg:NormalPSexample}.)  Thus, the Greedy
            Strategy reduces to playing each card as far to the left
            as possible, but this is exactly how cards are played
            under the Simplified Greedy Strategy.
            
            We have both proven the following Proposition and provided
            an independent proof of
            Theorem~\ref{thm:GreedyStrategyOptimalForPSviaSS}.
            
            \begin{Proposition}
            \label{prop:GreedyStrategiesEquivalentforPS}
                Strategies~\ref{alg:GreedyStrategyforPS} and
                \ref{alg:SimplifiedGreedyStrategyforPS} produce the
                same piles when applied to a permutation $\sigma \in
                \mathfrak{S}_{n}$.
            \end{Proposition}
        
        \section{A Look-Ahead Strategy for Floyd's Game}
        \label{sec:PSasGame:LookAhead}

            In this section, we outline yet another proof of
            Theorem~\ref{thm:GreedyStrategyOptimalForPSviaSS} by again
            constructing a strategy that yields the same piles as the
            Greedy Strategy (Strategy~\ref{alg:GreedyStrategyforPS}).
            Unlike the Greedy Strategy, though, this new strategy is
            actually a look-ahead strategy in the sense that one takes
            into account the entire structure of the deck of cards
            when forming piles.  Aside from providing a proof of
            Theorem~\ref{thm:GreedyStrategyOptimalForPSviaSS} that
            does not rely upon Lemma~\ref{lem:PSPilesBound} or
            ``strategy stealing'', this approach also has the
            advantage of resembling how one actually plays Klondike
            Solitaire.  Specifically, it consists of moving sub-piles
            around (in this case, strictly from right to left) in a
            manner that mimics the arbitrary rearrangement of piles in
            Klondike Solitaire.
            
            The following strategy essentially builds the
            left-to-right minima subsequences (see
            Definition~\ref{defn:LtoRminimaSubsequence}) of the deck
            of cards by reading through the permutation repeatedly
            from right to left.  The example that follows should make
            the usage clear.
            
            \begin{Definition}
                Given a subset $E \subset \mathbb{R}$, we define the
                \emph{minimum (positive) excluded integer} of $E$ to
                be
                \[
                    \mathrm{mex}^{+}(E) = \min(\mathbb{Z}_{+} \backslash E),
                \]
                where $\mathbb{Z}_{+}$ denotes the set of positive
                integers.
            \end{Definition}
            
            \begin{Strategy}[Look-ahead from Right Strategy for Floyd's Game]
            \label{alg:LookAheadStrategyForPS}
            
                Given a permutation $\sigma =
                \sigma_{1}\sigma_{2}\cdots\sigma_{n} \in
                \mathfrak{S}_{n}$, set $E = \{ \}$ and inductively
                refine the set of piles $\{p_{i}\}$ as follows:
                
                \begin{itemize}
                    \item (Base Step) Form initial piles $p_{1},
                    \ldots, p_{n}$, where the top card of each $p_{i}$ is
                    $\sigma_{i}$.
                
                    \item (Inductive Step) Suppose that we currently have
                    the $k$ piles
                    
                        \begin{center}
                           
                            \begin{tabular}{l l l l}
                                $p_{1} = \left\{ \begin{array}{l}
                                            \sigma_{1 s_{1}} \\
                                            \vdots \\
                                            \sigma_{1 1}
                                        \end{array}
                                        \right.$, &

                                $p_{2} = \left\{ \begin{array}{l}
                                            \sigma_{2 s_{2}} \\
                                            \vdots \\
                                            \sigma_{2 1}
                                        \end{array}
                                        \right.$, &
                                        
                                \ldots \ , &
                                
                                $p_{k} = \left\{ \begin{array}{l}
                                            \sigma_{k s_{k}} \\
                                            \vdots \\
                                            \sigma_{k 1}
                                        \end{array}
                                        \right.$,

                            \end{tabular}
                            
                        \end{center}
                        
                        \noindent and set $m = \mathrm{mex}^{+}(E)$.
                        
                        \begin{itemize}
                            
                            \item[(i)] If $m \geq n + 1$, then cease 
                            playing Floyd's Game.
                        
                            \item[(ii)] Otherwise, by construction,
                            there exists a pile $p_{l}$ having top
                            card $m$.  Thus, we iteratively combine
                            certain piles, starting at pile $p_{l}$,
                            as follows:
                            
                            \begin{enumerate}
                                \item If $\sigma_{l 1} > \sigma_{i
                                s_{i}}$ for $i = 1, 2, \ldots, l-1$
                                (i.e, the bottom card of the
                                $l^{\mathrm{th}}$ pile is larger than
                                the top card of each pile to its
                                left), then we redefine $E := E \cup
                                \{ \sigma_{1 \mu}\}$ and return to the
                                start of the inductive step (since
                                pile $p_{l}$ cannot legally be placed
                                on top of any pile to its left).
                            
                                \item Otherwise, take $p_{t}$ to be the
                                right-most pile that is both to the left
                                of $p_{l}$ and for which $\sigma_{l 1}$ is
                                larger than the bottom card of $p_{t}$. 
                                I.e.,
                                \begin{displaymath}
                                    t = \max_{1 \, \leq \, i \, < \,
                                    l} \{ i \ | \ \sigma_{i 1} <
                                    \sigma_{l 1}\}.
                                \end{displaymath}
                                    
                                \begin{itemize}
                                    
                                    \item[(a)] If $| t - l | \leq
                                    1$, then we redefine $E := E
                                    \cup \{ \sigma_{1 \mu}\}$ and
                                    return to the start of the
                                    induction step to avoid moving
                                    pile $p_{l}$ past pile
                                    $p_{t}$.
                                
                                    \item[(b)] Otherwise, we place
                                    pile $p_{l}$ atop the pile
                                    $p_{\mu}$ between piles
                                    $p_{t}$ and $p_{l}$ such that
                                    the current top card of
                                    $p_{\mu}$ is the smallest card
                                    greater than $\sigma_{l 1}$.
                                    I.e., we redefine the pile
                                    $p_{\mu}$ to be
                                    \begin{displaymath}
                                        p_{\mu} =
                                            \left\{
                                                \begin{array}{l}
                                                    p_{l}\\
                                                    p_{\mu}
                                                \end{array}
                                        \right.  \ \
                                        \mathrm{where} \ \ \mu=
                                        \min_{t \, < \, i \, < \,
                                        k} \{ i \ | \ m <
                                        \sigma_{i s_{i}}\},
                                    \end{displaymath}
                                    and then we redefine 
                                    $E := E \cup \{\sigma_{1 \mu}\}$.

                                \end{itemize}
                                        
                            \end{enumerate}

                        \end{itemize}
                               
                    \end{itemize}
                     
                \end{Strategy}

                \begin{Example}
                \label{eg:LookAheadExampleForPS}
                    Let $\sigma = 64518723 \in \mathfrak{S}_{8}$.  Then
                    one plays Floyd's Game under the Look-ahead from Right
                    Strategy (Strategy~\ref{alg:LookAheadStrategyForPS})
                    as follows:
                    
                    {\singlespacing
                    \begin{itemize}
                        
                        \item Start with $E = \{\}$ and the initial
                        piles\\
                        
                        \begin{center}
                            
                            \begin{tabular}{l l l l l l l l}
                                6 & 4 & 5 &  1 &  8 &  7 &  2 &  3 \\
                            \end{tabular}.
                            
                        \end{center}
                                      
                        \item $\mathrm{mex}^{+}(E) = 1$, so move the pile
                        ``$1$'' onto the pile ``$4$'' and then the resulting
                        pile onto the pile ``$6$''.  This results in $E = \{1,
                        4, 6\}$ and the piles\\
                        
                            \begin{center}
                                
                                \begin{tabular}{l l l l l l}
                                    1 &    &     &     &     &   \\
                                    4 &    &     &     &     &   \\
                                    6 & 5 &  8 &  7 &  2 &  3 \\
                                \end{tabular}. 
                                
                            \end{center}
                            
                        \item $\mathrm{mex}^{+}(E) = 2$, so move the pile
                        ``$2$'' onto the pile ``$5$''.  This results in $E =
                        \{1, 2, 4, 5, 6\}$ and the piles\\
                        
                            \begin{center}
                                
                                \begin{tabular}{l l l l l}
                                    1 &    &     &      &     \\
                                    4 & 2 &     &      &     \\
                                    6 & 5 &  8 &  7  &  3 \\
                                \end{tabular}.
                                
                            \end{center}

                        \item $\mathrm{mex}^{+}(E) = 3$, so move the pile
                        ``$3$'' onto the pile ``$7$'' and then the resulting
                        pile onto pile the ``$8$''.  This results in $E = \{1,
                        2, 3, 4, 5, 6, 7, 8\}$ and the piles\\
                        
                            \begin{center}
                                
                                \begin{tabular}{l l l}
                                    1 &    &  3\\
                                    4 & 2 &  7\\
                                    6 & 5 &  8\\
                                \end{tabular}.
                                
                            \end{center}

                        \item Finally, $\mathrm{mex}^{+}(E) = 9 > 8$, so
                        we cease playing Floyd's Game.
                        
                    \end{itemize}
                    }
         
            \end{Example}
            
            \begin{Remark}
                In view of how Strategy~\ref{alg:LookAheadStrategyForPS}
                was played out in the above example, it should be fairly
                clear that it is an optimal strategy for Floyd's Game.
                Though we do not include a proof, we nonetheless note that
                Geometric Patience Sorting (see
                Section~\ref{sec:ExtendingPS:GeometricPS}) can be used to
                show that Strategies~\ref{alg:GreedyStrategyforPS} and
                \ref{alg:LookAheadStrategyForPS} always result in the same
                piles.
            \end{Remark}
        
        \section{The Greedy Strategy for Two-color Patience Sorting}
        \label{sec:PSasGame:TCPS}

            The following is a natural generalization of Floyd's Game:

            \begin{CardGame}[Two-color Patience Sorting]
            \label{game:TCPS}
                Given a deck of (uniquely labeled) bi-colored cards
                $c_{1}, c_{2}, \ldots, c_{2 n}$,

                \begin{itemize}
                    \item place the first card $c_{1}$ from the deck
                    into a pile by itself.

                    \item Then, for each card $c_{i}$ ($i = 2,
                    \ldots, 2 n$), either

                        \begin{itemize}

                            \item put $c_{i}$ into a new pile by
                            itself

                            \item or play $c_{i}$ on top of any pile
                            whose current top card is larger than
                            $c_{i}$ and of opposite color.

                        \end{itemize}

                    \item The object of the game is to end with as
                    few piles as possible.

                \end{itemize}

            \end{CardGame}
            
            \noindent In other words, Two-color Patience Sorting
            (TCPS) is played exactly like Floyd's Game except that
            each card is one of two colors and the colors of the cards
            in each pile must alternate.  As discussed in
            Section~\ref{sec:PSasGame:KlondikeModel}, this
            introduction of card color results in piles that are
            significantly more like those produced while playing
            Klondike Solitaire.
                
            In order to model a bi-colored deck of cards, we introduce
            the following conventions.  Given a positive integer $n
            \in \mathbb{Z}_{+}$, we denote $[n] = \{1, 2, \ldots, n\}$
            and by $\mathfrak{S}_{n}$ the set of all permutations on
            $[n]$.  Then, as a natural extension of these conventions,
            we define the overbarred integers $[\bar{n}] = \{\bar{1},
            \bar{2}, \ldots, \bar{n}\}$ by analogy and use
            $\mathfrak{S}_{\bar{n}}$ to denote the set of all
            permutations on $[\bar{n}]$.  Letting $\mathcal{N} =
            [n]\cup[\bar{n}]$, we will use $i^{\pm}$ to denote $i \in
            \mathcal{N}$ when it is unimportant whether or not $i$ has
            an overbar.
            
            With this notation, we can now define the elements in our 
            bi-colored deck of cards as follows.
            
            \begin{Definition}
                A \emph{2-permutation} $w = \left(
                    \begin{array}{l l l l l l l}
                        1 & 2 & 3 & 4 & \cdots & 2n-1 & 2n \\
                        w_{1} & w_{2} & w_{3} & w_{4} & \cdots & w_{2n-1}
                        & w_{2n}
                    \end{array}
                \right)$ is any permutation on the set $\mathcal{N}$.
            \end{Definition}
            
            As with normal permutations, we will often also denote the
            2-permutation $w$ by the bottom row of the two-line array.
            That is, by abuse of notation, we write $w = w_{1} w_{2}
            \cdots w_{2n}$.  Furthermore, we denote by
            $\mathfrak{S}_{\mathcal{N}}$ the set of all 2-permutations
            on $\mathcal{N}$.  Note that $\mathfrak{S}_{\mathcal{N}}$
            is isomorphic to $\mathfrak{S}_{2n}$ as a set, but we
            distinguish $\mathfrak{S}_{\mathcal{N}}$ from
            $\mathfrak{S}_{2n}$ in order to emphasize the difference
            in how we will treat the former combinatorially.  In
            particular, there is no transitive order relation between
            barred and unbarred numbers.
                  
            Finally, we introduce a useful function $\mathcal{N}
            \times \mathcal{N} \rightarrow \{0, 1\}$ that makes
            dealing with a mixture of barred and unbarred numbers more
            convenient:
            
            \begin{displaymath}
                bar(i,j) = \left\{
                    \begin{array}{l @{, \ } l}
                        1 & \mathrm{\ if \ exactly \ one \ of \ } i 
                        \mathrm{\ and \ } j
                        \mathrm{\ is \ overbarred} \\
                        0 & \mathrm{\ otherwise}
                    \end{array}\right. .
            \end{displaymath}

            Now that we have established notation for our bi-colored
            deck of cards, we show that the Greedy Strategy for
            Floyd's Game can be naturally extended in order to form an
            optimal strategy for TCPS while the Simplified Greedy
            Strategy cannot.  More precisely, we begin by defining the
            following strategy.
            
            \begin{Strategy}[Na\"{\i}ve Greedy Strategy for Two-color Patience Sorting]
            \label{alg:NaiveGreedyStrategyForTCPS}
                Given a random 2-permutation $w = w_{1}w_{2}\cdots
                w_{2 n} \in \mathfrak{S}_{\mathcal{N}}$, build the set
                of piles $\{p_{i}\}$ by
                
                \begin{enumerate}
                    \item first forming a new pile $p_{1}$ with top
                    card $w_{1}$.
                    
                    \item Then, for each $l = 2, \ldots, 2n$, suppose
                    that $w_{1}, w_{2}, \cdots, w_{l-1}$, have been
                    used to form piles
                        
                    \begin{center}
                       
                        \begin{tabular}{l l l l}
                            $p_{1} = \left\{ \begin{array}{l}
                                        w_{1 s_{1}} \\
                                        \vdots \\
                                        w_{1 1}
                                    \end{array}
                                    \right.$, &

                            $p_{2} = \left\{ \begin{array}{l}
                                        w_{2 s_{2}} \\
                                        \vdots \\
                                        w_{2 1}
                                    \end{array}
                                    \right.$, &
                                    
                            \ldots \ , &
                            
                            $p_{k} = \left\{ \begin{array}{l}
                                        w_{k s_{k}} \\
                                        \vdots \\
                                        w_{k 1}
                                    \end{array}
                                    \right.$.

                        \end{tabular}
                        
                    \end{center}
                    
                    \begin{enumerate}
                        
                        \item If $w^{\pm}_{l} \geq w^{\pm}_{j s_{j}}$
                        for each $j = 1, \ldots, k$ such that
                        bar$(w_{l}, w_{j s_{j}}) = 1$, then form a new
                        pile $p_{k+1}$ with top card $w_{l}$.
                        
                        \item Otherwise, redefine pile $p_{m}$ to be
                        
                            \begin{displaymath}
                                p_{m} =
                                    \left\{
                                        \begin{array}{l} w_{l} \\
                                            w_{m s_{m}} \\
                                            \vdots \\
                                            w_{m 1}
                                        \end{array}
                                    \right.  \ \ \mathrm{where} \ \ m
                                    = \min_{1 \, \leq \, j \, \leq \,
                                    k} \{ j \ | \ w^{\pm}_{l} <
                                    w^{\pm}_{j s_{j}} , \mathrm{\ bar
                                    }(w_{l}, w_{j s_{j}}) = 1 \} .
                            \end{displaymath}                    
                        
                    \end{enumerate}
                    
                \end{enumerate}
            \end{Strategy}
            
            \noindent In other words, we play each card as far to the
            left as possible (up to card color).  We both illustrate
            this algorithm and show that it is not an optimal strategy
            for TCPS in the following example.
            
            \begin{Example}
            \label{eg:TCPSexample}
                Let $w = \bar{3}\bar{2}3\bar{1}12 \in 
                \mathfrak{S}_{[3]\cup[\bar{3}]}$.  Then
                one applies
                Strategy~\ref{alg:NaiveGreedyStrategyForTCPS} to $w$
                as follows:
            \end{Example}
            
            \begin{center}
                
                \begin{tabular}{l p{80pt} l p{80pt}}
                    After playing $\mathbf{\bar{3}}$ : &
                    \begin{tabular}{l}
                        $\mathbf{\bar{3}}$
                    \end{tabular}&
                    After playing $\mathbf{\bar{2}}$: &
                    \begin{tabular}{l l}
                        $\bar{3}$ & $\mathbf{\bar{2}}$
                    \end{tabular}
                \end{tabular}\\[12pt]

                \begin{tabular}{l p{80pt} l p{80pt}}
                    After playing $\mathbf{3}$: &
                    \begin{tabular}{l l l}
                        $\bar{3}$ & $\bar{2}$ & $\mathbf{3}$
                    \end{tabular}&
                    After playing $\mathbf{\bar{1}}$: &
                    \begin{tabular}{l l l}
                          &  & $\mathbf{\bar{1}}$ \\
                        $\bar{3}$ & $\bar{2}$ & 3
                    \end{tabular}
                \end{tabular}\\[12pt]
                    
                \begin{tabular}{l p{80pt} l p{80pt}}
                    After playing $\mathbf{1}$: &
                    \begin{tabular}{l l l}
                        $\mathbf{1}$ &  & $\bar{1}$ \\
                        $\bar{3}$ & $\bar{2}$ & 3
                    \end{tabular}&
                    After playing $\mathbf{2}$: &
                    \begin{tabular}{l l l l}
                         1 &  & $\bar{1}$  &   \\
                        $\bar{3}$ & $\bar{2}$ & 3 & $\mathbf{2}$
                    \end{tabular}
                \end{tabular}\\

            \end{center}
        
            \noindent However, note that one could also play to obtain
            the following piles (as given by
            Strategy~\ref{alg:GreedyStrategyForTCPS} below):
            
            \begin{center}
                
                \begin{tabular}{l p{80pt} l p{80pt}}
                    After playing the $\mathbf{\bar{3}}$ : &
                    \begin{tabular}{l}
                        $\mathbf{\bar{3}}$
                    \end{tabular}&
                    After playing the $\mathbf{\bar{2}}$: &
                    \begin{tabular}{l l}
                        $\bar{3}$ & $\mathbf{\bar{2}}$
                    \end{tabular}
                \end{tabular}\\[12pt]

                \begin{tabular}{l p{80pt} l p{80pt}}
                    After playing the $\mathbf{3}$: &
                    \begin{tabular}{l l l}
                        $\bar{3}$ & $\bar{2}$ & $\mathbf{3}$
                    \end{tabular}&
                    After playing the $\mathbf{\bar{1}}$: &
                    \begin{tabular}{l l l}
                          &  & $\mathbf{\bar{1}}$ \\
                        $\bar{3}$ & $\bar{2}$ & 3
                    \end{tabular}
                \end{tabular}\\[12pt]
                    
                \begin{tabular}{l p{80pt} l p{80pt}}
                    After playing the $\mathbf{1}$: &
                    \begin{tabular}{l l l}
                           & $\mathbf{1}$ & $\bar{1}$ \\
                        $\bar{3}$ & $\bar{2}$ & 3
                    \end{tabular}&
                    After playing the $\mathbf{2}$: &
                    \begin{tabular}{l l l}
                         $\mathbf{2}$ & 1 & $\bar{1}$ \\
                        $\bar{3}$ & $\bar{2}$ & 3
                    \end{tabular}
                \end{tabular}\\

            \end{center}
            
            \begin{Strategy}[Greedy Strategy for Two-color Patience Sorting]
            \label{alg:GreedyStrategyForTCPS}
                Given a random 2-permutation $w = w_{1}w_{2}\cdots
                w_{2 n} \in \mathfrak{S}_{\mathcal{N}}$, build the set
                of piles $\{p_{i}\}$ by
            
                \begin{itemize}
                    \item first forming a new pile $p_{1}$ with top
                    card $w_{1}$.
                    
                    \item Then, for $l = 2, \ldots, 2n$, suppose that
                    $w_{1}, w_{2}, \cdots, w_{l-1}$, have been used to
                    form piles
                        
                    \begin{center}
                       
                        \begin{tabular}{l l l l}
                            $p_{1} = \left\{ \begin{array}{l}
                                        w_{1 s_{1}} \\
                                        \vdots \\
                                        w_{1 1}
                                    \end{array}
                                    \right.$, &

                            $p_{2} = \left\{ \begin{array}{l}
                                        w_{2 s_{2}} \\
                                        \vdots \\
                                        w_{2 1}
                                    \end{array}
                                    \right.$, &
                                    
                            \ldots &
                            
                            $p_{k} = \left\{ \begin{array}{l}
                                        w_{k s_{k}} \\
                                        \vdots \\
                                        w_{k 1}
                                    \end{array}
                                    \right.$.

                        \end{tabular}
                        
                    \end{center}
                    
                    \begin{enumerate}
                        
                        \item If $w^{\pm}_{l} \geq w^{\pm}_{j s_{j}}$
                        for each $j = 1, \ldots, k$ such that
                        bar$(w_{l}, w_{j s_{j}}) = 1$, then form a new
                        pile $p_{k+1}$ with top card $w_{l}$.
                        
                        \item  Otherwise, redefine pile $p_{m}$ to be
                        
                            \begin{displaymath}
                                \!\!\!\!\!\!
                                p_{m} =
                                    \left\{
                                        \begin{array}{l} w_{l} \\
                                            w_{m s_{m}} \\
                                            \vdots \\
                                            w_{m 1}
                                        \end{array}
                                    \right.  \ \ \mathrm{where} \ \
                                    w_{m s_{m}} = \min_{1 \, \leq \, j
                                    \, \leq \, k} \{ w_{j s_{j}} \ | \
                                    w^{\pm}_{l} < w^{\pm}_{j s_{j}} ,
                                    \mathrm{\ bar }(w_{l}, w_{j
                                    s_{j}}) = 1 \} .
                            \end{displaymath}
                            
                    \end{enumerate}
                    
                \end{itemize}
            \end{Strategy}
            
            \noindent In other words, just as for Floyd's Game, the
            above Greedy Strategy
            (Strategy~\ref{alg:GreedyStrategyForTCPS}) starts with a
            single pile consisting of the first card from the deck.
            Then, if possible, one plays each remaining card on top of
            the pre-existing pile having smallest top card that is
            both larger than the given card and of opposite color.
            Otherwise, if no such pile exists, a new pile is formed.
            Moreover, just as with Floyd's Game, this strategy will
            again be optimal for TCPS since it forms new piles only
            when absolutely necessary.
            
            \begin{Theorem}
            \label{thm:GreedyStrategyOptimalForTCPSviaSS}
                The Greedy Strategy
                (Strategy~\ref{alg:GreedyStrategyForTCPS}) is an optimal
                strategy for Two-color Patience Sorting in the sense that
                it forms the fewest number of piles possible under any
                strategy.
            \end{Theorem}
            
            \begin{proof}   
                We use an inductive \emph{strategy stealing} argument
                to show that the position in which each card is played
                under the Greedy Strategy cannot be improved upon so
                that fewer piles are formed: Suppose that, at a given
                moment in playing according to the Greedy Strategy,
                card $c$ will be played atop pile $p$; suppose further
                that, according to some optimal strategy $S$, card $c$
                is played atop pile $q$.  We will show that it is
                optimal to play $c$ onto $p$ by ``stealing'' the
                latter strategy.
                
                Denote by $c_p$ the card atop pile $p$ and by $c_q$
                the card atop pile $q$.  Since the Greedy Strategy
                plays each card atop the pile having top card that is
                larger than $c$, of opposite color, and smaller than
                all other top cards that are both larger than $c$ and
                of opposite color, we have that $c < p \leq q $ with
                $\textrm{bar}(c, c_p) = \textrm{bar}(c, c_q) = 0$ and
                $\textrm{bar}(c_p, c_q) = 1$.  Thus, if we were to
                play $c$ atop $c_q$, then any card playable atop the
                modified pile $q$ from that point onward can also be
                played atop $c_p$.  As such, we can construct a new
                optimal strategy $T$ that mimics $S$ verbatim but with
                the roles of piles $p$ and $q$ interchanged from the
                moment that card $c$ is played.  It is therefore
                optimal to play card $c$ atop pile $p$.
                
                Since we can apply the above argument to each card in
                the deck, it follows that no other strategy can form
                fewer piles than are formed under the Greedy Strategy.
            \end{proof}
                        
            We conclude by making the following useful observation,
            which, unlike Floyd's Game, does not have an affirmative
            converse (as illustrated in the example that follows).
            
            \begin{Proposition}
            \label{prop:TCPSPilesBoundFromBelow}
                The number of piles that results from playing
                Two-color Patience Sorting under any strategy on the
                2-permutation $w = w_{1}w_{2}\cdots w_{2n} \in
                \mathfrak{S}_{\mathcal{N}}$ is bounded from below by
                the length of the longest weakly increasing
                subsequence in the word $w^{\pm}$.
            \end{Proposition}
            
            \begin{proof}
                The proof is identically to that of
                Lemma~\ref{lem:PSPilesBound}.
            \end{proof}
            
            \begin{Example}
                Let $w = 2\bar{3}\bar{1}43\bar{2}\bar{4}1 \in
                \mathfrak{S}_{[4]\cap[\bar{4}]}$.  Then, upon applying
                the Greedy Strategy
                (Strategy~\ref{alg:GreedyStrategyForTCPS}) to $w$,
                one obtains the following five piles:
                \begin{center}
                    \begin{tabular}{l l l l l}
                        $\bar{1}$ & 1 & $\bar{2}$ &  &  \\
                        2 & $\bar{3}$ & 4 & 3 & $\bar{4}$
                    \end{tabular}.
                \end{center}
                However, the length of the longest
                weakly increasing subsequence of $w^{\pm} = 23143241$
                is four, which corresponds to the unique longest
                weakly increasing subsequence $2334$.
            \end{Example}

            In fact, one can easily construct pathological examples in
            which the ``gap'' between the length of the longest weakly
            increasing subsequence and the number of piles formed
            under the Greedy Strategy grows arbitrarily large.  E.g.,
            consider
            \[
                w = n, n - 1, \ldots, 2, 1,
                \overline{n}, \overline{n - 1}, \ldots, \overline{2}, \overline{1}
                \in \mathfrak{S}_{\mathcal{N}}.
            \]
            Unlike Floyd's Game, no known combinatorial statistic on
            2-permutations is equidistributed with the number of piles
            formed when TCPS is played under the Greedy Strategy.
            Since Proposition~\ref{prop:TCPSPilesBoundFromBelow}
            yields a lower bound for the number of piles formed,
            though, one could conceivably study the asymptotic ``gap''
            between the length of the longest weakly increasing
            subsequence and the number of piles formed in order to
            gain insight into the expected number of piles formed.

    %
    %

    \newchapter{Patience as an Algorithm: Mallows' Patience Sorting Procedure}{Patience as an Algorithm: Mallows' Patience Sorting Procedure}{Patience as an Algorithm: Mallows' Patience Sorting Procedure}
    \label{sec:PSasAlgorithm}
    
        \section[Pile Configurations and Northeast Shadow Diagrams]{Pile Configurations and\\ Northeast Shadow Diagrams}
        \label{sec:PSasAlgorithm:PilesAndNEshadows}

            Given a positive integer $n \in \mathbb{Z}_{+}$, we begin
            by explicitly characterizing the objects that result when
            Patience Sorting (Algorithm~\ref{alg:MallowsPSprocedure})
            is applied to a permutation $\sigma \in \mathfrak{S}_{n}$:

            \begin{Lemma}
            \label{lem:PileConfigurationAreSetPartitions}
                Let $\sigma \in \mathfrak{S}_{n}$ be a permutation and
                $R(\sigma) = \{r_{1}, r_{2}, \ldots, r_{k}\}$ be the
                pile configuration associated to $\sigma$ under
                Algorithm~\ref{alg:MallowsPSprocedure}.  Then
                $R(\sigma)$ is a partition of the set $[n] = \{1, 2,
                \dots, n\}$ such that, denoting $r_{j} = \{ r_{j 1}>
                r_{j 2} > \dots > r_{j s_{j}}\}$,
                \begin{equation}
                \label{eqn:PileConfigurationCondition}
                    r_{j s_{j}} < r_{i s_{i}} \quad \mathrm{if}
                    \quad j < i.
                \end{equation}
                Moreover, for every set partition $T =
                \{t_{1}, t_{2}, \ldots, t_{k}\}$ satisfying
                Equation~\eqref{eqn:PileConfigurationCondition}, there
                is a permutation $\tau \in \mathfrak{S}_{n}$ such that
                $R(\tau) = T$.
            \end{Lemma}

            \begin{proof}
                Given the pile configuration $R(\sigma) = \{r_{1},
                r_{2}, \ldots, r_{k}\}$ associated to $\sigma \in
                \mathfrak{S}_{n}$, suppose that, for some pair of
                indices $i, j \in [k]$, we have that $j < i$ but $r_{j
                s_{j}} > r_{i s_{i}}$.  Then the card $r_{i s_{i}}$
                was put atop pile $r_i$ when pile $r_j$ had top card
                $d_{j} \geq r_{j s_{j}}$ so that $d_{j} > r_{i s_{i}}$
                as well.  However, it then follows that the card $r_{i
                s_{i}}$ would actually have been placed atop either
                pile $r_j$ or atop some pile to the left of $r_j$
                instead of atop pile $r_i$.  This resulting
                contradiction implies that $r_{j s_{j}} < r_{i s_{i}}$
                for each $j<i$.

                Conversely, let $T = \{t_{1}, t_{2}, \dots, t_{k}\}$
                be any set partition of $[n]$, with the block $t_{j} =
                \{ t_{j 1}> t_{j 2} > \dots > t_{j s_{j}} \}$ for each
                $j \in [k]$ and with $t_{j s_{j}} < t_{i s_{i}}$ for
                every pair of indices $i, j \in [k]$ such that $j <
                i$.  Then, setting
                \[
                    \tau= t_{1 1} t_{1 2} \cdots t_{1 s_{1}} \, t_{2 1}
                    t_{2 2} \dots t_{2 s_{2}} \; \cdots \; t_{k 1} t_{k
                    2} \cdots t_{k s_{k}},
                \]
                it is easy to see that $\tau \in  \mathfrak{S}_{n}$ 
                and that $R(\tau)=T$.
            \end{proof}

            According to
            Lemma~\ref{lem:PileConfigurationAreSetPartitions}, pile
            configurations formed from $[n]$ are set partitions of
            $[n]$ in which the constituent blocks have been ordered by
            their minimal element.
            (Cf.~Example~\ref{eg:PileConfigurationExample}.)  We
            devote the remainder of this section to an alternate
            characterization involving the \emph{(northeast) shadow
            diagram} of a permutation.  The construction of shadow
            diagrams was first used by Viennot \cite{refViennot1977}
            to study properties of the RSK Correspondence
            (Algorithm~\ref{alg:RSKAlgorithm}) for permutations.  (See
            Section~\ref{sec:ExtendingPS:GeometricPS}.)

            \begin{Definition}
            \label{defn:NEshadow}
                Given a lattice point $(m, n) \in \mathbb{Z}^{2}$, we
                define the \emph{northeast shadow} of $(m, n)$ to be the
                quarter space
                \[
                    S_{NE}(m, n)
                    =
                    \{ (x, y) \in \mathbb{R}^{2} \ | \ x \geq m, \ y \geq n\}.
                \]
            \end{Definition}

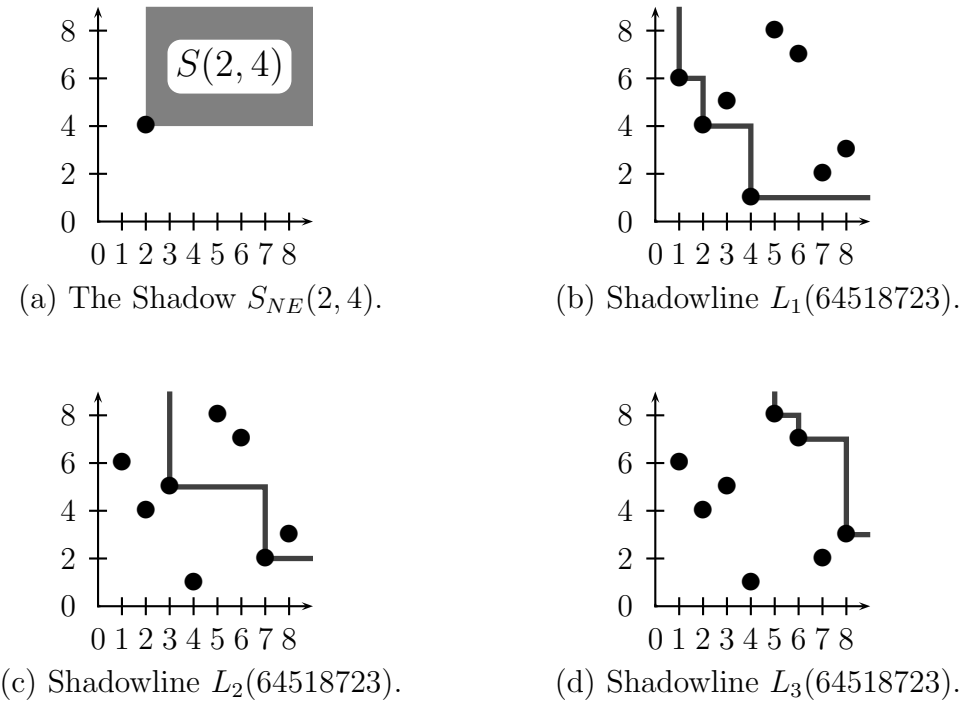
\begin{figure}[t]
            \label{fig:ShadowExample}
                \centering
                \begin{tabular}{cccc}

                    \begin{minipage}[c]{2.75in}
                        \centering

                            \psset{xunit=0.125in,yunit=0.125in}

                            \begin{pspicture}(0,0)(9,9)

                                \psaxes[Dy=2]{->}(9,9)

                                \psframe[fillcolor=gray,fillstyle=solid,linecolor=gray](2,4)(9,9)

                                \rput(5.5,6.5){\psframebox*[framearc=.5]{{\large
                                $S(2,4)$}}}

                                \rput(2,4){{\Large $\bullet$}}

                            \end{pspicture}\\[0.25in]

                          (a) The Shadow $S_{NE}(2,4)$.
                    \end{minipage}
                    &
                    \begin{minipage}[c]{2.75in}
                        \centering

                            \psset{xunit=0.125in,yunit=0.125in}

                            \begin{pspicture}(0,0)(9,9)

                                \psaxes[Dy=2]{->}(9,9)

                                \psline[linecolor=darkgray,linewidth=2pt](1,9)(1,6)(2,6)(2,4)(4,4)(4,1)(9,1)

                                \rput(1,6){{\Large $\bullet$}}%
                                \rput(2,4){{\Large $\bullet$}}%
                                \rput(3,5){{\Large $\bullet$}}%
                                \rput(4,1){{\Large $\bullet$}}%
                                \rput(5,8){{\Large $\bullet$}}%
                                \rput(6,7){{\Large $\bullet$}}%
                                \rput(7,2){{\Large $\bullet$}}%
                                \rput(8,3){{\Large $\bullet$}}%

                            \end{pspicture}\\[0.25in]

                        (b) Shadowline $L_{1}(64518723)$.
                    \end{minipage}
                      \\
                      & \vspace{0.5cm}
                      \\
                    \begin{minipage}[c]{2.75in}
                        \centering

                            \psset{xunit=0.125in,yunit=0.125in}

                            \begin{pspicture}(0,0)(9,9)

                                \psaxes[Dy=2]{->}(9,9)

                                \psline[linecolor=darkgray,linewidth=2pt](3,9)(3,5)(7,5)(7,2)(9,2)

                                \rput(1,6){{\Large $\bullet$}}%
                                \rput(2,4){{\Large $\bullet$}}%
                                \rput(3,5){{\Large $\bullet$}}%
                                \rput(4,1){{\Large $\bullet$}}%
                                \rput(5,8){{\Large $\bullet$}}%
                                \rput(6,7){{\Large $\bullet$}}%
                                \rput(7,2){{\Large $\bullet$}}%
                                \rput(8,3){{\Large $\bullet$}}%

                            \end{pspicture}\\[0.25in]

                        (c) Shadowline $L_{2}(64518723)$.
                    \end{minipage}
                    &
                    \begin{minipage}[c]{2.75in}
                        \centering

                            \psset{xunit=0.125in,yunit=0.125in}

                            \begin{pspicture}(0,0)(9,9)

                                \psaxes[Dy=2]{->}(9,9)

                                \psline[linecolor=darkgray,linewidth=2pt](5,9)(5,8)(6,8)(6,7)(8,7)(8,3)(9,3)

                                \rput(1,6){{\Large $\bullet$}}%
                                \rput(2,4){{\Large $\bullet$}}%
                                \rput(3,5){{\Large $\bullet$}}%
                                \rput(4,1){{\Large $\bullet$}}%
                                \rput(5,8){{\Large $\bullet$}}%
                                \rput(6,7){{\Large $\bullet$}}%
                                \rput(7,2){{\Large $\bullet$}}%
                                \rput(8,3){{\Large $\bullet$}}%

                            \end{pspicture}\\[0.25in]

                        (d) Shadowline $L_{3}(64518723)$.
                    \end{minipage}

                \end{tabular}\\
                
                \caption{Examples of Northeast Shadow and Shadowline Construction}

            \end{figure}

            \noindent See Figure~\ref{fig:ShadowExample}(a) for an
            example of a point's shadow.
            
            By itself, the notion of shadow doesn't come across as
            particularly exciting.  However, one can use these shadows
            in order to associate a lattice path to any (finite)
            collection of lattice points.

            \begin{Definition}
            \label{defn:NEshadowline}
                Given lattice points $(m_{1}, n_{1}), (m_{2}, n_{2}),
                \ldots, (m_{k}, n_{k}) \in \mathbb{Z}^{2}$, we define
                their \emph{northeast shadowline} to be the boundary
                of the union of the northeast shadows $S_{NE}(m_{1},
                n_{1}), S_{NE}(m_{2}, n_{2}), \ldots, S_{NE}(m_{k},
                n_{k})$.
            \end{Definition}

            In particular, we wish to associate to every permutation a
            certain collection of shadowlines (as illustrated in
            Figure~\ref{fig:ShadowExample}(b)--(d)):

            \begin{Definition}
            \label{defn:NEshadowDiagram}
                Given a permutation $\sigma =
                \sigma_{1}\sigma_{2}\cdots\sigma_{n} \in
                \mathfrak{S}_{n}$, the \emph{northeast shadow diagram}
                $D_{NE}(\sigma) = D_{NE}^{(0)}(\sigma)$ of $\sigma$
                consists of the shadowlines $L_{1}(\sigma),
                L_{2}(\sigma), \ldots, L_{k}(\sigma)$ formed as
                follows:
                \begin{itemize}
                    \item $L_{1}(\sigma)$ is the northeast shadowline
                    for the set of lattice points
                    \[
                        \{ (1, \sigma_{1}),
                           (2,\sigma_{2}),
                           \ldots,
                           (n, \sigma_{n})
                        \}.
                    \]

                    \item Then, while at least one of the points $(1,
                    \sigma_{1}), (2, \sigma_{2}), \ldots, (n,
                    \sigma_{n})$ is not contained in the shadowlines
                    $L_{1}(\sigma), L_{2}(\sigma), \ldots,
                    L_{j}(\sigma)$, define $L_{j+1}(\sigma)$ to be the
                    northeast shadowline for the points
                    \[
                        \{(i, \sigma_{i}) \ | \ i \in [n],
                          (i, \sigma_{i}) \notin \bigcup^{j}_{k=1} L_{k}(\sigma)
                        \}.
                    \]
                \end{itemize}
            \end{Definition}

            In other words, the shadow diagram $D_{NE}(\sigma) = \{
            L_{1}(\sigma), L_{2}(\sigma), \ldots, L_{k}(\sigma) \}$ of
            the permutation $\sigma \in \mathfrak{S}_{n}$ is defined
            inductively by first taking $L_{1}(\sigma)$ to be the
            shadowline for the so-called \emph{diagram} $\{(1,
            \sigma_{1}), (2, \sigma_{2}), \ldots, (n, \sigma_{n})\}$
            of the permutation $\sigma \in \mathfrak{S}_{n}$.  Then we
            ignore the points whose shadows were actually used in
            building $L_{1}(\sigma)$ and define $L_{2}(\sigma)$ to be
            the shadowline of the resulting subset of the permutation
            diagram.  We then build $L_{3}(\sigma)$ as the shadowline
            for the points not yet used in constructing either
            $L_{1}(\sigma)$ or $L_{2}(\sigma)$, and this process
            continues until each of the points in the permutation's
            diagram has been exhausted.\medskip

            One can show that the shadow diagram for a permutation
            $\sigma \in \mathfrak{S}_{n}$ encodes the top row of the
            RSK Correspondence insertion tableau $P(\sigma)$
            (resp.~recording tableau $Q(\sigma)$) as the smallest
            ordinates (resp.  smallest abscissae) of all points along
            the constituent shadowlines $L_{1}(\sigma), L_{2}(\sigma),
            \ldots, L_{k}(\sigma)$.  (See Sagan \cite{refSagan2000}
            for a proof.)  In particular, if $\sigma$ has pile
            configuration $R(\sigma) = \{r_{1}, r_{2}, \ldots,
            r_{m}\}$, then $m = k$ since the number of piles is equal
            to the length of the top row of $P(\sigma)$ (as both are
            the length of the longest increasing subsequence of
            $\sigma$; see Sections~\ref{sec:Intro:Motivation:RSK} and
            \ref{sec:Intro:Motivation:PS}).  We can say even more
            about the relationship between $D_{NE}(\sigma)$ and
            $R(\sigma)$ when both are viewed in terms of
            \emph{left-to-right minima subsequences}
            (a.k.a.~\emph{basic subsequences}).

            \begin{Definition}
            \label{defn:LtoRminimaSubsequence}
                Let $\pi = \pi_{1} \pi_{2} \cdots \pi_{l}$ be a
                partial permutation on $[n] = \{1, 2, \ldots, n\}$.
                Then the \emph{left-to-right minima subsequence} of
                $\pi$ consists of those components $\pi_{j}$ of $\pi$
                such that
                \[
                    \pi_{j} = \min_{1 \, \leq \, i \, \leq \, j}\{ \pi_{i} \}.
                \]
            \end{Definition}

            \noindent We then inductively define the left-to-right
            minima subsequences $s_{1}, s_{2}, \ldots, s_{k}$ of a
            permutation $\sigma \in \mathfrak{S}_{n}$ by first taking
            $s_{1}$ to be the left-to-right minima subsequence for
            $\sigma$ itself.  Then each subsequent $i^{\mathrm{th}}$
            \emph{left-to-right minima subsequence} $s_{i}$ is defined
            to be the left-to-right minima subsequence for the partial
            permutation obtained by removing the elements of $s_{1},
            s_{2}, \ldots, s_{i-1}$ from $\sigma$.
            
            We are now in a position to give a particularly nice
            correspondence between the piles formed under Patience
            Sorting and the shadowlines that constitute the shadow
            diagram of a permutation via these left-to-right minima
            subsequences.  We will rely heavily upon this
            correspondence in the sections that follow.

            \begin{Lemma}
            \label{lem:ShadowDiagramPileCorrespondence}
                Suppose $\sigma \in \mathfrak{S}_{n}$ has shadow
                diagram $D_{NE}(\sigma) = \{ L_{1}(\sigma), \ldots,
                L_{k}(\sigma) \}$.  Then the ordinates of the
                southwest corners of each $L_{i}$ are exactly the
                cards in the $i^{\mathrm{th}}$ pile $r_{i} \in
                R(\sigma)$ formed by applying Patience Sorting
                (Algorithm~\ref{alg:MallowsPSprocedure}) to $\sigma$.
                In other words, the $i^{\mathrm{th}}$ pile $r_{i}$ is
                exactly the $i^{\mathrm{th}}$ left-to-right minima
                subsequence of $\sigma$.
            \end{Lemma}

            \begin{proof}
                The $i^{\mathrm{th}}$ left-to-right minima subsequence
                $s_{i}$ of $\sigma$ consists of the entries in
                $\sigma$ that appear at the end of an increasing
                subsequence of length $i$ but not at the end of an
                increasing subsequence of length $i+1$.  Thus, since
                each element added to a pile must be smaller than all
                other elements already in the pile, $s_{1} = r_{1}$.
                It then follows similarly by induction that $s_{i} =
                r_{i}$ for each $i = 2, \ldots, k$.

                The proof that the ordinates of the southwest corners
                of the shadowlines $L_{i}$ are also exactly the
                elements of the left-to-right minima subsequences
                $s_{i}$ is similar.
            \end{proof}
            
            We conclude this section with an example.
            
            \begin{Example}
            \label{eg:NESalientPointExample}
                Consider $\sigma = 64518723 \in \mathfrak{S}_{8}$.  From
                Figure~\ref{fig:ShadowExample}, we see that
                \[
                    D_{NE}(\sigma) = \{ L_{1}(\sigma), L_{2}(\sigma), L_{3}(\sigma) \},
                \]
                where $L_{1}(\sigma)$ has southwest corners $\{ (1,6),
                (2,4), (4,1) \}$, $L_{2}(\sigma)$ has southwest
                corners $\{ (3,5), (7,2) \}$, and $L_{3}(\sigma)$ has
                southwest corners $\{ (5,8), (6,7), (8,3) \}$.
                Moreover, one can check that $\sigma$ has
                left-to-right minima subsequence $s_{1} = 641$
                (corresponding to the ordinates of the southwest
                corners for $L_{1}(\sigma)$), $s_{2} = 52$
                (corresponding to the ordinates of the southwest
                corners for $L_{2}(\sigma)$), and $s_{3} = 873$
                (corresponding to the ordinates of the southwest
                corners for $L_{3}(\sigma)$).
                
                Similarly, $R(\sigma) = \{ s_{1}, s_{2}, s_{3} \}$ as
                in Example~\ref{eg:PileConfigurationExample}.
            \end{Example}
        
        \section[Reverse Patience Words and Patience Sorting Equivalence]{Reverse Patience Words and\\ Patience Sorting Equivalence}
        \label{sec:PSasAlgorithm:PSEquivalence}

            In the proof of
            Lemma~\ref{lem:PileConfigurationAreSetPartitions}, we gave
            the construction of a special permutation that can be used
            to generate any set partition under Patience Sorting
            (Algorithm~\ref{alg:MallowsPSprocedure}).  At the same
            time, though, it should be clear that there are in general
            many permutations resulting in a given set partition.  In
            this section, we characterize the corresponding
            equivalence relation on the symmetric group
            $\mathfrak{S}_{n}$.  We also characterize the most natural
            choice of generators for the resulting equivalence
            classes.

            \begin{Definition}
            \label{defn:ReversePatienceWords}
                Given a pile configuration $R = \{ r_{1}, \ldots,
                r_{k} \}$, the \emph{reverse patience word} $RPW(R)$
                of $R$ is the permutation formed by concatenating the
                piles $r_{1}, r_{2}, \ldots, r_{k}$ together, with
                each pile $r_{j} = \{ r_{j 1} > r_{j 2} > \cdots >
                r_{j s_{j}} \}$ written in decreasing order.  Using
                the notation of
                Lemma~\ref{lem:PileConfigurationAreSetPartitions},
                \[
                    RPW(R) = r_{1 1} r_{1 2} \cdots r_{1 s_{1}} \
                    r_{2 1} r_{2 2} \cdots r_{2 s_{2}} \; \cdots \;
                    r_{k 1} r_{k 2} \cdots r_{k s_{k}}.
                \]
            \end{Definition}

            \begin{Example}
                The pile configuration $R = \{\{6 > 4 > 1\}, \{5 >
                2\}, \{8 > 7 > 3\}\}$ from
                Example~\ref{eg:PileConfigurationExample} is
                represented by the piles
                
                {\singlespacing
                \begin{center}
                    \begin{tabular}{l l l}
                        1 &   & 3 \\
                        4 & 2 & 7 \\
                        6 & 5 & 8 \end{tabular}\\
                \end{center}
                }
                
                \noindent and has reverse patience word $RPW(R) =
                64152873$.  Moreover, as in the proof of
                Lemma~\ref{lem:PileConfigurationAreSetPartitions},
                \[
                    R(RPW(R))=R(64152873) = R.
                \]
                This illustrates the following Lemma.
            \end{Example}

            \begin{Lemma}
            \label{lem:ReversePatienceWordEquivalence}
                Given a permutation $\sigma \in \mathfrak{S}_{n}$,
                $R(RPW(R(\sigma))) = R(\sigma)$.
            \end{Lemma}

            \begin{proof}
                Suppose that $T$ is a set partition of $[n]$
                satisfying
                Equation~\eqref{eqn:PileConfigurationCondition}, and
                let $\tau = RPW(T)$ as in the proof of
                Lemma~\ref{lem:PileConfigurationAreSetPartitions}.
                Then $T = R(\tau) = R(RPW(T))$.  In particular, given
                $\sigma \in \mathfrak{S}_{n}$, $R(\sigma) = T =
                R(RPW(R(\sigma)))$.
            \end{proof}
                
            As with the column word operation on standard Young
            tableaux (from Example~\ref{eg:TableauColumnWordExample}),
            Lemma~\ref{lem:ReversePatienceWordEquivalence} can also be
            recast from an algebraic point of view.  Denote by
            $\mathfrak{P}_{n}$ the set of all pile configurations with
            some composition shape $\gamma \models n$.  Patience
            Sorting and the reverse patience word operation can then
            be viewed as maps $R: \mathfrak{S}_{n} \to
            \mathfrak{P}_{n}$ and $RPW: \mathfrak{P}_{n} \to
            \mathfrak{S}_{n}$, respectively.  With this notation,
            Lemma~\ref{lem:ReversePatienceWordEquivalence} becomes
            
            \begin{ReuseTheorem}{Lemma}{lem:ReversePatienceWordEquivalence}
                The composition $R \circ RPW$ is the
                identity map on the set $\mathfrak{P}_{n}$.
            \end{ReuseTheorem}
            
            \noindent In particular, even though
            $RPW(R(\mathfrak{S}_{n})) = RPW(\mathfrak{P}_{n})$ is a
            proper subset of the symmetric group $\mathfrak{S}_{n}$,
            we nonetheless have that $R(RPW(R(\mathfrak{S}_{n}))) =
            R(\mathfrak{S}_{n}) = \mathfrak{P}_{n}$.  As such, it
            makes sense to define the following non-trivial
            equivalence relation on $\mathfrak{S}_{n}$, with each
            element of $RPW(\mathfrak{P}_{n})$ being the most natural
            choice of representative for the distinct equivalence
            class to which it belongs.

            \begin{Definition}
            \label{defn:PSequivalence}
                Two permutations $\sigma, \tau \in \mathfrak{S}_{n}$
                are said to be \emph{patience sorting equivalent},
                written $\sigma \stackrel{PS}{\sim} \tau$, if they
                yield the same pile configuration $R(\sigma) =
                R(\tau)$ under Patience Sorting
                (Algorithm~\ref{alg:MallowsPSprocedure}).  We denote
                the equivalence class generated by $\sigma$ as
                $\widetilde{\sigma}$.
            \end{Definition}

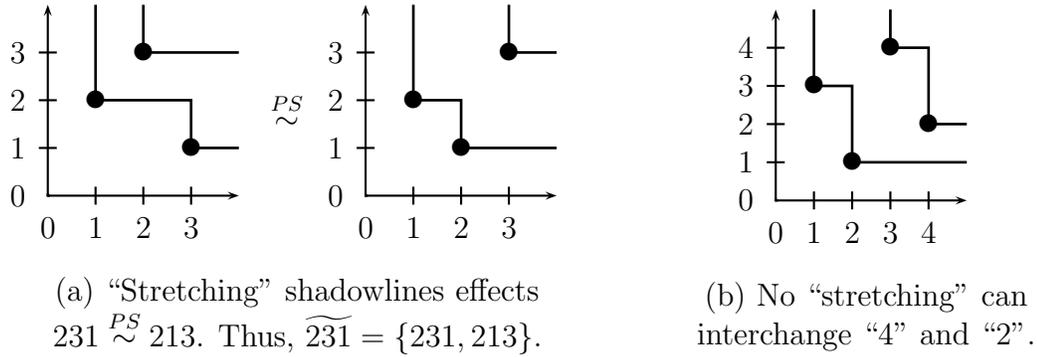
\begin{figure}[t]
            \label{fig:PSequivalenceExamples}
                \centering
                \begin{tabular}{cc}

                    \begin{minipage}[c]{3.5in}
                        \centering
                        \begin{tabular}{c c c}

                            \psset{xunit=0.25in,yunit=0.25in}

                            \begin{pspicture}(0,0)(4,4)

                                \psaxes{->}(4,4)

                                \psline[linecolor=black,linewidth=1pt](1,4)(1,2)(3,2)(3,1)(4,1)

                                \psline[linecolor=black,linewidth=1pt](2,4)(2,3)(4,3)

                                 \rput(1,2){{\Large $\bullet$}}%
                                 \rput(2,3){{\Large $\bullet$}}%
                                 \rput(3,1){{\Large $\bullet$}}%

                            \end{pspicture}

                             & \raisebox{0.35in}{$\stackrel{PS}{\sim}$ \
                             } &

                            \psset{xunit=0.25in,yunit=0.25in}

                            \begin{pspicture}(0,0)(4,4)

                                \psaxes{->}(4,4)

                                \psline[linecolor=black,linewidth=1pt](1,4)(1,2)(2,2)(2,1)(4,1)

                                \psline[linecolor=black,linewidth=1pt](3,4)(3,3)(4,3)

                                 \rput(1,2){{\Large $\bullet$}}%
                                 \rput(2,1){{\Large $\bullet$}}%
                                 \rput(3,3){{\Large $\bullet$}}%

                            \end{pspicture}\\[0.35in]

                        \end{tabular}

                        (a) ``Stretching'' shadowlines effects \\ $231
                        \stackrel{PS}{\sim} 213$.  Thus, $\widetilde{231} =
                        \{231, 213\}$.
                    \end{minipage}
                    &
                    \begin{minipage}[c]{2.125in}
                        \centering

                            \psset{xunit=0.2in,yunit=0.2in}

                            \begin{pspicture}(0,0)(5,5)

                                \psaxes{->}(5,5)

                                \psline[linecolor=black,linewidth=1pt](1,5)(1,3)(2,3)(2,1)(5,1)

                                \psline[linecolor=black,linewidth=1pt](3,5)(3,4)(4,4)(4,2)(5,2)

                                 \rput(1,3){{\Large $\bullet$}}%
                                 \rput(2,1){{\Large $\bullet$}}%
                                 \rput(3,4){{\Large $\bullet$}}%
                                 \rput(4,2){{\Large $\bullet$}}%

                            \end{pspicture}\\[0.35in]

                            (b) No ``stretching'' can interchange ``4''
                            and ``2''.
                    \end{minipage}

                    \end{tabular}

                \caption{Examples of patience sorting equivalence and non-equivalence}

            \end{figure}

            From Lemma~\ref{lem:ShadowDiagramPileCorrespondence}, we
            know that the pile configurations $R(\sigma)$ and
            $R(\tau)$ correspond to the shadow diagrams of $\sigma,
            \tau \in \mathfrak{S}_{n}$, respectively.  Thus, it should
            be intuitively clear that preserving a given pile
            configuration is equivalent to preserving the ordinates
            for the southwest corners of the shadowlines.  In
            particular, this means that we are limited to horizontally
            ``stretching'' shadowlines up to the point of not allowing
            them to cross.  This is illustrated in
            Figure~\ref{fig:PSequivalenceExamples} and the following
            examples.

        \begin{Example}
            \label{eg:PSequivalenceExample}
                ~
                \begin{enumerate}
                    \item The only non-singleton patience sorting
                    equivalence class for $\mathfrak{S}_{3}$ consists
                    of $\widetilde{231} = \{231, 213\}$.  We
                    illustrate $231 \stackrel{PS}{\sim} 213$ in
                    Figure~\ref{fig:PSequivalenceExamples}(a).
                
                    \item One can similarly check that $6 \underline{4
                    5 1} 8 7 2 3 6 \stackrel{PS}{\sim} 6 \underline{4 1
                    5} 8 7 2 3$.
                \end{enumerate}
            \end{Example}

            Note, in particular, that the actual values of the
            elements interchanged in
            Example~\ref{eg:PSequivalenceExample} are immaterial so
            long as they have the same relative magnitudes as the
            literal values in $231 \in \mathfrak{S}_{3}$.  (I.e., they
            have to be order-isomorphic to the interchange $231
            \leadsto 213$ as in $451 \leadsto 415$ from
            Example~\ref{eg:PSequivalenceExample}(2).)  Moreover, it
            should also be clear that any value greater than the
            element playing the role of ``1'' can be inserted between
            the elements playing the roles of ``2'' and ``3'' in
            ``231'' without affecting the ability to interchange the
            ``1'' and ``3'' elements.  Problems with this interchange
            only start to arise when a value smaller than the element
            playing the role of ``1'' is inserted between the elements
            playing the roles of ``2'' and ``3''.  More formally, one
            describes this idea using the language of generalized
            permutation patterns
            (Definition~\ref{defn:GeneralizedPattern}):

            \begin{Theorem}
                \label{thm:PSequivalence}
               Let $\sigma, \tau \in \mathfrak{S}_{n}$.  Then $\sigma
               \stackrel{PS}{\sim} \tau$ if and only if $\sigma$ and
               $\tau$ can be transformed into the same permutation by
               changing one or more occurrences of the pattern
               $2\mathrm{-}31$ into occurrences of the pattern
               $2\mathrm{-}13$ such that none of these $2\mathrm{-}31$
               patterns are contained within an occurrence of a
               $3\mathrm{-}1\mathrm{-}42$ pattern.

               In other words, $\stackrel{PS}{\sim}$ is the
               equivalence relation generated by changing
               $3\mathrm{-}\overline{1}\mathrm{-}42$ patterns into
               $3\mathrm{-}\overline{1}\mathrm{-}24$ patterns.
            \end{Theorem}

            \begin{proof}
                By Lemma~\ref{lem:ReversePatienceWordEquivalence}, it
                suffices to show that $\sigma$ can be transformed into
                the reverse patience word $RPW(R(\sigma))$ via a
                sequence of pattern interchanges
                \[
                    \sigma = \sigma^{(0)} \leadsto \sigma^{(1)} \leadsto
                    \sigma^{(2)} \leadsto \cdots \leadsto \sigma^{(\ell)} =
                    RPW(R(\sigma)),
                \]
                where each ``$\leadsto$'' denotes a pattern
                interchange and each $\sigma^{(i)} \stackrel{PS}{\sim}
                \sigma^{(i+1)}$.  However, this should be clear by the
                interpretation of pile configurations in terms
                shadowlines as given by
                Lemma~\ref{lem:ShadowDiagramPileCorrespondence}.
            \end{proof}

            \begin{Example}
            \label{eg:PSequivalenceNonExample}
                ~
                \begin{enumerate}
                   \item Notice that $2431$ contains exactly one instance of
                   a $2\mathrm{-}31$ pattern as the bold underlined
                   subsequence
                   $\mathbf{\underline{2}}4\mathbf{\underline{3}}\mathbf{\underline{1}}$.
                   (Conversely,
                   $\mathbf{\underline{2}}\mathbf{\underline{4}}3\mathbf{\underline{1}}$
                   is an instance of $23\mathrm{-}1$ but not of
                   $2\mathrm{-}31$.)  Moreover, it is clear that $2431
                   \stackrel{PS}{\sim} 2413$.

                   \item Even though $3142$ contains a $2\mathrm{-}31$
                   pattern (as the subsequence
                   $\mathbf{\underline{3}}1\mathbf{\underline{4}}\mathbf{\underline{2}}$),
                   we cannot interchange ``4'' and ``2'', and so
                   $R(3142) \neq R(3124)$.  As illustrated in
                   Figure~\ref{fig:PSequivalenceExamples}(b), this is
                   because ``4'' and ``2'' are on the same shadowline.
                \end{enumerate}
            \end{Example}

            \begin{Remark}
                It follows from Theorem~\ref{thm:PSequivalence} that
                Examples~\ref{eg:PSequivalenceExample}(1) and
                \ref{eg:PSequivalenceNonExample}(2) sufficiently
                characterize $\stackrel{PS}{\sim}$.  It is worth
                pointing out, though, that these examples also begin
                to illustrate an infinite sequence of generalized
                permutation patterns (all of them containing either
                $2\mathrm{-}13$ or $2\mathrm{-}31$) with the following
                property: an interchange of the pattern
                $2\mathrm{-}13$ with the pattern $2\mathrm{-}31$ is
                allowed within an odd-length pattern in this sequence
                unless the elements used to form the odd-length
                pattern can also be used as part of a longer
                even-length pattern in this sequence.
            \end{Remark}

            \begin{Example}
                Even though the permutation $34152$ contains a
                $3\mathrm{-}1\mathrm{-}42$ pattern in the suffix
                ``4152", one can still directly interchange the ``5"
                and the ``2" because of the ``3" prefix (or via the
                sequence of interchanges $34152 \leadsto 31452
                \leadsto 31425 \leadsto 34125$).
            \end{Example}
        
        \section[Enumerating $S_n(3\mathrm{-}\bar{1}\mathrm{-}42)$ and Related Avoidance Sets]{Enumerating $S_n(3\mathrm{-}\bar{1}\mathrm{-}42)$ and\\ Related Avoidance Sets}
        \label{sec:PSasAlgorithm:Enumerating3142}

            In this section, we use results from
            Sections~\ref{sec:PSasAlgorithm:PilesAndNEshadows} and
            \ref{sec:PSasAlgorithm:PSEquivalence} to enumerate and
            characterize the permutations that avoid the generalized
            permutation pattern $2\mathrm{-}31$ unless it occurs as
            part of an occurrence of the generalized pattern
            $3\mathrm{-}1\mathrm{-}42$.  As in
            Section~\ref{sec:Intro:Motivation:Patterns}, we call this
            restricted form of the generalized pattern $2\mathrm{-}31$
            a \emph{barred (generalized) permutation pattern} and
            denote it by $3\mathrm{-}\bar{1}\mathrm{-}42$.

            \begin{Theorem}
            \label{thm:EnumeratingSn3142}
                ~
                \begin{enumerate}
                    \item $S_n(3\mathrm{-}\bar{1}\mathrm{-}42) = \{
                    RPW(R(\sigma)) \mid \sigma \in \mathfrak{S}_{n} \}$.
                    
                    In particular, given $\sigma \in
                    S_n(3\mathrm{-}\bar{1}\mathrm{-}42)$, the entries
                    in each column of $R(\sigma)$ (when read from
                    bottom to top) occupy successive positions in the
                    permutation $\sigma$.
                    
                    \item $S_n(\bar{2}\mathrm{-}41\mathrm{-}3) = \{
                    RPW(R(\sigma))^{-1} \mid \sigma \in \mathfrak{S}_{n}
                    \}$.
                    
                    In particular, given $\sigma \in
                    S_n(\bar{2}\mathrm{-}41\mathrm{-}3)$, the columns
                    of $R(\sigma)$ (when read from top to bottom)
                    contain successive values.
                    
                    \item The size of
                    $S_{n}(3\mathrm{-}\bar{1}\mathrm{-}42)$ is given by
                    the $n^{\mathrm{th}}$ Bell number ${\displaystyle B_{n} = 
                    \frac{1}{e}\sum_{k=0}^{\infty}\frac{k^{n}}{k!}}$.
                \end{enumerate}
           \end{Theorem}

           \begin{proof}
               ~
               \begin{enumerate}
                   \item Let $\sigma \in
                   S_{n}(3\mathrm{-}\bar{1}\mathrm{-}42)$.  Then, for
                   each $i = 1, 2, \ldots, n-1$, define
                   \[
                       \sigma_{m_{i}}
                       =
                       \min_{i \, \leq \, j \, \leq \, n}\{ \sigma_{j} \}.
                   \]
                   Since $\sigma$ avoids
                   $3\mathrm{-}\bar{1}\mathrm{-}42$, the
                   subpermutation
                   $\sigma_{i}\sigma_{i+1}\cdots\sigma_{m_{i}}$ is a
                   decreasing subsequence of $\sigma$.  (Otherwise,
                   $\sigma$ would necessarily contain an occurrence of
                   a $2\mathrm{-}31$ pattern that is not part of an
                   occurrence of a $3\mathrm{-}1\mathrm{-}42$ pattern.)
                   It follows that the left-to-right minima
                   subsequences $s_{1}, s_{2}, \ldots, s_{k}$ of
                   $\sigma$ must be disjoint and satisfy
                   Equation~\eqref{eqn:PileConfigurationCondition}.
                   The result then follows by
                   Lemmas~\ref{lem:PileConfigurationAreSetPartitions}
                   and \ref{lem:ShadowDiagramPileCorrespondence}.
                   
                   \item This follows immediate by taking inverses in
                   Part (1).

                   \item Recall that the Bell number $B_{n}$
                   enumerates all set partitions $[n] = \{1, 2,
                   \ldots, n\}$.  (See \cite{refStanley1997}.)  From
                   Part (1), the elements of
                   $S_{n}(3\mathrm{-}\bar{1}\mathrm{-}42)$ are in
                   bijection with pile configurations.  Thus, since
                   pile configurations are themselves set partitions
                   by
                   Lemma~\ref{lem:PileConfigurationAreSetPartitions},
                   we need only show that every set partition is also
                   a pile configuration.  But this follows by ordering
                   the components of a given set partition by their
                   smallest element so that
                   Equation~\eqref{eqn:PileConfigurationCondition} is
                   satisfied.
               \end{enumerate}
           ~\\[-55pt] 
           \end{proof}

           Even though the set $S_{n}(3\mathrm{-}\bar{1}\mathrm{-}42)$
           is enumerated by such a well-known sequence as the Bell
           numbers, it cannot be described in a simpler way using
           classical pattern avoidance.  This means that there is no
           countable set of non-generalized (a.k.a.~classical)
           permutation patterns $\pi_{1}, \pi_{2}, \ldots$ such that
           \[
               S_{n}(3\mathrm{-}\bar{1}\mathrm{-}42) = S_{n}(\pi_{1},
               \pi_{2}, \ldots) = \bigcap_{i \geq 1} S_{n}(\pi_{i}).
           \]
           There are two very important reasons that this
           cannot happen.

           First of all, the Bell numbers satisfy $\log B_{n} = n
           (\log n - \log\log n + O(1))$ and so exhibit
           superexponential growth.  However, in light of the
           Stanley-Wilf ex-Conjecture (proven by Marcus and Tardos
           \cite{refMT2004} in 2004), the set of permutations
           $S_{n}(\pi)$ avoiding any classical pattern $\pi$ can only
           grow at most exponentially in $n$.

           Second, the so-called \emph{avoidance class} with basis 
           $\{3\mathrm{-}\bar{1}\mathrm{-}42\}$,
           \[
                   Av(3\mathrm{-}\bar{1}\mathrm{-}42) = \bigcup_{n
                   \geq 1} S_{n}(3\mathrm{-}\bar{1}\mathrm{-}42),
           \]
           is not closed under taking order-isomorphic
           subpermutations, whereas it is easy to see that classes of
           permutations defined by classical pattern avoidance must be
           closed.  (See Chapter 5 of \cite{refBona2004}.)  In
           particular, $3142 \in Av(3\mathrm{-}\bar{1}\mathrm{-}42)$,
           but $231 \notin Av(3\mathrm{-}\bar{1}\mathrm{-}42)$.
           
           At the same time, Theorem~\ref{thm:EnumeratingSn3142}
           implies that $3\mathrm{-}\bar{1}\mathrm{-}42$ belongs to
           the so-called \emph{Wilf Equivalence class} for the
           generalized pattern $1\mathrm{-}23$.  That is, if
           \[
               \pi
               \in
               \{ 1\mathrm{-}23, 3\mathrm{-}21, 12\mathrm{-}3, 32\mathrm{-}1,
                  1\mathrm{-}32, 3\mathrm{-}12, 21\mathrm{-}3, 23\mathrm{-}1
               \},
           \]
           then $|S_{n}(\pi)| = B_{n}$.  In particular, Claesson
           \cite{refClaesson2001} showed that $|S_{n}(23\mathrm{-}1)|
           = B_{n}$ using a direct bijection between permutations
           avoiding $23\mathrm{-}1$ and set partitions.  Furthermore,
           given $\sigma \in S_{n}(3\mathrm{-}\bar{1}\mathrm{-}42)$,
           each segment between consecutive \emph{right-to-left
           minima} must form a decreasing subsequence (when read from
           left to right), so it is easy to see that
           $S_{n}(3\mathrm{-}\bar{1}\mathrm{-}42) =
           S_{n}(23\mathrm{-}1)$.  Thus, the barred pattern
           $3\mathrm{-}\bar{1}\mathrm{-}42$ and the generalized
           pattern $23\mathrm{-}1$ are not just in the same Wilf
           equivalence class.  They also have identical avoidance
           classes.
           
           We collect this and similar results together in the
           following Theorem.

           \begin{Theorem}
           \label{thm:23-1and3-1-42Equivalence}
               Let $B_{n}$ denote the $n^{\rm th}$ Bell number.  Then
               \begin{enumerate}
                   \item $S_n(3\mathrm{-}\bar{1}\mathrm{-}42) =
                   S_n(3\mathrm{-}\bar{1}\mathrm{-}4\mathrm{-}2) =
                   S_n(23\mathrm{-}1)$.
                   
                   \item $S_n(31\mathrm{-}\bar{4}\mathrm{-}2) =
                   S_n(3\mathrm{-}1\mathrm{-}\bar{4}\mathrm{-}2) =
                   S_n(3\mathrm{-}12)$.

                   \item $S_n(\bar{2}\mathrm{-}41\mathrm{-}3) =
                   S_n(\bar{2}\mathrm{-}4\mathrm{-}1\mathrm{-}3) =
                   S_n(2\mathrm{-}4\mathrm{-}1\mathrm{-}\bar{3}) =
                   S_n(2\mathrm{-}41\mathrm{-}\bar{3})$.

                   \item $|S_n(\bar{2}\mathrm{-}41\mathrm{-}3)| =
                   |S_n(31\mathrm{-}\bar{4}\mathrm{-}2)| =
                   |S_n(3\mathrm{-}\bar{1}\mathrm{-}42)| = B_n$.
               \end{enumerate}
            \end{Theorem}

            \begin{proof}
                ~
                \begin{enumerate}
                    \item This is proven above.
                
                    \item This follows from Part (1) by taking the
                    reverse complement (as defined in
                    \cite{refBona2004}) of each element in
                    $S_{n}(3\mathrm{-}\bar{1}\mathrm{-}42)$.
                
                    \item The proof is similar to that in Part (2).
                    (This part is also proven in \cite{refALR2005}.)
                
                    \item This follows from the fact that the patterns
                    $3\mathrm{-}1\mathrm{-}\bar{4}\mathrm{-}2$ and
                    $\bar{2}\mathrm{-}4\mathrm{-}1\mathrm{-}3$ are
                    inverses of each other (as classical permutation
                    patterns).
                \end{enumerate}
           ~\\[-55pt] 
            \end{proof}

           \begin{Remark}
               Even though $S_{n}(3\mathrm{-}\bar{1}\mathrm{-}42) =
               S_{n}(23\mathrm{-}1)$, it is more natural to use
               avoidance of $3\mathrm{-}\bar{1}\mathrm{-}42$ when
               studying Patience Sorting.  Fundamentally, this lets us
               look at $S_{n}(3\mathrm{-}\bar{1}\mathrm{-}42)$ as the
               set of equivalence classes in $\mathfrak{S}_{n}$ modulo
               $3\mathrm{-}\bar{1}\mathrm{-}42 \stackrel{PS}{\sim}
               3\mathrm{-}\bar{1}\mathrm{-}24$, where each equivalence
               class corresponds to a unique pile configuration.  The
               same equivalence relation is not easy to describe when
               starting with an occurrence of $23\mathrm{-}1$.
              
               In other words, both
               \[
                   23\mathrm{-}1 \sim 2\mathrm{-}13
                   \mbox{ \ and \ }
                   23\mathrm{-}1 \sim 21\mathrm{-}3
               \]
               are wrong since they incorrectly suggest
               \[
                   2431 \sim 2314
                   \mbox{ \ and \ }
                   2431 \sim 2134,
               \]
               respectively, instead of the correct $2431
               \stackrel{PS}{\sim} 2413$.
           \end{Remark}
           
            We conclude this section with an immediate corollary of
            Theorem~\ref{thm:23-1and3-1-42Equivalence} that
            characterizes an important category of classical
            permutation patterns.

            \begin{Definition}
                Given a composition $\gamma = (\gamma_{1}, \gamma_{2},
                \ldots, \gamma_{m}) \models n$, the
                \emph{(classical) layered permutation pattern}
                $\pi_{\gamma} \in \mathfrak{S}_{n}$ is the permutation
                \[
                    \gamma_{1}(\gamma_{1} - 1) \cdots 1
                    (\gamma_{1} + \gamma_{2}) (\gamma_{1} + \gamma_{2} - 1) \cdots (\gamma_{1} + 1)
                    \; \cdots \;
                    n (n - 1) \cdots (\gamma_{1} + \gamma_{2} + \cdots + \gamma_{m - 1} + 1).
                \]
            \end{Definition}
            
            \begin{Example}
                Given $\gamma = (3, 2, 3) \models 8$, the
                corresponding layered pattern (following the notation
                in \cite{refPrice1997}) is $\pi_{(3,2,3)} =
                \widehat{321} \widehat{54} \widehat{876} \in
                \mathfrak{S}_{8}$.
            \end{Example}

            \begin{Corollary}
                \label{cor:LayeredPatternCharacterization}
                $S_{n}(3\mathrm{-}\bar{1}\mathrm{-}42,
                \bar{2}\mathrm{-}41\mathrm{-}3)$ is the set of layered
                patterns in $\mathfrak{S}_{n}$.
            \end{Corollary}

            \begin{proof}
                By Theorem~\ref{thm:23-1and3-1-42Equivalence}, we have that
                $S_{n}(3\mathrm{-}\bar{1}\mathrm{-}42,
                \bar{2}\mathrm{-}41\mathrm{-}3) = S_{n}(23\mathrm{-}1,
                31\mathrm{-}2)$, the latter being a characterization
                for layered patterns given in \cite{refCM2002}.
            \end{proof}
        
        \section[Invertibility of Patience Sorting]{Invertibility of Patience Sorting}
        \label{sec:PSasAlgorithm:Invertibility}

            As discussed in
            Section~\ref{sec:PSasAlgorithm:PSEquivalence}, many
            different permutations can correspond to the same pile
            configuration under Patience Sorting
            (Algorithm~\ref{alg:MallowsPSprocedure}).  E.g., $R(3142)
            = R(3412)$.  In this section, we use barred permutation
            patterns to characterize permutations for which this does
            not hold.  We then establish a non-trivial enumeration for
            the resulting avoidance sets.

            \begin{Theorem}
            \label{thm:PermutationsWithUniquePileConfigurations}
                A pile configuration pile $R$ has a unique preimage
                $\sigma \in \mathfrak{S}_{n}$ under Patience Sorting
                if and only if $\sigma \in
                S_{n}(3\mathrm{-}\bar{1}\mathrm{-}42,
                3\mathrm{-}\bar{1}\mathrm{-}24)$.
            \end{Theorem}

            \begin{proof}
                From the proof of
                Lemma~\ref{lem:PileConfigurationAreSetPartitions}, we
                know that every pile configuration $R$ has at least
                its reverse patience word $RPW(R)$ as a preimage under
                Patience Sorting, and, by
                Theorem~\ref{thm:EnumeratingSn3142}, $RPW(R) \in
                S_{n}(3\mathrm{-}\bar{1}\mathrm{-}42)$.  Furthermore,
                by Theorem~\ref{thm:PSequivalence}, two permutations
                yield the same pile configurations under Patience
                Sorting if and only if one can be obtained from the
                other by a sequence of order-isomorphic exchanges of
                the form
                \[
                    3\mathrm{-}\bar{1}\mathrm{-}24 \leadsto 3\mathrm{-}\bar{1}\mathrm{-}42
                    \mbox{ \ or \ }
                    3\mathrm{-}\bar{1}\mathrm{-}42\leadsto 3\mathrm{-}\bar{1}\mathrm{-}24.
                \]
                (I.e., the occurrence of one pattern is reordered to
                form an occurrence of the other pattern.)  Thus, it is
                easy to see that $R$ has the unique preimage $RPW(R)$
                if and only if $RPW(R)$ avoids both
                $3\mathrm{-}\bar{1}\mathrm{-}42$ and
                $3\mathrm{-}\bar{1}\mathrm{-}24$.
            \end{proof}
            
            Given this pattern avoidance characterization of
            invertibility of Patience Sorting, we have the following
            recurrence relation for the size of the avoidance sets in
            Theorem~\ref{thm:PermutationsWithUniquePileConfigurations}.

            \begin{Lemma}
            \label{lem:RecurrenceRelation}
                Set $f(n) = |S_n(3\mathrm{-}\bar{1}\mathrm{-}42,
                3\mathrm{-}\bar{1}\mathrm{-}24)|$ and, for $k \le n$,
                denote by $f(n, k)$ the cardinality
                \[
                    f(n, k)
                    =
                    \#
                    \{
                    \sigma \in S_{n}(3\mathrm{-}\bar{1}\mathrm{-}42, 3\mathrm{-}\bar{1}\mathrm{-}24)
                    \mid
                    \sigma(1) = k
                    \}.
                \]
                Then ${\displaystyle f(n) = \sum_{k = 1}^{n}{f(n,
                k)}}$, and $f(n, k)$ satisfies the four part recurrence
                relation
                \begin{eqnarray}
                    f(n, 0) = 0
                    & \textrm{for}
                    & n \geq 1
                    \label{eqn:RecurrenceRelation:Eqn1}\\
                    f(n, 1) = f(n, n) = f(n - 1)
                    & \textrm{for}
                    & n \geq 1
                    \label{eqn:RecurrenceRelation:Eqn2}\\
                    f(n, 2) = 0
                    & \textrm{for}
                    & n \geq 3
                    \label{eqn:RecurrenceRelation:Eqn3}\\
                    f(n, k) = f(n, k - 1) + f(n - 1, k - 1) + f(n - 2, k - 2)
                    & \textrm{for}
                    & n \geq 3
                    \label{eqn:RecurrenceRelation:Eqn4}
                \end{eqnarray}
                subject to the initial conditions $f(0, 0) = f(0) = 1$.

            \end{Lemma}

            \begin{proof}
                Note first that Equation
                \eqref{eqn:RecurrenceRelation:Eqn1} is the obvious
                boundary condition for $k = 0$.

                Now, suppose that the first component of $\sigma \in
                S_{n}(3\mathrm{-}\bar{1}\mathrm{-}42,
                3\mathrm{-}\bar{1}\mathrm{-}24)$ is either $\sigma(1)
                = 1$ or $\sigma(1) = n$.  Then $\sigma(1)$ cannot be
                part of any occurrence of
                $3\mathrm{-}\bar{1}\mathrm{-}42$ or
                $3\mathrm{-}\bar{1}\mathrm{-}24$ in $\sigma$.  Thus,
                upon removing $\sigma(1)$ from $\sigma$, and
                subtracting one from each component if $\sigma(1) =
                1$, a bijection is formed with
                $S_{n-1}(3\mathrm{-}\bar{1}\mathrm{-}42,
                3\mathrm{-}\bar{1}\mathrm{-}24)$.  Therefore,
                Equation~\eqref{eqn:RecurrenceRelation:Eqn2} follows.

                Next, suppose that the first component of $\sigma \in
                S_{n}(3\mathrm{-}\bar{1}\mathrm{-}42,
                3\mathrm{-}\bar{1}\mathrm{-}24)$ is $\sigma(1) = 2$.
                Then the first column of $R(\sigma)$ must be $r_{1} =
                21$ regardless of where 1 occurs in $\sigma$.
                Therefore, $R(\sigma)$ has the unique preimage
                $\sigma$ if and only if $\sigma = 21 \in
                \mathfrak{S}_{2}$, and from this Equation
                \eqref{eqn:RecurrenceRelation:Eqn3} follows.

                Finally, suppose that $\sigma \in
                S_{n}(3\mathrm{-}\bar{1}\mathrm{-}42,
                3\mathrm{-}\bar{1}\mathrm{-}24)$ with $3\le k\le n$.
                Since $\sigma$ avoids
                $3\mathrm{-}\bar{1}\mathrm{-}42$, $\sigma$ is a RPW by
                Theorem~\ref{thm:EnumeratingSn3142}, and hence the
                left prefix of $\sigma$ from $k$ to $1$ is a
                decreasing subsequence.  Let $\sigma'$ be the
                permutation obtained by interchanging the values $k$
                and $k-1$ in $\sigma$.  Then the only instances of the
                patterns $3\mathrm{-}\bar{1}\mathrm{-}42$ and
                $3\mathrm{-}\bar{1}\mathrm{-}24$ in $\sigma'$ must
                involve both $k$ and $k-1$.  Note that the number of
                $\sigma$ for which no instances of these patterns are
                created by interchanging $k$ and $k - 1$ is $f(n, k -
                1)$.

                There are now two cases in which an instance of the
                barred pattern $3\mathrm{-}\bar{1}\mathrm{-}42$ or
                $3\mathrm{-}\bar{1}\mathrm{-}24$ will be created in
                $\sigma'$ by this interchange:

                \emph{Case 1.} If $k - 1$ occurs between $\sigma(1) =
                k$ and $1$ in $\sigma$, then $\sigma(2) = k - 1$, so
                interchanging $k$ and $k-1$ will create an instance of
                the pattern $23\mathrm{-}1$ via the subsequence
                $(k-1,k,1)$ in $\sigma'$.  Thus, by
                Theorem~\ref{thm:23-1and3-1-42Equivalence}, $\sigma'$
                contains $3\mathrm{-}\bar{1}\mathrm{-}42$ from which
                $\sigma' \in S_{n}(3\mathrm{-}\bar{1}\mathrm{-}42)$ if
                and only if $k-1$ occurs after $1$ in $\sigma$.  Note
                also that if $\sigma(2) = k - 1$, then removing $k$
                from $\sigma$ yields a bijection with permutations in
                $S_{n-1}(3\mathrm{-}\bar{1}\mathrm{-}42,
                3\mathrm{-}\bar{1}\mathrm{-}24)$ that start with
                $k-1$.  Therefore, the number of permutations counted
                in Case 1 is $f(n-1,k-1)$.

                \emph{Case 2.} If $k-1$ occurs to the right of $1$ in
                $\sigma$, then $\sigma'$ both contains the subsequence
                $(k-1,1,k)$ and avoids the pattern
                $3\mathrm{-}\bar{1}\mathrm{-}42$, so it must also
                contain the pattern $3\mathrm{-}\bar{1}\mathrm{-}24$.
                If an instance of $3\mathrm{-}\bar{1}\mathrm{-}24$ in
                $\sigma'$ involves both $k-1$ and $k$, then $k-1$ and
                $k$ must play the roles of ``3'' and ``4'',
                respectively.  Moreover, if the value $\ell$ preceding
                $k$ is not $1$, then the subsequence $(k-1,1,\ell,k)$
                is an instance of $3\mathrm{-}1\mathrm{-}24$, so
                $(k-1,\ell,k)$ is not an instance of
                $3\mathrm{-}\bar{1}\mathrm{-}24$.  Therefore, for
                $\sigma'$ to contain $3\mathrm{-}\bar{1}\mathrm{-}24$,
                $k$ must follow $1$ in $\sigma'$, and so $k-1$ follows
                $1$ in $\sigma$.  Similarly, if the letter preceding
                $1$ is some $m < k$, then the subsequence $(m,1,k-1)$
                is an instance of $3\mathrm{-}\bar{1}\mathrm{-}24$ in
                $\sigma$, which is impossible.  Therefore, $k$ must
                precede $1$ in $\sigma$, from which $\sigma$ must
                start with the initial segment $(k,1,k-1)$.  It
                follows that removing the values $k$ and $1$ from
                $\sigma$ and then subtracting $1$ from each component
                yields a bijection with permutations in
                $S_{n-2}(3\mathrm{-}\bar{1}\mathrm{-}42,
                3\mathrm{-}\bar{1}\mathrm{-}24)$ that start with
                $k-2$.  Thus, the number of permutations counted in
                Case 2 is then exactly $f(n - 2, k - 2)$, which yields
                Equation \eqref{eqn:RecurrenceRelation:Eqn4}.
            \end{proof}
             
            If we denote by
            \[
                \Phi(x, y) = \sum_{n = 0}^{\infty} \sum_{k=0}^{n}
                f(n,k)x^{n}y^{k}
            \]
            the bivariate generating
            function for the sequence $\{f(n,k)\}_{n \geq k\geq 0}$, then
            Equation \eqref{eqn:RecurrenceRelation:Eqn4} implies that
            \[
                (1-y-xy-x^{2}y^{2})\Phi(x,y) =
                1-y-xy+xy^{2}-xy^{2}\Phi(xy,1)+xy(1-y-xy)\Phi(x,1).
            \]

            We conclude this section with the following enumerative result.

            \begin{Theorem}
            \label{thm:EnumerateInvertibility}
                Denote by $F_{n}$ the $n^{\rm th}$ Fibonacci number (with
                $F_{0} = F_{1} = 1$) and by
                \[
                    a(n,k) \hspace{0.25cm} =
                    \hspace{-2cm}
                    \sum_{\stackrel{n_1,\dots,n_{k+1}\ge
                    0}{n_1+\dots+n_{k+1}=(n-2)-(k+1)=n-k-3}}
                    \hspace{-2cm}
                    {F_{n_1}F_{n_2}\dots
                    F_{n_{k+1}}}
                \]
                the convolved Fibonacci numbers for $n \geq k + 2$
                (with $a(n,k) := 0$ otherwise).  Then, defining
                \[
                    X =
                        \begin{bmatrix}
                            f(0)  \\
                            f(1)  \\
                            f(2)  \\
                            f(3)  \\
                            f(4)  \\
                            \vdots
                        \end{bmatrix},
                    \quad
                    F =
                        \begin{bmatrix}
                            1  \\
                            F_{0} \\
                            F_{1} \\
                            F_{2} \\
                            F_{3} \\
                            \vdots
                        \end{bmatrix},
                    \ \ \mathrm{and} \quad
                    \mathbf{A} = (a(n,k))_{n,k\ge 0},
                \]
                we have that $X = (\mathbf{I} - \mathbf{A})^{-1}F$,
                where $\mathbf{I}$ is the infinite identity matrix and
                $\mathbf{A}$ is, by definition, lower triangular with
                zero main diagonal.
            \end{Theorem}

            \begin{proof}
                From Equations
                \eqref{eqn:RecurrenceRelation:Eqn1}--\eqref{eqn:RecurrenceRelation:Eqn4},
                one can conjecture an equivalent recurrence in which
                Equations~\eqref{eqn:RecurrenceRelation:Eqn3} and
                \eqref{eqn:RecurrenceRelation:Eqn4} are replaced by
                the following equation (and where $\delta_{nk}$
                denotes the Kronecker delta function):
                \begin{equation}
                \label{eqn:RecurrenceRelation:Eqn5}
                    f(n,k)=\sum_{m=0}^{k-3}{c(k,m)f(n-k+m)} +
                    \delta_{nk}F_{k-2}, \quad n\ge k\ge 2.
                \end{equation}
                For this relation to hold, the coefficients
                $c(k,m)$ must satisfy the recurrence relation
                \[
                    c(k,m) = c(k-1,m-1) + c(k-1,m) + c(k-2,m), \quad k\ge
                    2,
                \]
                or, equivalently,
                \[
                    c(k-1,m-1) = c(k,m) - c(k-1,m) - c(k-2,m), \quad k\ge
                    2,
                \]
                with $c(2,0)=1$ and $c(k,m)=0$ in the case that
                $k < 2$, $m < 0$, or $m > k - 2$.  This implies that the
                generating function for the sequence $\{c(k,m)\}_{k\ge 0}$
                (for each $m\ge 0$) is
                \[
                \sum_{n\ge 0}{c(k,m)x^k} = \frac{x^{m+2}}{(1-x-x^2)^{m+1}}.
                \]
                Thus, the coefficients $c(k,m) = a(k,m)$ in Equation
                \eqref{eqn:RecurrenceRelation:Eqn5} are the convolved
                Fibonacci numbers \cite{refOEIS} and form the
                so-called skew Fibonacci-Pascal triangle in the matrix
                $\mathbf{A}=(a(k,m))_{k,m\ge 0}$.  In particular, the
                sequence of nonzero entries in column $m\ge 0$ of
                $\mathbf{A}$ is the $m^{\rm th}$ convolution of the
                sequence $\{F_n\}_{n\ge 0}$.

                Finally, upon combining the expansion of $f(n,n)$ from
                Equation \eqref{eqn:RecurrenceRelation:Eqn5} with
                Equation \eqref{eqn:RecurrenceRelation:Eqn2},
                \[
                    f(n) = \sum_{m=0}^{n-2}{a(n,m)f(m)} + F_{n-1},
                \]
                which is equivalent to the matrix equation $X =
                \mathbf{A}X + F$.  Therefore, since $\mathbf{I} -
                \mathbf{A}$ is clearly invertible, the result follows.
            \end{proof}
            
            \begin{Remark}
                Since $\mathbf{A}$ is strictly lower triangular with
                zero main diagonal and zero sub-diagonal, it follows
                that multiplication of a matrix $\mathbf{B}$ by
                $\mathbf{A}$ shifts the position of the highest
                nonzero diagonal in $\mathbf{B}$ down by two rows.
                Thus, $(\mathbf{I} - \mathbf{A})^{-1}=\sum_{n\ge
                0}{\mathbf{A}^n}$ as a Neumann series, and all nonzero
                entries of $(\mathbf{I} - \mathbf{A})^{-1}$ are
                positive integers.
            \end{Remark}

            In particular, one can explicitly compute
            
            {\singlespacing
            \[
                \mathbf{A} =
                \begin{bmatrix}
                    0 \\
                    0 &  0 \\
                    1 &  0 &  0 \\
                    1 &  1 &  0 &  0 \\
                    2 &  2 &  1 &  0 &  0 \\
                    3 &  5 &  3 &  1 &  0 &  0 \\
                    5 & 10 &  9 &  4 &  1 &  0 &  0 \\
                    8 & 20 & 22 & 14 &  5 &  1 &  0 &  0 \\
                    \vdots & \vdots & \vdots & \vdots & \vdots &
                    \vdots & \vdots & \vdots & \ddots \\
                \end{bmatrix}
            \]
            }
            
            \noindent from which
            
            {\singlespacing
            \[
                (\mathbf{I} - \mathbf{A})^{-1} =
                \begin{bmatrix}
                     1 \\
                     0 &  1 \\
                     1 &  0 &  1 \\
                     1 &  1 &  0 &  1 \\
                     3 &  2 &  1 &  0 &  1 \\
                     7 &  6 &  3 &  1 &  0 &  1 \\
                    21 & 16 & 10 &  4 &  1 &  0 &  1 \\
                    66 & 50 & 30 & 15 &  5 &  1 &  0 & 1 \\
                    \vdots &  \vdots & \vdots &  \vdots &  \vdots & \vdots & \vdots & \vdots & \ddots\\
                \end{bmatrix}.
            \]
            }
            
            \bigskip
            
            \noindent It follows that the first few values of the sequence
            $\{f(n)\}_{n \ge 0}$ are
            \[
                1, 1, 2, 4, 9, 23, 66, 209, 718,
                2645, 10373, 43090, 188803, 869191, 4189511, \ldots.
            \]
        
    %
    %

    \newchapter{Bijectively Extending Patience Sorting}{Bijectively Extending Patience Sorting}{Bijectively Extending Patience Sorting}
    \label{sec:ExtendingPS}
    
        \section[Characterizing ``Stable Pairs'']{Characterizing ``Stable Pairs''}
        \label{sec:ExtendingPS:StablePairs}

            According to Theorem~\ref{thm:EnumeratingSn3142}, the
            number of pile configurations that can be formed from
            $[n]$ is given by the Bell number $B_{n}$.  Comparing this
            to the number of standard Young tableau
            $|\mathfrak{T}_{n}|$ (see, e.g., \cite{refKnuth1998}), it
            is clear that there are significantly more possible pile
            configurations than standard Young tableau.  Consequently,
            not every ordered pair of pile configurations with the
            same shape can result from Extended Patience Sorting
            (Algorithm~\ref{alg:ExtendedPSalgorithm}).  In this
            section, we characterize the ``stable pairs'' of pile
            configurations that result from applying Extended Patience
            Sorting to a permutation.
            
            The following example, though very small, illustrates the
            most generic behavior that must be avoided in constructing
            these ``stable pairs''.  As in
            Example~\ref{eg:PileReflectionExample}, we denote by $S'$
            the ``reversed pile configuration'' corresponding to $S$
            (which has all piles listed in reverse order).

            \begin{Example}
            \label{eg:SmallBadPilesExample}
                Even though the pile configuration $R = \{\{3 > 1\},
                \{2\}\}$ cannot result as the insertion piles of an
                involution under Extended Patience Sorting, we can
                nonetheless look at the shadow diagram for the
                pre-image of the pair $(R, R)$ under
                Algorithm~\ref{alg:ExtendedPSalgorithm}:

                {\singlespacing
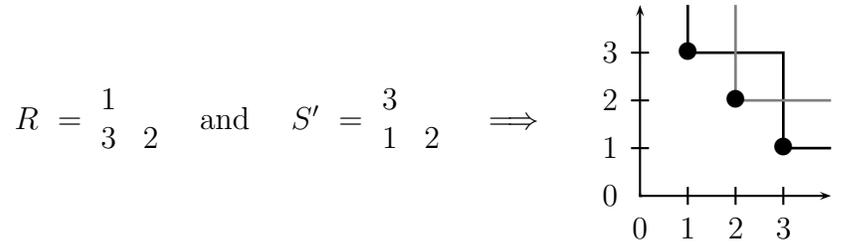
                \begin{center}
                    \begin{tabular}{lr}
                        \raisebox{.35in}{
                                \begin{minipage}[c]{3in}
                                        $R \ = \
                                          \begin{matrix}
                                              1 & \\
                                              3 & 2
                                          \end{matrix}$
                                          \quad \text{and} \quad
                                          $S' \ = \
                                          \begin{matrix}
                                              3 & \\
                                              1 & 2
                                          \end{matrix}
                                          \quad \implies $
                                \end{minipage}
                        } &

                        \psset{xunit=0.25in,yunit=0.25in}

                        \begin{pspicture}(0,0)(4,4)

                            \psaxes{->}(4,4)

                            \psline[linecolor=black,linewidth=1pt](1,4)(1,3)(3,3)(3,1)(4,1)

                            \psline[linecolor=gray,linewidth=1pt](2,4)(2,2)(4,2)

                             \rput(1,3){{\Large $\bullet$}}%
                             \rput(2,2){{\Large $\bullet$}}%
                             \rput(3,1){{\Large $\bullet$}}%

                        \end{pspicture}\\[0.25in]

                    \end{tabular}.
                \end{center}
                }

                \noindent Note that there are two competing
                constructions here.  On the one hand, we have the
                diagram $\{(1, 3), (2, 2), (3, 1)\}$ of the
                permutation $321 \in \mathfrak{S}_{3}$ given by the
                entries in the pile configurations.  (In particular,
                the values in $R$ specify the ordinates, and the
                values in the corresponding boxes of $S'$ specify the
                abscissae.)  On the other hand, the piles in $R$ also
                specify shadowlines with respect to this permutation
                diagram.  Here, the pair $(R, S) = (R, R)$ of pile
                configurations is ``unstable'' because their
                combination yields crossing shadowlines --- which is
                clearly not allowed.
                
                Similar considerations lead to crossings of the form\\

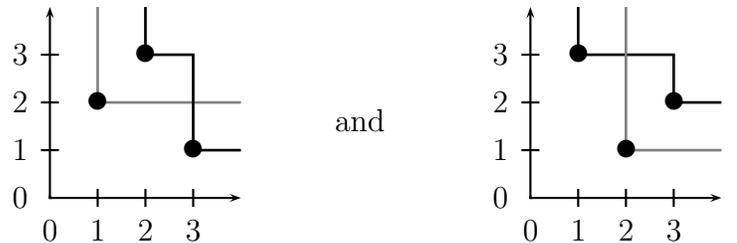
                \begin{center}

                    \begin{tabular}{rlr}

                        \psset{xunit=0.25in,yunit=0.25in}

                        \begin{pspicture}(0,0)(4,4)

                            \psaxes{->}(4,4)

                            \psline[linecolor=gray,linewidth=1pt](1,4)(1,2)(4,2)

                            \psline[linecolor=black,linewidth=1pt](2,4)(2,3)(3,3)(3,1)(4,1)

                             \rput(1,2){{\Large $\bullet$}}%
                             \rput(2,3){{\Large $\bullet$}}%
                             \rput(3,1){{\Large $\bullet$}}%

                        \end{pspicture}

                        &

                        \raisebox{.35in}{\qquad and \qquad\quad}

                        &

                        \psset{xunit=0.25in,yunit=0.25in}

                        \begin{pspicture}(0,0)(4,4)

                            \psaxes{->}(4,4)

                            \psline[linecolor=black,linewidth=1pt](1,4)(1,3)(3,3)(3,2)(4,2)

                            \psline[linecolor=gray,linewidth=1pt](2,4)(2,1)(4,1)

                             \rput(1,3){{\Large $\bullet$}}%
                             \rput(2,1){{\Large $\bullet$}}%
                             \rput(3,2){{\Large $\bullet$}}%

                        \end{pspicture}\\[0.25in]

                    \end{tabular}.

                \end{center}

                \noindent Note also that these latter two crossings can
                also be used together to build something like the first
                crossing but with ``extra'' elements on the boundary of
                the polygon formed:

                \begin{center}

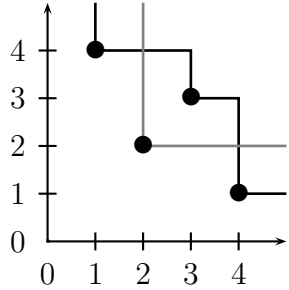
                    \psset{xunit=0.25in,yunit=0.25in}

                    \begin{pspicture}(0,0)(5,5)

                        \psaxes{->}(5,5)

                        \psline[linecolor=black,linewidth=1pt](1,5)(1,4)(3,4)(3,3)(4,3)(4,1)(5,1)

                        \psline[linecolor=gray,linewidth=1pt](2,5)(2,2)(5,2)

                         \rput(1,4){{\Large $\bullet$}}%
                         \rput(2,2){{\Large $\bullet$}}%
                         \rput(3,3){{\Large $\bullet$}}%
                         \rput(4,1){{\Large $\bullet$}}%

                    \end{pspicture}\\[0.4in]

                \end{center}

            \end{Example}

            We are now in a position to make the following fundamental
            definitions:

            \begin{Definition}
                Given a composition $\gamma$ of $n$ (denoted $\gamma
                \models n$), we define $\mathfrak{P}_{n}(\gamma)$ to
                be the set of all pile configurations $R$ having shape
                $\mathrm{sh}(R) = \gamma$ and put
                \begin{displaymath}
                    \mathfrak{P}_{n} = \bigcup_{\gamma
                    \ \models \ n}\mathfrak{P}_{n}(\gamma).
                \end{displaymath}
            \end{Definition}

            \begin{Definition}
            \label{defn:StablePairs}
                Define the set $\Sigma(n) \subset \mathfrak{P}_{n}
                \times \mathfrak{P}_{n}$ to consist of all ordered
                pairs $(R, S)$ with $\mathrm{sh}(R) = \mathrm{sh}(S)$
                such that the ordered pair $(RPW(R), RPW(S'))$ avoids
                simultaneous occurrences of the pairs of patterns
                $(31\mathrm{-}2,13\mathrm{-}2)$,
                $(31\mathrm{-}2,32\mathrm{-}1)$ and
                $(32\mathrm{-}1,13\mathrm{-}2)$ at the same positions
                in $RPW(R)$ and $RPW(S')$.
            \end{Definition}

            In other words, if $RPW(R)$ contains an occurrence of the
            first pattern in any of the above pairs, then $RPW(S')$
            cannot contain an occurrence at the same positions of the
            second pattern in the same pair, and vice versa.  In
            effect, Definition~\ref{defn:StablePairs} characterizes
            ``stable pairs'' of pile configurations $(R, S)$ by
            forcing $R$ and $S$ to avoid certain sub-pile pattern
            pairs.  As in Example~\ref{eg:SmallBadPilesExample}, we
            are characterizing when the induced shadowlines cross.

            \begin{Theorem}
            \label{thm:ExtendedPSbijection}
                Extended Patience Sorting
                (Algorithm~\ref{alg:ExtendedPSalgorithm}) gives a
                bijection between the symmetric group
                $\mathfrak{S}_{n}$ and the ``stable pairs'' set
                $\Sigma(n)$.
            \end{Theorem}

            \begin{proof}
                We show that, for any ``stable pair'' $(R, S) \in
                \Sigma(n)$ and any permutation $\sigma \in
                \mathfrak{S}_{n}$, $(R, S)=(R(\sigma), S(\sigma))$ if
                and only if
                \[
                    \sigma = \begin{pmatrix} RPW(S') \\ RPW(R) \end{pmatrix}
                    \; \text{(in the two-line notation)}.
                \]
                Clearly, if $(R, S)=(R(\sigma), S(\sigma))$ for some
                $\sigma \in \mathfrak{S}_{n}$, then $\sigma =
                \begin{pmatrix} RPW(S') \\ RPW(R) \end{pmatrix}$.
                Thus, we need only to prove that $(R, S)\in\Sigma(n)$.
                In particular, if $(R, S) \notin \Sigma(n)$, then
                $RPW(R)$ and $RPW(S')$ contain simultaneous
                occurrences of (at least) one of the three forbidden
                pattern pairs given in
                Definition~\ref{defn:StablePairs}.

                Suppose that $RPW(R)$ contains an occurrence $(r_{3},
                r_{1}, r_{2})$ of $31\mathrm{-}2$, and suppose also
                that $RPW(S')$ contains an occurrence
                $(s'_1,s'_3,s'_2)$ of $13\mathrm{-}2$, with both
                occurrences at the same positions.  Since $r_3 > r_1$
                and since $r_3$ and $r_1$ are consecutive entries in
                $RPW(R)$, it follows that $r_3$ and $r_1$ must be in
                the same column $c_i(R)$ of $R$ (in fact, $r_1$ is
                immediately on top of $r_3$).  Furthermore, since $r_1
                < r_2$ and since $r_2$ is to the right of $r_1$ in
                $R$, it follows that the column $c_j(R)$ of $R$
                containing $r_2$ must be to the right of the column
                containing $r_1$ atop $r_3$.  Therefore, $s'_2$ must
                also be in a column $c_i(S')$ of $S'$ to the right of
                the column $c_j(S')$ containing $s'_3$ atop $s'_1$.

                Now, consider the subpermutation $\tau$ of $\sigma$
                formed by removing all components of $RPW(R)$ and
                $RPW(S')$ that are not in these two columns.
                Alternatively, let $R_\ast$ and $S'_\ast$ consist only
                of the columns $(c_i(R), c_j(R))$ of $R$ and
                $(c_i(S'),c_j(S'))$ of $S'$, respectively.  Then
                \[
                    \tau
                    =
                    \begin{pmatrix}RPW(S'_\ast)\\RPW(R_\ast)\end{pmatrix}.
                \]
                Note that the values $r_3$ and $r_1$ in $c_i(S')$ are
                consecutive left-to-right minima of $\tau$, whereas
                $r_2$ is not a left-to-right minimum of $\tau$.  Since
                $r_1<r_2<r_3$, it follows that $r_2$ cannot occur
                between $r_1$ and $r_3$ in $\tau$.  However, since
                $\begin{pmatrix}s'_1 & s'_3 & s'_2\\r_3 & r_1 &
                r_2\end{pmatrix}$ is a subpermutation of $\tau$ and
                since $s'_1<s'_2<s'_3$, it follows that $r_2$ does
                occur between $r_1$ and $r_3$, which is a
                contradiction.

                A similar argument applies to both $(31\mathrm{-}2,
                32\mathrm{-}1)$ and $(32\mathrm{-}1, 13\mathrm{-}2)$,
                which then implying $(R, S) \in \Sigma(n)$.

                Conversely, given $(R, S) \in \Sigma(n)$, set $\sigma
                = \begin{pmatrix}RPW(S')\\RPW(R)\end{pmatrix}$.  Then,
                since the pattern avoidance conditions defining
                $\Sigma(n)$ forbid intersections in the northeast
                shadow diagram corresponding to $\sigma$ (as
                illustrated in Example~\ref{eg:SmallBadPilesExample}),
                it follows by
                Lemma~\ref{lem:ShadowDiagramPileCorrespondence} that
                $(R, S) = (R(\sigma), S(\sigma))$.
            \end{proof}

            \begin{Example}
            \label{eg:GoodPilesExample}
                The pair of piles
                
                {\singlespacing
                \[
                    (R,S)=
                    \left(
                    \begin{tabular}{c c c}
                        1 &   & 3 \\
                        4 & 2 & 7 \\
                        6 & 5 & 8
                    \end{tabular}
                    \; , \;
                    \begin{tabular}{c c c}
                        1 &   & 5 \\
                        2 & 3 & 6 \\
                        4 & 7 & 8
                    \end{tabular}
                    \right)
                    \in\Sigma(8)
                \]
                }
                
                \noindent corresponds to the permutation
                \[
                    \sigma=\begin{pmatrix}RPW(S')\\RPW(R)\end{pmatrix}=
                    \begin{pmatrix}
                    1&2&4&3&7&5&6&8\\
                    6&4&1&5&2&8&7&3
                    \end{pmatrix}
                    =
                    64518723\in \mathfrak{S}_8.
                \]
            \end{Example}
            
            The similarities between Extended Patience Sorting
            (Algorithm~\ref{alg:ExtendedPSalgorithm}) and the RSK
            Correspondence (Algorithm~\ref{alg:RSKAlgorithm}) are
            perhaps most observable in the following simple
            Proposition.

            \begin{Proposition}
            \label{prop:PSmonotonePatterns}
                Let $\textrm{\emph{\i}}_k =
                1\mathrm{-}2\mathrm{-}\cdots\mathrm{-}k$ and
                $\textrm{\emph{\j}}_k =
                k\mathrm{-}\cdots\mathrm{-}2\mathrm{-}1$ be the
                classical monotone permutation patterns.  Then there
                is a bijection

                \begin{enumerate}

                    \item between $S_n(\textrm{\emph{\i}}_{k+1})$ and
                    ``stable pairs'' of pile configurations having the same
                    composition shape $\gamma = (\gamma_{1},
                    \gamma_{2}, \ldots, \gamma_{m}) \models n$ but
                    with at most $k$ piles (i.e., $m \leq k$),

                    \item as well as a bijection between
                    $S_n(\textrm{\emph{\j}}_{k+1})$ and ``stable
                    pairs'' of pile configurations having the same
                    composition shape $\gamma = (\gamma_{1},
                    \gamma_{2}, \ldots, \gamma_{m}) \models n$ but
                    with no pile having more than $k$ cards in it
                    (i.e., $\gamma_{i} \leq k$ for each $i = 1, 2,
                    \ldots, m$).

                \end{enumerate}
            \end{Proposition}

            \begin{proof}
                ~
                \begin{enumerate}
                    \item Given $\sigma \in \mathfrak{S}_{n}$, the
                    proof of
                    Proposition~\ref{prop:GreedyStrategyUpperBoundPSviaLIS}
                    yields a bijection between the set of piles
                    $R(\sigma) = \{r_{1}, r_{2}, \ldots, r_{k}\}$
                    formed under Patience Sorting and the components
                    of a particular longest increasing subsequence in
                    $\sigma$.  Since avoiding the monotone pattern
                    $\textrm{\i}_{k+1}$ is equivalent to restricting
                    the length of the longest increasing subsequence
                    in a permutation, the result then follows.

                    \item Follows from Part~(1) by reversing the order
                    of the components in each of the permutations in
                    $S_n(\textrm{\i}_{k+1})$ in order to form
                    $S_n(\textrm{\j}_{k+1})$.\\[-55pt]
                \end{enumerate}
            \end{proof}

            Proposition~\ref{prop:PSmonotonePatterns} states that
            Patience Sorting can be used to efficiently compute the
            length of both the longest increasing and longest
            decreasing subsequences in a given permutation.  In
            particular, one can compute these lengths without
            examining every subsequence of a permutation, just as with
            the RSK Correspondence.  However, while both the RSK
            Correspondence and Patience Sorting can be used to
            implement this computation in $O(n\log(n))$ time, an
            extension is given in \cite{refBS2000} that also
            simultaneously tabulates all of the longest increasing or
            decreasing subsequences without incurring any additional
            asymptotic computational cost.

        \section[Sch\"{u}tzenberger-type Symmetry and a Bijection with Involutions]{Sch\"{u}tzenberger-type Symmetry and\\ a Bijection with Involutions}
        \label{sec:ExtendingPS:Involutions}

            We are now in a position to prove that Extended Patience
            Sorting (Algorithm~\ref{alg:ExtendedPSalgorithm}) has the
            same form of symmetry as does the RSK Correspondence
            (Algorithm~\ref{alg:RSKAlgorithm}).

            \begin{Theorem}
            \label{thm:TwiddlePilesInverseResult}
                Let $(R(\sigma), S(\sigma))$ be the insertion and
                recording piles, respectively, formed by applying
                Algorithm~\ref{alg:ExtendedPSalgorithm} to $\sigma \in
                \mathfrak{S}_{n}$.  Then, applying
                Algorithm~\ref{alg:ExtendedPSalgorithm} to the inverse
                permutation $\sigma^{-1}$, one obtains the pair
                $(S(\sigma), R(\sigma))$.
            \end{Theorem}

            \begin{proof}
                Construct $S'(\sigma)$ from $S(\sigma)$ as discussed
                in Example~\ref{eg:PileReflectionExample}, and form
                the $n$ ordered pairs $(r_{i j}, s'_{i j})$, with $i$
                indexing the individual piles and $j$ indexing the
                cards in the $i^{\mathrm{th}}$ piles.  Then these $n$
                points correspond to the diagram of a permutation
                $\tau \in \mathfrak{S}_{n}$.  However, since
                reflecting these points through the line ``$y = x$''
                yields the diagram for $\sigma$, it follows that $\tau
                = \sigma^{-1}$.
            \end{proof}

            Proposition~\ref{thm:TwiddlePilesInverseResult} suggests
            that Extended Patience Sorting is the right generalization
            of Patience Sorting
            (Algorithm~\ref{alg:MallowsPSprocedure}) since we obtain
            the same symmetry property as for the RSK Correspondence
            (Theorem~\ref{thm:RSKSchuetzenbergerSymmetryProperty}).
            Moreover, Proposition~\ref{thm:TwiddlePilesInverseResult}
            also implies that there is a bijection between involutions
            and pile configurations that avoid simultaneously
            containing the symmetric sub-pile patterns corresponding
            to the patterns given in
            Definition~\ref{defn:StablePairs}.  This corresponds to
            the reverse patience word for a pile configuration
            simultaneously avoiding a symmetric pair of the
            generalized patterns $31\mathrm{-}2$ and $32\mathrm{-}1$,
            etc.  As such, it is interesting to compare this
            construction to the following results obtained by Claesson
            and Mansour \cite{refCM2005}:

            \begin{enumerate}
                \item The size of $S_{n}(3\mathrm{-}12, 3\mathrm{-}21)$ is
                equal to the number of involutions $|\mathfrak{I}_{n}|$ in
                $\mathfrak{S}_{n}$.\smallskip

                \item The size of $S_{n}(31\mathrm{-}2, 32\mathrm{-}1)$ is
                $2^{n-1}$.
            \end{enumerate}

            \noindent The first result suggests that there should be a
            way to relate the result in
            Theorem~\ref{thm:ExtendedPSbijection} to simultaneous
            avoidance of the similar looking patterns $3\mathrm{-}12$
            and $3\mathrm{-}21$.  The second result suggests that
            restricting to complete avoidance of all simultaneous
            occurrences of $31\mathrm{-}2$ and $32\mathrm{-}1$ will
            yield a natural bijection between $S_{n}(31\mathrm{-}2,
            32\mathrm{-}1)$ and a subset $\mathfrak{N} \subset
            \mathfrak{P}_{n}$ such that $\mathfrak{N} \cap
            \mathfrak{P}_{n}(\gamma)$ contains exactly one pile
            configuration of each shape $\gamma$.  A natural family
            for this collection of pile configurations consists of
            what we call \emph{non-crossing pile configurations};
            namely, for the composition $\gamma = (\gamma_{1},
            \gamma_{2}, \ldots, \gamma_{k}) \models n$,
            \[
                \mathfrak{N} \cap
                \mathfrak{P}_{n}(\gamma) = \{\{ \{\gamma_{1} > \cdots > 1\},
                \{\gamma_{1} + \gamma_{2} > \cdots > \gamma_{1} + 1 \},
                \ldots, \{n > \cdots > n - \gamma_{k-1}\}\}\}
            \]
            so that there are exactly $2^{n-1}$ such pile
            configurations.  One can also show that $\mathfrak{N}$ is
            the image $R(S_{n}(3\mathrm{-}1\mathrm{-}2))$ of all
            permutations avoiding the classical pattern
            $3\mathrm{-}1\mathrm{-}2$ under the Patience Sorting map
            $R:\mathfrak{S}_{n} \to \mathfrak{P}_{n}$.
        
        \section[Geometric Form for Extended Patience Sorting]{A Geometric Form for the Extended Patience Sorting Algorithm}
        \label{sec:ExtendingPS:GeometricPS}

            Viennot introduced the shadow diagram of a permutation
            while studying Sch\"{u}t\-zen\-ber\-ger Symmetry for the RSK
            Correspondence
            (Theorem~\ref{thm:RSKSchuetzenbergerSymmetryProperty},
            which was first proven using a direct combinatorial
            argument in \cite{refSchutzenberger1963}).  Specifically,
            using a particular labelling of the constituent ``shadow
            lines'' in recursively defined shadow diagrams, one
            recovers successive rows in the usual RSK insertion and
            recording tableaux.  Sch\"{u}tzenberger Symmetry for the
            RSK Correspondence then immediately follows since
            reflecting these shadow diagrams through the line ``$y =
            x$'' both inverts the permutation and exactly interchanges
            these labellings.

            We review Viennot's Geometric RSK Algorithm in
            Section~\ref{sec:ExtendingPS:GeometricPS:NEshadows} below.
            Then, in
            Section~\ref{sec:ExtendingPS:GeometricPS:SWshadows}, we
            define a natural dual to Viennot's construction that
            similarly produces a geometric characterization for
            Extended Patience Sorting.  As with the RSK
            Correspondence, the analog of Sch\"{u}tzenberger Symmetry
            follows as an immediate consequence.  Unlike Geometric
            RSK, though, the lattice paths formed under Geometric
            Patience Sorting are allowed to intersect.  Thus, having
            defined these two algorithms, we classify in
            Section~\ref{sec:ExtendingPS:GeometricPS:TypesOfCrossings}
            the types of intersections that can occur under Geometric
            Patience Sorting and then characterize them in
            Section~\ref{sec:ExtendingPS:GeometricPS:CharacterizingCrossings}.
        
            \subsection[Northeast Shadowlines and Geometric RSK]{Northeast Shadowlines and Geometric RSK}
            \label{sec:ExtendingPS:GeometricPS:NEshadows}

            In this section, we briefly review Viennot's geometric
            form for the RSK Correspondence in order to motivate the
            geometric form for Extended Patience Sorting that is given
            in Section~\ref{sec:ExtendingPS:GeometricPS:SWshadows}.
            (Viennot's Geometric RSK Algorithm was first introduced in
            \cite{refViennot1977}; an English version can also be
            found in \cite{refViennot1984} and in
            \cite{refSagan2000}.)

            In Section~\ref{sec:PSasAlgorithm:PilesAndNEshadows}, the
            northeast shadow diagram for a collection of lattice
            points was defined by inductively taking northeast
            shadowlines for those lattice points not employed in
            forming the previous shadowlines.  In particular, given a
            permutation $\sigma \in \mathfrak{S}_{n}$, we form the
            northeast shadow diagram $D_{NE}(\sigma) =
            \{L_{1}(\sigma), \ldots, L_{k}(\sigma)\}$ for $\sigma$ by
            first forming the northeast shadowline $L_{1}(\sigma)$ for
            $\{(1, \sigma_{1}), (2, \sigma_{2}), \ldots, (n,
            \sigma_{n})\}$.  Then we ignore the lattice points whose
            northeast shadows were used in building $L_{1}(\sigma)$
            and define $L_{2}(\sigma)$ to be the northeast shadowline
            of the resulting subset of the permutation diagram.  We
            then take $L_{3}(\sigma)$ to be the northeast shadowline
            for the points not yet used in constructing either
            $L_{1}(\sigma)$ or $L_{2}(\sigma)$, and this process
            continues until all points in the permutation diagram are
            exhausted.

            We can characterize the points whose shadows define the
            shadowlines at each stage of this process as follows: they
            are the smallest collection of unused points whose shadows
            collectively contain all other remaining unused points
            (and hence also contain the shadows of those points).  As
            a consequence of this shadow containment property, the
            shadowlines in a northeast shadow diagram will never
            cross.  However, as we will see in
            Section~\ref{sec:ExtendingPS:GeometricPS:SWshadows} below,
            the dual construction to the definition of northeast
            shadow diagrams will allow for crossing shadowlines, which
            are then classified and characterized in
            Section~\ref{sec:ExtendingPS:GeometricPS:TypesOfCrossings}
            and
            \ref{sec:ExtendingPS:GeometricPS:CharacterizingCrossings},
            respectively.  This distinction results from the reversal
            of the above shadow containment property.  
              
            As simple as northeast shadowlines were to define in
            Section~\ref{sec:PSasAlgorithm:PilesAndNEshadows}, a great
            deal of information can still be gotten from them.  One of
            the most basic properties of the northeast shadow diagram
            $D_{NE}^{(0)}(\sigma)$ for a permutation $\sigma \in
            \mathfrak{S}_{n}$ is that it encodes the top row of the
            RSK insertion tableau $P(\sigma)$ (resp.~recording tableau
            $Q(\sigma)$) as the smallest ordinates (resp.  smallest
            abscissae) of all points belonging to the constituent
            shadowlines $L_{1}(\sigma), L_{2}(\sigma), \ldots,
            L_{k}(\sigma)$.  One proves this by comparing the use of
            Schensted Insertion on the top row of the insertion
            tableau with the intersection of vertical lines having the
            form ``$x = a$''.  In particular, as $a$ increases from
            $0$ to $n$, the line ``$x = a$'' intersects the lattice
            points in the permutation diagram in the order that they
            are inserted into the top row, and so shadowlines connect
            elements of $\sigma$ to those smaller elements that will
            eventually bump them.  (See Sagan \cite{refSagan2000} for
            more details.)  

            \begin{figure}[t]\label{fig:GeometricRSKexample}
                  \centering
                  \begin{tabular}{cc}

                      \begin{minipage}[c]{2.75in}
                          \centering
                              \psset{xunit=0.125in,yunit=0.125in}

                              \begin{pspicture}(0,0)(9,9)

                                  \psaxes[Dy=2]{->}(9,9)

                                  \psline[linecolor=black,linewidth=1pt](5,9)(5,8)(6,8)(6,7)(8,7)(8,3)(9,3)

                                  \psline[linecolor=black,linewidth=1pt](3,9)(3,5)(7,5)(7,2)(9,2)

                                  \psline[linecolor=black,linewidth=1pt](1,9)(1,6)(2,6)(2,4)(4,4)(4,1)(9,1)

                                   \rput(1,6){{\large $\bullet$}}%
                                   \rput(2,4){{\large $\bullet$}}%
                                   \rput(3,5){{\large $\bullet$}}%
                                   \rput(4,1){{\large $\bullet$}}%
                                   \rput(5,8){{\large $\bullet$}}%
                                   \rput(6,7){{\large $\bullet$}}%
                                   \rput(7,2){{\large $\bullet$}}%
                                   \rput(8,3){{\large $\bullet$}}%

                                   \rput(2,6){{\large $\odot$}}%
                                   \rput(4,4){{\large $\odot$}}%
                                   \rput(7,5){{\large $\odot$}}%
                                   \rput(6,8){{\large $\odot$}}%
                                   \rput(8,7){{\large $\odot$}}%

                              \end{pspicture}\\[0.25in]

                            (a) Salient points for $D_{NE}^{(0)}(64518723)$.\\[0.325in]
                      \end{minipage}
                      &
                      \begin{minipage}[c]{2.75in}
                          \centering
                              \psset{xunit=0.125in,yunit=0.125in}

                              \begin{pspicture}(0,0)(9,9)

                                  \psaxes[Dy=2]{->}(9,9)

                                  \psline[linecolor=black,linewidth=1pt](9,7)(8,7)(8,9)

                                  \psline[linecolor=black,linewidth=1pt](6,9)(6,8)(7,8)(7,5)(9,5)

                                  \psline[linecolor=black,linewidth=1pt](2,9)(2,6)(4,6)(4,4)(9,4)

                                   \rput(2,6){{\large $\bullet$}}%
                                   \rput(4,4){{\large $\bullet$}}%
                                   \rput(7,5){{\large $\bullet$}}%
                                   \rput(6,8){{\large $\bullet$}}%
                                   \rput(8,7){{\large $\bullet$}}%

                              \end{pspicture}\\[0.25in]

                          (b) Shadow Diagram $D_{NE}^{(1)}(64518723)$.
                      \end{minipage}
                      \\
                      & \vspace{0.5cm}
                      \\
                      \begin{minipage}[c]{2.75in}
                          \centering

                              \psset{xunit=0.125in,yunit=0.125in}

                              \begin{pspicture}(0,0)(9,9)

                                  \psaxes[Dy=2]{->}(9,9)

                                  \psline[linecolor=black,linewidth=1pt](9,7)(8,7)(8,9)

                                  \psline[linecolor=black,linewidth=1pt](6,9)(6,8)(7,8)(7,5)(9,5)

                                  \psline[linecolor=black,linewidth=1pt](2,9)(2,6)(4,6)(4,4)(9,4)

                                   \rput(2,6){{\large $\bullet$}}%
                                   \rput(4,4){{\large $\bullet$}}%
                                   \rput(7,5){{\large $\bullet$}}%
                                   \rput(6,8){{\large $\bullet$}}%
                                   \rput(8,7){{\large $\bullet$}}%

                                   \rput(4,6){{\large $\odot$}}%
                                   \rput(7,8){{\large $\odot$}}%

                              \end{pspicture}\\[0.25in]

                            (c) Salient points for $D_{NE}^{(1)}(64518723)$.
                      \end{minipage}
                      &
                      \begin{minipage}[c]{2.75in}
                          \centering

                              \psset{xunit=0.125in,yunit=0.125in}

                              \begin{pspicture}(0,0)(9,9)

                                  \psaxes[Dy=2]{->}(9,9)

                                  \psline[linecolor=black,linewidth=1pt](7,9)(7,8)(9,8)

                                  \psline[linecolor=black,linewidth=1pt](4,9)(4,6)(9,6)

                                   \rput(4,6){{\large $\bullet$}}%
                                   \rput(7,8){{\large $\bullet$}}%

                              \end{pspicture}\\[0.25in]

                          (d) Shadow Diagram $D_{NE}^{(2)}(64518723)$.
                      \end{minipage}

                  \end{tabular}\\

                  \caption{The northeast shadow diagrams for the permutation $64518723 \in \mathfrak{S}_{8}$.}

              \end{figure}

            Remarkably, one can then use the southwest corners (called
            the \emph{salient points}) of $D_{NE}^{(0)}(\sigma)$ to
            form a new shadow diagram $D_{NE}^{(1)}(\sigma)$ that
            similarly gives the second rows of $P(\sigma)$ and
            $Q(\sigma)$.  Then, inductively, the salient points of
            $D_{NE}^{(1)}(\sigma)$ can be used to give the third rows
            of $P(\sigma)$ and $Q(\sigma)$, and so on.  As such, one
            can view this recursive formation of shadow diagrams as a
            geometric form for the RSK correspondence.  We illustrate
            this process in Figure~4.1 for the following permutation:
            \begin{displaymath}
                \sigma =  64518723 \stackrel{RSK}{\longleftrightarrow}
                \left(~
                    \young(123,457,68)
                    \raisebox{-0.5cm}{,}\
                    \young(135,268,47)
                ~\right)
            \end{displaymath}
            
            \subsection[Southwest Shadowlines and Geometric Patience Sorting]{Southwest Shadowlines and\\ Geometric Patience Sorting}
            \label{sec:ExtendingPS:GeometricPS:SWshadows}

            In this section, we introduce a natural dual to Viennot's
            Geometric RSK construction as given in
            Section~\ref{sec:ExtendingPS:GeometricPS:NEshadows}.  We
            begin with the following fundamental definition.

            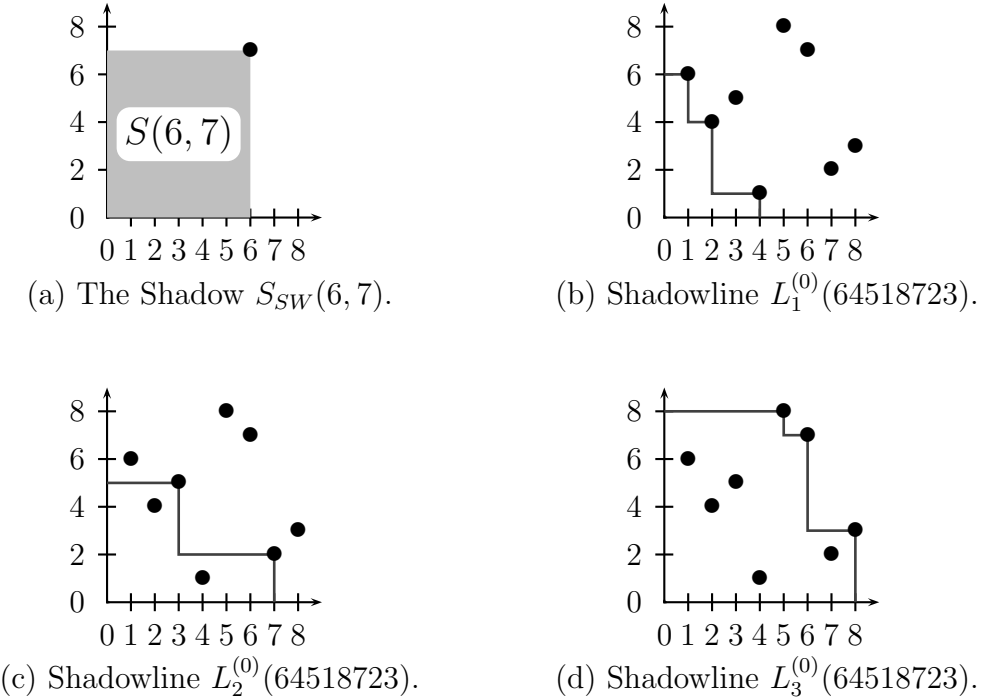
\begin{figure}[t]
                \label{fig:PSshadowExample}
                \centering

                \begin{tabular}{cccc}

                    \begin{minipage}[c]{2.75in}
                        \centering

                            \psset{xunit=0.125in,yunit=0.125in}

                            \begin{pspicture}(0,0)(9,9)

                                \psaxes[Dy=2]{->}(9,9)

                                \psframe[fillcolor=lightgray,fillstyle=solid,linecolor=lightgray](6,7)(0,0)

                                \rput(3,3.5){\psframebox*[framearc=.5]{{\large
                                $S(6,7)$}}}

                                \rput(6,7){{\large $\bullet$}}

                            \end{pspicture}\\[0.25in]

                          (a) The Shadow $S_{SW}(6,7)$.
                    \end{minipage}
                    &
                    \begin{minipage}[c]{2.75in}
                        \centering

                            \psset{xunit=0.125in,yunit=0.125in}

                          \begin{pspicture}(0,0)(9,9)

                              \psaxes[Dy=2]{->}(9,9)
                              \psline[linecolor=darkgray,
                              linewidth=1pt](0,6)(1,6)(1,4)(2,4)(2,1)(4,1)(4,0)

                               \rput(1,6){{\large $\bullet$}}%
                               \rput(2,4){{\large $\bullet$}}%
                               \rput(3,5){{\large $\bullet$}}%
                               \rput(4,1){{\large $\bullet$}}%
                               \rput(5,8){{\large $\bullet$}}%
                               \rput(6,7){{\large $\bullet$}}%
                               \rput(7,2){{\large $\bullet$}}%
                               \rput(8,3){{\large $\bullet$}}%

                          \end{pspicture}\\[0.25in]

                        (b) Shadowline $L_{1}^{(0)}(64518723)$.
                    \end{minipage}
                      \\
                      & \vspace{0.5cm}
                      \\
                    \begin{minipage}[c]{2.75in}
                        \centering

                          \psset{xunit=0.125in,yunit=0.125in}

                          \begin{pspicture}(0,0)(9,9)

                              \psaxes[Dy=2]{->}(9,9)
                              
                              \psline[linecolor=darkgray,
                              linewidth=1pt](0,5)(3,5)(3,2)(7,2)(7,0)

                               \rput(1,6){{\large $\bullet$}}%
                               \rput(2,4){{\large $\bullet$}}%
                               \rput(3,5){{\large $\bullet$}}%
                               \rput(4,1){{\large $\bullet$}}%
                               \rput(5,8){{\large $\bullet$}}%
                               \rput(6,7){{\large $\bullet$}}%
                               \rput(7,2){{\large $\bullet$}}%
                               \rput(8,3){{\large $\bullet$}}%

                          \end{pspicture}\\[0.25in]

                        (c) Shadowline $L_{2}^{(0)}(64518723)$.
                    \end{minipage}
                    &
                    \begin{minipage}[c]{2.75in}
                        \centering

                          \psset{xunit=0.125in,yunit=0.125in}

                          \begin{pspicture}(0,0)(9,9)

                              \psaxes[Dy=2]{->}(9,9)
                              
                              \psline[linecolor=darkgray,
                              linewidth=1pt](0,8)(5,8)(5,7)(6,7)(6,3)(8,3)(8,0)

                               \rput(1,6){{\large $\bullet$}}%
                               \rput(2,4){{\large $\bullet$}}%
                               \rput(3,5){{\large $\bullet$}}%
                               \rput(4,1){{\large $\bullet$}}%
                               \rput(5,8){{\large $\bullet$}}%
                               \rput(6,7){{\large $\bullet$}}%
                               \rput(7,2){{\large $\bullet$}}%
                               \rput(8,3){{\large $\bullet$}}%

                          \end{pspicture}\\[0.25in]

                        (d) Shadowline $L_{3}^{(0)}(64518723)$.
                    \end{minipage}

                \end{tabular}\\

                \caption{Examples of Southwest Shadow and Shadowline
                Constructions}

            \end{figure}

            \begin{Definition}
            \label{defn:PSshadow}
                Given a lattice point $(m, n) \in \mathbb{Z}^{2}$, we
                define the \emph{southwest shadow} of $(m, n)$ to be
                the quarter space
                \[
                    S_{SW}(m, n) = \{ (x, y) \in \mathbb{R}^{2} \ | \
                    x \leq m, \ y \leq n\}.
                \]
            \end{Definition}

            \noindent See Figure~4.2(a) for an example of a point's
            southwest shadow.

            As with their northeast counterparts, the most important use
            of these shadows is in building \emph{southwest shadowlines}.

            \begin{Definition}
            \label{defn:PSshadowline}
                Given lattice points $(m_{1}, n_{1}), (m_{2}, n_{2}),
                \ldots, (m_{k}, n_{k}) \in \mathbb{Z}^{2}$, we define
                their \emph{southwest shadowline} to be the boundary of
                the union of the shadows $S_{SW}(m_{1}, n_{1})$,
                $S_{SW}(m_{2}, n_{2})$, $\ldots$, $S_{SW}(m_{k}, n_{k})$.
            \end{Definition}

            In particular, we wish to associate to each permutation a
            specific collection of southwest shadowlines.  However,
            unlike the northeast case, these shadowlines are allowed
            to cross (as illustrated in Figures~4.2(b)--(d) and
            Figures~4.3(a)--(b)).

            \begin{Definition}
            \label{defn:PSshadowDiagram}
                Given $\sigma = \sigma_{1}\sigma_{2}\cdots\sigma_{n}
                \in \mathfrak{S}_{n}$, the \emph{southwest shadow
                diagram} $D_{SW}^{(0)}(\sigma)$ of $\sigma$ consists
                of the southwest shadowlines $L_{1}^{(0)}(\sigma),
                L_{2}^{(0)}(\sigma), \ldots, L_{k}^{(0)}(\sigma)$
                formed as follows:\medskip

                \begin{itemize}
                    \item $L_{1}^{(0)}(\sigma)$ is the shadowline for
                    those points $(x, y) \in \{(1, \sigma_{1}), (2,
                    \sigma_{2}), \ldots, (n, \sigma_{n})\}$ such that
                    $S_{SW}(x, y)$ does not contain any other lattice
                    points.\medskip

                    \item Then, while at least one of the points $(1,
                    \sigma_{1}), (2, \sigma_{2}), \ldots, (n,
                    \sigma_{n})$ is not contained in the shadowlines
                    $L_{1}^{(0)}(\sigma), L_{2}^{(0)}(\sigma), \ldots,
                    L_{j}^{(0)}(\sigma)$, define
                    $L_{j+1}^{(0)}(\sigma)$ to be the shadowline for
                    the points
                    \[
                        (x, y) \in \{(i, \sigma_{i}) \ | \ i \in [n],
                        (i, \sigma_{i}) \notin \bigcup^{j}_{k=1}
                        L_{k}^{(0)}(\sigma)\}
                    \]
                    such that $S_{SW}(x, y)$ does not contain any other
                    lattice points in the same set.
                \end{itemize}

            \end{Definition}

            In other words, we again define a shadow diagram by
            recursively eliminating certain points in the permutation
            diagram until every point has been used to define a
            shadowline.  Here, however, we are reversing both the
            direction of the shadows and the shadow containment
            property used in the northeast case.  It is in this sense
            that the geometric form for Extended Patience Sorting
            given below can be viewed as ``dual'' to Viennot's
            geometric form for the RSK Correspondence.

            \begin{figure}[t]\label{fig:GeometricPSExample}
                \centering

                \begin{tabular}{cc}

                    \begin{minipage}[c]{2.75in}
                        \centering

                          \psset{xunit=0.125in,yunit=0.125in}

                          \begin{pspicture}(0,0)(9,9)

                              \psline[linecolor=darkgray,
                              linewidth=1pt](0,5)(3,5)(3,2)(7,2)(7,0)

                              \psaxes[Dy=2]{->}(9,9)

                              \psline[linecolor=darkgray,
                              linewidth=1pt](0,6)(1,6)(1,4)(2,4)(2,1)(4,1)(4,0)

                              \psline[linecolor=darkgray,
                              linewidth=1pt](0,5)(3,5)(3,2)(7,2)(7,0)

                              \psline[linecolor=darkgray,
                              linewidth=1pt](0,8)(5,8)(5,7)(6,7)(6,3)(8,3)(8,0)

                              \rput(1,6){{\large $\bullet$}}%
                              \rput(2,4){{\large $\bullet$}}%
                              \rput(3,5){{\large $\bullet$}}%
                              \rput(4,1){{\large $\bullet$}}%
                              \rput(5,8){{\large $\bullet$}}%
                              \rput(6,7){{\large $\bullet$}}%
                              \rput(7,2){{\large $\bullet$}}%
                              \rput(8,3){{\large $\bullet$}}%

                              \rput(1,4){{\large $\odot$}}%
                              \rput(2,1){{\large $\odot$}}%
                              \rput(3,2){{\large $\odot$}}%
                              \rput(5,7){{\large $\odot$}}%
                              \rput(6,3){{\large $\odot$}}%

                          \end{pspicture}\\[0.25in]

                          (a) Salient points for $D_{SW}^{(0)}(64518723)$.
                    \end{minipage}
                    &
                    \begin{minipage}[c]{2.75in}
                        \centering

                          \psset{xunit=0.125in,yunit=0.125in}

                          \begin{pspicture}(0,0)(9,9)

                              \psaxes[Dy=2]{->}(9,9)

                              \psline[linecolor=darkgray,
                              linewidth=1pt](0,4)(1,4)(1,1)(2,1)(2,0)

                              \psline[linecolor=darkgray,
                              linewidth=1pt](0,2)(3,2)(3,0)

                              \psline[linecolor=darkgray,
                              linewidth=1pt](0,7)(5,7)(5,3)(6,3)(6,0)

                               \rput(1,4){{\large $\bullet$}}%
                               \rput(2,1){{\large $\bullet$}}%
                               \rput(3,2){{\large $\bullet$}}%
                               \rput(5,7){{\large $\bullet$}}%
                               \rput(6,3){{\large $\bullet$}}%

                          \end{pspicture}\\[0.25in]

                        (b) Shadow Diagram $D_{SW}^{(1)}(64518723)$.
                    \end{minipage}
                      \\
                      & \vspace{0.5cm}
                      \\
                    \begin{minipage}[c]{2.75in}
                        \centering

                          \psset{xunit=0.125in,yunit=0.125in}

                          \begin{pspicture}(0,0)(9,9)

                              \psaxes[Dy=2]{->}(9,9)

                              \psline[linecolor=darkgray,
                              linewidth=1pt](0,4)(1,4)(1,1)(2,1)(2,0)

                              \psline[linecolor=darkgray,
                              linewidth=1pt](0,2)(3,2)(3,0)

                              \psline[linecolor=darkgray,
                              linewidth=1pt](0,7)(5,7)(5,3)(6,3)(6,0)

                               \rput(1,4){{\large $\bullet$}}%
                               \rput(2,1){{\large $\bullet$}}%
                               \rput(3,2){{\large $\bullet$}}%
                               \rput(5,7){{\large $\bullet$}}%
                               \rput(6,3){{\large $\bullet$}}%

                               \rput(1,1){{\large $\odot$}}%
                               \rput(5,3){{\large $\odot$}}%

                          \end{pspicture}\\[0.25in]

                          (c) Salient points for $D_{SW}^{(1)}(64518723)$.
                    \end{minipage}
                    &
                    \begin{minipage}[c]{2.75in}
                        \centering

                          \psset{xunit=0.125in,yunit=0.125in}

                          \begin{pspicture}(0,0)(9,9)

                              \psaxes[Dy=2]{->}(9,9)

                              \psline[linecolor=darkgray,
                              linewidth=1pt](0,1)(1,1)(1,0)

                              \psline[linecolor=darkgray,
                              linewidth=1pt](0,3)(5,3)(5,0)

                               \rput(1,1){{\large $\bullet$}}%
                               \rput(5,3){{\large $\bullet$}}%

                          \end{pspicture}\\[0.25in]

                        (d) Shadow Diagram $D_{SW}^{(2)}(64518723)$.
                    \end{minipage}

                \end{tabular}\\

                \caption{The southwest shadow diagrams for the permutation $64518723 \in \mathfrak{S}_{8}$.}

            \end{figure}

            As with northeast shadow diagrams, one can also produce a
            sequence
            \[
                D_{SW}(\sigma) = (D_{SW}^{(0)}(\sigma),
                D_{SW}^{(1)}(\sigma), D_{SW}^{(2)}(\sigma), \ldots)
            \]
            of southwest shadow diagrams for a given permutation
            $\sigma \in \mathfrak{S}_{n}$ by recursively applying
            Definition~\ref{defn:PSshadowDiagram} to salient points,
            with the restriction that new shadowlines can only connect
            points that were on the same shadowline in the previous
            iteration.  (The reason for this important distinction
            from Geometric RSK is discussed further in
            Section~\ref{sec:ExtendingPS:GeometricPS:TypesOfCrossings}
            below.)  The salient points in this case are naturally
            defined to be the northeast corner points of a given set
            of shadowlines.  See Figure~4.3 for an example.

            \begin{Definition}
                We call $D_{SW}^{(k)}(\sigma)$ the $k^{\rm th}$
                \emph{iterate} of the \emph{exhaustive shadow diagram}
                $D_{SW}(\sigma)$ for $\sigma \in \mathfrak{S}_{n}$.
            \end{Definition}

            As mentioned above, the resulting sequence of shadow
            diagrams can be used to reconstruct the pair of pile
            configurations given by Extended Patience Sorting
            (Algorithm~\ref{alg:ExtendedPSalgorithm}).  To accomplish
            this, index the cards in a pile configuration using the
            French convention for tableaux (see \cite{refFulton1997})
            so that the row index increases from bottom to top and the
            column index from left to right.  (In other words, we are
            labelling boxes as we would lattice points in the first
            quadrant of $\mathbb{R}^2$).  Then, for a given
            permutation $\sigma \in \mathfrak{S}_{n}$, the elements of
            the $i$th row of the insertion piles $R(\sigma)$
            (resp.~recording piles $S(\sigma)$) are given by the
            largest ordinates (resp.~abscissae) of the shadowlines
            that comprise $D_{SW}^{(i)}$.

            The main difference between this process and Viennot's
            Geometric RSK is that care must be taken to assemble each
            row in its proper order.  Unlike the entries of a standard
            Young tableau, the elements in the rows of a pile
            configuration do not necessarily increase from left to
            right, and they do not have to be contiguous.  As such,
            the components of each row should be recorded in the order
            that the shadowlines are formed.  The rows can then
            uniquely be assembled into a legal pile configuration
            since the elements in the columns of a pile configuration
            must both decrease (when read from bottom to top) and
            appear in the leftmost pile possible.

            To prove that this process works, one argues along the
            same lines as with Viennot's Geometric RSK. In other
            words, one thinks of the shadowlines produced by
            Definition~\ref{defn:PSshadowDiagram} as a visual record
            for how cards are played atop each other under
            Algorithm~\ref{alg:ExtendedPSalgorithm}.  In particular,
            it should be clear that, given a permutation $\sigma \in
            \mathfrak{S}_{n}$, the shadowlines in both of the shadow
            diagrams $D_{SW}^{(0)}(\sigma)$ and $D_{NE}^{(0)}(\sigma)$
            are defined by the same lattice points from the
            permutation diagram for $\sigma$.  By
            Lemma~\ref{lem:ShadowDiagramPileCorrespondence}, the
            points along a given northeast shadowline correspond
            exactly to the elements in some column of $R(\sigma)$ (as
            both correspond to one of the left-to-right minima
            subsequences of $\sigma$).  Thus, by reading the lattice
            points in the permutation diagram in increasing order of
            their abscissae, one can uniquely reconstruct both the
            piles in $R(\sigma)$ and the exact order in which cards
            are added to these piles (which implicitly yields
            $S(\sigma)$).  In this sense, both $D_{SW}^{(0)}(\sigma)$
            and $D_{NE}^{(0)}(\sigma)$ encode the bottom rows of
            $R(\sigma)$ and $S(\sigma)$.
            
            It is then easy to see by induction that the salient points of
            $D_{SW}^{(k-1)}(\sigma)$ yield the $k^{\rm th}$ rows of
            $R(\sigma)$ and $S(\sigma)$, and so this gives the 
            following Theorem.

            \begin{Theorem}
            \label{thm:GeometricPS}
                The process described above for creating a pair of
                pile configurations $(R, S)$ from the Geometric
                Patience Sorting construction yields the same pair of
                pile configurations $(R(\sigma), S(\sigma))$ as
                Extended Patience Sorting
                (Algorithm~\ref{alg:ExtendedPSalgorithm}) applied to
                $\sigma \in \mathfrak{S}_{n}$.
            \end{Theorem}
            
            \subsection[Types of Crossings in Geometric Patience Sorting]{Types of Crossings in Geometric Patience Sorting}
            \label{sec:ExtendingPS:GeometricPS:TypesOfCrossings}

            As discussed in Section~\ref{sec:Intro:Motivation:PS},
            Extended Patience Sorting
            (Algorithm~\ref{alg:ExtendedPSalgorithm}) can be viewed as
            a ``non-bumping'' version of the RSK Correspondence
            (Algorithm~\ref{alg:RSKAlgorithm}) in that cards are
            permanently placed into piles and are covered by other
            cards rather being displaced by them.  In this sense, one
            of the main differences between their geometric
            realizations lies in how and in what order (when read from
            left to right) the salient points of their respective
            shadow diagrams are determined.  In particular, as playing
            a card atop a pre-existing pile under Patience Sorting is
            essentially like non-recursive Schensted Insertion,
            certain particularly egregious ``multiple bumps'' that
            occur under Schensted Insertion prove to be too
            complicated to be properly modeled by the ``static
            insertions'' of Patience Sorting.

            At the same time, it is also easy to see that, for a given
            $\sigma \in \mathfrak{S}_{n}$, the cards atop the piles in
            the pile configurations $R(\sigma)$ and $S(\sigma)$ (as
            given by Algorithm~\ref{alg:ExtendedPSalgorithm}) are
            exactly the cards in the top rows of the RSK insertion
            tableau $P(\sigma)$ and recording tableau $Q(\sigma)$,
            respectively.  Thus, this raises the question of when the
            remaining rows of $P(\sigma)$ and $Q(\sigma)$ can likewise
            be recovered from $R(\sigma)$ and $S(\sigma)$.  While this
            appears to be directly related to the order in which
            salient points are read (as illustrated in
            Example~\ref{eg:2431RSKvsPSexample} below), one would
            ultimately hope to characterize the answer in terms of
            generalized pattern avoidance similar to the description
            of reverse patience words for pile configurations (as
            given in Section~\ref{sec:ExtendingPS:StablePairs}).

            \begin{Example}
            \label{eg:2431RSKvsPSexample}
                Consider the northeast and southwest shadow diagrams

                \begin{center}

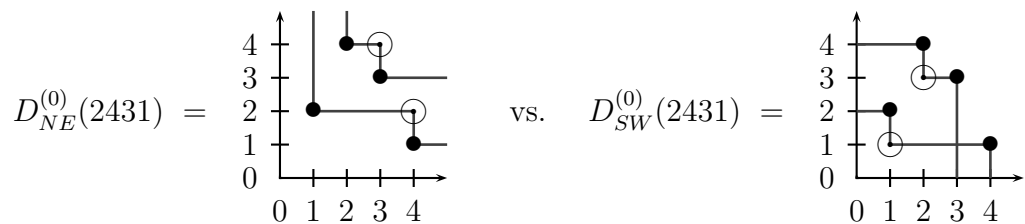
                    \begin{tabular}{ccccc}

                        \raisebox{-.25cm}{$D_{NE}^{(0)}(2431) \ = $}
                        &
                        \begin{minipage}[c]{1.25in}
                            \centering

                              \psset{xunit=0.175in,yunit=0.175in}

                              \begin{pspicture}(0,0)(5,5)

                                  \psaxes{->}(5,5)

                                  \psline[linecolor=darkgray,
                                  linewidth=1pt](1,5)(1,2)(4,2)(4,1)(5,1)

                                  \psline[linecolor=darkgray,
                                  linewidth=1pt](2,5)(2,4)(3,4)(3,3)(5,3)

                                  \rput(1,2){{\large $\bullet$}}%
                                  \rput(2,4){{\large $\bullet$}}%
                                  \rput(3,3){{\large $\bullet$}}%
                                  \rput(4,1){{\large $\bullet$}}%

                                  \rput(4,2){{\large $\odot$}}%
                                  \rput(3,4){{\large $\odot$}}%

                              \end{pspicture}
                        \end{minipage}
                        &
                        \raisebox{-.25cm}{vs.~\quad $D_{SW}^{(0)}(2431) \ = $}
                        &
                        \begin{minipage}[c]{1.25in}
                            \centering

                              \psset{xunit=0.175in,yunit=0.175in}

                              \begin{pspicture}(0,0)(5,5)

                                  \psaxes{->}(5,5)

                                  \psline[linecolor=darkgray,
                                  linewidth=1pt](0,2)(1,2)(1,1)(4,1)(4,0)

                                  \psline[linecolor=darkgray,
                                  linewidth=1pt](0,4)(2,4)(2,3)(3,3)(3,0)

                                  \rput(1,2){{\large $\bullet$}}%
                                  \rput(2,4){{\large $\bullet$}}%
                                  \rput(3,3){{\large $\bullet$}}%
                                  \rput(4,1){{\large $\bullet$}}%

                                  \rput(1,1){{\large $\odot$}}%
                                  \rput(2,3){{\large $\odot$}}%

                              \end{pspicture}
                        \end{minipage}\\[.65in]

                    \end{tabular}.
                \end{center}

                \noindent In particular, note that the order in which
                the salient points are formed (when read from left to
                right) is reversed.  Such reversals serve to
                illustrate one of the inherent philosophical
                differences between the RSK Correspondence and the
                Extended Patience Sorting.

                As previously mentioned, another fundamental
                difference between Geometric RSK and Geometric
                Patience Sorting is that the latter allows certain
                crossings to occur in the lattice paths formed during
                the same iteration of the algorithm.  We classify
                these crossings below and then characterize those
                permutations that yield entirely non-intersecting
                lattice paths in
                Section~\ref{sec:ExtendingPS:GeometricPS:CharacterizingCrossings}.
                
            \end{Example}

            Given $\sigma \in \mathfrak{S}_{n}$, we can classify the
            basic types of crossings in $D_{SW}^{(0)}(\sigma)$ as
            follows: First note that each southwest shadowline in
            $D_{SW}^{(0)}(\sigma)$ corresponds to a pair of decreasing
            sequences of the same length, namely a column from the
            insertion piles $R(\sigma)$ and its corresponding column
            from the recording piles $S(\sigma)$.  Then, given two
            different pairs of such columns in $R(\sigma)$ and
            $S(\sigma)$, the shadowline corresponding to the rightmost
            (resp.~leftmost) pair --- under the convention that new
            columns are always added to the right of all other columns
            in Algorithm~\ref{alg:ExtendedPSalgorithm} --- is called
            the \emph{upper} (resp.~\emph{lower}) shadowline.  More
            formally:

            \begin{Definition}
                Given shadowlines $L^{(m)}_{i}(\sigma),
                L^{(m)}_{j}(\sigma)\in D_{SW}^{(m)}(\sigma)$ with
                $i<j$, we call $L^{(m)}_{i}(\sigma)$ the \emph{lower}
                shadowline and $L^{(m)}_{j}(\sigma)$ the \emph{upper}
                shadowline.  Moreover, should $L^{(m)}_{i}(\sigma)$
                and $L^{(m)}_{j}(\sigma)$ intersect, then we call this
                a \emph{vertical crossing} (resp.~\emph{horizontal
                crossing}) if it involves a vertical
                (resp.~horizontal) segment of $L^{(m)}_{j}(\sigma)$.
            \end{Definition}

             We illustrate these crossings in the following example.
             In particular, note that the only permutations $\sigma
             \in \mathfrak{S}_{3}$ of length three having
             intersections in their $0^{\textrm{th}}$ iterate shadow
             diagram $D_{SW}^{(0)}(\sigma)$ are $312, 231 \in
             \mathfrak{S}_{3}$.

            \begin{figure}[t]\label{fig:PSvsRSKexamples}
                \centering

                \begin{tabular}{cc}

                    \begin{minipage}[c]{2.75in}
                        \centering

                          \psset{xunit=0.125in,yunit=0.125in}

                          \begin{pspicture}(0,0)(9,9)

                              \psaxes[dx=2,dy=2]{->}(9,9)

                              \psline[linecolor=darkgray,
                              linewidth=1pt](0,6)(2,6)(2,2)(4,2)(4,0)

                              \psline[linecolor=darkgray,
                              linewidth=1pt](0,4)(6,4)(6,0)

                              \rput(2,6){{\large $\bullet$}}%
                              \rput(4,2){{\large $\bullet$}}%
                              \rput(6,4){{\large $\bullet$}}%

                          \end{pspicture}\\[0.25in]

                          (a) Shadow Diagram $D_{SW}^{(0)}(312)$.
                    \end{minipage}
                    &
                    \begin{minipage}[c]{2.75in}
                        \centering

                          \psset{xunit=0.125in,yunit=0.125in}

                          \begin{pspicture}(0,0)(9,9)

                              \psaxes[dx=2,dy=2]{->}(9,9)

                              \psline[linecolor=darkgray,
                              linewidth=1pt](0,4)(2,4)(2,2)(6,2)(6,0)

                              \psline[linecolor=darkgray,
                              linewidth=1pt](0,6)(4,6)(4,0)

                              \rput(2,4){{\large $\bullet$}}%
                              \rput(4,6){{\large $\bullet$}}%
                              \rput(6,2){{\large $\bullet$}}%

                          \end{pspicture}\\[0.25in]

                          (b) Shadow Diagram $D_{SW}^{(0)}(231)$.
                    \end{minipage}
                      \\
                      & \vspace{0.5cm}
                      \\
                    \begin{minipage}[c]{2.75in}
                        \centering

                          \psset{xunit=0.125in,yunit=0.125in}

                          \begin{pspicture}(0,0)(9,9)

                              \psaxes[dx=2,dy=2]{->}(9,9)

                              \psline[linecolor=darkgray,
                              linewidth=1pt](0,8)(2,8)(2,4)(4,4)(4,2)(8,2)(8,0)

                              \psline[linecolor=darkgray,
                              linewidth=1pt](0,6)(6,6)(6,0)

                              \psline[linecolor=lightgray, linestyle=solid,
                              linewidth=1pt](0,4)(2,4)(2,2)(4,2)(4,0)

                              \rput(2,8){{\large $\bullet$}}%
                              \rput(4,4){{\large $\bullet$}}%
                              \rput(6,6){{\large $\bullet$}}%
                              \rput(8,2){{\large $\bullet$}}%

                          \end{pspicture}\\[0.25in]

                          (c) {\small Shadow Diagrams} $D_{SW}^{(0)},
                          D_{SW}^{(1)}(4231)$.
                    \end{minipage}
                    &
                    \begin{minipage}[c]{2.875in}
                        \centering

                          \psset{xunit=0.175in,yunit=0.175in}

                          \begin{pspicture}(0,0)(6,6)

                              \psaxes{->}(6,6)

                              \psline[linecolor=darkgray,
                              linewidth=1pt](0,4)(1,4)(1,3)(3,3)(3,1)(4,1)(4,0)

                              \psline[linecolor=darkgray,
                              linewidth=1pt](0,5)(2,5)(2,2)(5,2)(5,0)

                              \psline[linecolor=lightgray, linestyle=solid,
                              linewidth=1pt](0,3)(1,3)(1,1)(3,1)(3,0)

                              \psline[linecolor=lightgray, linestyle=dashed,
                              linewidth=1pt](0,2)(2,2)(2,0)

                              \rput(1,4){{\large $\bullet$}}%
                              \rput(2,5){{\large $\bullet$}}%
                              \rput(3,3){{\large $\bullet$}}%
                              \rput(4,1){{\large $\bullet$}}%
                              \rput(5,2){{\large $\bullet$}}%

                          \end{pspicture}\\[0.25in]

                          (d) {\small Shadow Diagrams} $D_{SW}^{(0)},
                          D_{SW}^{(1)}(45312)$.
                    \end{minipage}

                \end{tabular}\\

                \caption{Shadow diagrams with different types of crossings.}

            \end{figure}

             \begin{Example}
                 ~
                 \begin{enumerate}
                     \item The smallest permutation for which
                     $D_{SW}^{(0)}(\sigma)$ contains a horizontal
                     crossing is $\sigma = 312$ as illustrated in
                     Figure~4.4(a).  The upper
                     shadowline involved in this crossing is the one
                     with only two segments.\smallskip

                     \item The smallest permutation for which
                     $D_{SW}^{(0)}(\sigma)$ has a vertical
                     crossing is $\sigma = 231$ as illustrated in
                     Figure~4.4(b).  As in part
                     (1), the upper shadowline involved in this
                     crossing is again the one with only two
                     segments.\smallskip

                     \item Consider $\sigma = 4231 \in
                     \mathfrak{S}_{4}$.  From Figure~4.4(c),
                     $D_{SW}^{(0)}(\sigma)$ contains exactly two
                     southwest shadowlines, and these shadowlines form
                     a horizontal crossing followed by a vertical
                     crossing.  We call a configuration like this a
                     ``polygonal crossing.''  Note, in particular,
                     that $D_{SW}^{(1)}(\sigma)$ (trivially) has no
                     crossings.\smallskip

                     \item Consider $\sigma = 45312 \in
                     \mathfrak{S}_{5}$.  From Figure~4.4(d),
                     $D_{SW}^{(0)}(\sigma)$ not only has a ``polygonal
                     crossing'' (as two shadowlines with a vertical
                     crossing followed by a horizontal one) but
                     $D_{SW}^{(1)}(\sigma)$ does as well.\smallskip

                 \end{enumerate}
             \end{Example}

             \noindent Polygonal crossings are what make it
             necessary to read only the salient points along the
             same shadowline in the order in which shadowlines are
             formed (as opposed to constructing the subsequent
             shadowlines using the entire partial permutation of
             salient points as in Viennot's Geometric RSK).

             \begin{Example}
                 Consider the shadow diagram of $\sigma = 45312 \in
                 \mathfrak{S}_{5}$ as illustrated in Figure~4.4(d).
                 The $0^{\textrm{th}}$ iterate shadow diagram
                 $D_{SW}^{(0)}$ contain a polygonal crossing, and so
                 the $1^{\textrm{st}}$ iterate shadow diagram
                 $D_{SW}^{(1)}$ needs to be formed as indicated in
                 order to properly describe the pile configurations
                 $R(\sigma)$ and $S(\sigma)$ since
                 
                 {\singlespacing
                 \begin{displaymath}
                     \sigma =  45312 \stackrel{XPS}{\longleftrightarrow}
                     \left(~
                         \begin{minipage}[c]{32pt}
                             $\begin{array}{cc}
                                 1 &   \\
                                 3 & 2  \\
                                 4 & 5
                             \end{array}$
                         \end{minipage}
                         \raisebox{-0.5cm}{,}\
                         \begin{minipage}[c]{32pt}
                             $\begin{array}{cc}
                                 1 &   \\
                                 3 & 2  \\
                                 4 & 5
                             \end{array}$
                         \end{minipage}
                     ~\right)
                 \end{displaymath}
                 }
                 
                 \noindent under Extended Patience Sorting.
             \end{Example}

            \subsection[Characterizing Crossings in Geometric Patience Sorting]{Characterizing Crossings in\\ Geometric Patience Sorting}
            \label{sec:ExtendingPS:GeometricPS:CharacterizingCrossings}

            Unlike the rows of standard Young tableaux, the values in
            the rows of a pile configuration need not increase when
            read from left to right.  As we show below, descents in
            the rows of pile configurations are closely related to the
            crossings given by Geometric Patience Sorting.
            
            As noted in
            Section~\ref{sec:ExtendingPS:GeometricPS:SWshadows} above,
            Geometric Patience Sorting is ostensibly simpler than
            Geometric RSK in that one can essentially recover both the
            insertion piles $R(\sigma)$ and the recording piles
            $S(\sigma)$ from the $0^{\textrm{th}}$ iterate shadow
            diagram $D_{SW}^{(0)}$.  The fundamental use, then, of the
            iterates $D_{SW}^{(i+1)}, D_{SW}^{(i+2)}, \ldots$ is in
            understanding the intersections in the $i^{\textrm{th}}$
            iterate shadow diagram $D_{SW}^{(i)}$.  In particular,
            each shadowline $L^{(m)}_i(\sigma)\in
            D_{SW}^{(m)}(\sigma)$ corresponds to the pair of segments
            of the $i^{\rm th}$ columns of $R(\sigma)$ and $S(\sigma)$
            that are above the $m^{\rm th}$ row (or are the $i^{\rm
            th}$ columns if $m=0$), where rows are numbered from
            bottom to top.

            \begin{Theorem}
            \label{thm:NoncrossingPilesCondition}
                Each iterate $D_{SW}^{(m)}(\sigma)$ ($m\ge 0$) of $\sigma \in
                \mathfrak{S}_{n}$ is free from crossings if and only if
                every row in both $R(\sigma)$ and $S(\sigma)$ is monotone
                increasing from left to right.
            \end{Theorem}

            \begin{proof}
                Since each shadowline $L^{(m)}_{i} =
                L^{(m)}_i(\sigma)$ in the shadow diagram
                $D_{SW}^{(m)}(\sigma)$ depends only on the $i^{\rm
                th}$ columns of $R = R(\sigma)$ and $S = S(\sigma)$
                above row $m$, we may assume, without loss of
                generality, that $R$ and $S$ have the same shape with
                exactly two columns.

                Let $m+1$ be the highest row where a descent occurs in
                either $R$ or $S$.  If this descent occurs in $R$,
                then $L^{(m)}_2$ is the upper shadowline in a
                horizontal crossing since $L^{(m)}_2$ has ordinate
                below that of $L^{(m)}_1$, which is the lower
                shadowline in this crossing (as in $312$).  If this
                descent occurs in $S$, then $L^{(m)}_2$ is the upper
                shadowline in a vertical crossing since $L^{(m)}_2$
                has abscissa to the left of $L^{(m)}_1$, which is the
                lower shadowline in this crossing (as in $231$).  Note
                that both types of descents may occur simultaneously
                (as in $4231$ or $45312$).

                Conversely, suppose $m$ is the last iterate at which a
                crossing occurs in $D_{SW}(\sigma)$ (i.e.,
                $D_{SW}^{(\ell)}(\sigma)$ has no crossings for
                $\ell>m$).  We prove that the shadowline $L^{(m)}_2$
                can only form a crossing using its first or last
                segment.  This, in turn, implies that row $m$ in $R$
                or $S$ is decreasing.  Note that a crossing occurs
                when there is a vertex of $L^{(m)}_1$ that is not in
                the shadow of any point of $L^{(m)}_2$.  Thus, we need
                only show that this can only involve the first or last
                vertex.  Let $\{(s_1,r_1), (s_2,r_2),\dots\}$ and
                $\{(u_1,t_1), (u_2,t_2),\dots\}$ be the vertices that
                define $L^{(m)}_1$ and $L^{(m)}_2$, respectively.
                Then $\{r_i\}_{i\ge 1}$ and $\{t_i\}_{i\ge 1}$ are
                decreasing while $\{s_i\}_{i\ge 1}$ and $\{u_i\}_{i\ge
                1}$ are increasing.  Write $(a,b)\le (c,d)$ if $(a,b)$
                is in the shadow of $(c,d)$ (i.e., if $a\le b$ and
                $c\le d$), and consider $L^{(m+1)}_1$ and
                $L^{(m+1)}_2$.  By hypothesis, these are noncrossing
                shadowlines defined by the salient points
                $\{(s_1,r_2), (s_2,r_3),\dots\}$ and $\{(u_1,t_2),
                (u_2,t_3),\dots\}$, respectively.  Moreover, given any
                index $i$, there is an index $j$ such that
                $(s_i,r_{i+1})\le (u_j,t_{j+1})$.  Suppose, in
                particular, that $(s_i,r_{i+1})\le (u_j,t_{j+1})$ and
                $(s_{i+1},r_{i+2})\le (u_k,t_{k+1})$ for some $j<k$.
                Each upper shadowline vertex must contain some lower
                shadowline vertex in its shadow, so, for all indices
                $\ell$ satisfying $j \leq \ell \leq k$, either
                $(s_i,r_{i+1})\le(u_\ell,t_{\ell+1})$ or
                $(s_{i+1},r_{i+2})\le(u_\ell,t_{\ell+1})$.  Let $\ell$
                be the smallest such index such that
                $(s_{i+1},r_{i+2})\le(u_\ell,t_{\ell+1})$.  If
                $(s_i,r_{i+1})\le(u_\ell,t_{\ell+1})$, then
                $(s_{i+1},r_{i+1})\le (u_\ell,t_{\ell+1})\le
                (u_\ell,t_\ell)$.  Similarly, if
                $(s_i,r_{i+1})\nleq(u_\ell,t_{\ell+1})$, then
                $(s_i,r_{i+1})\le(u_{\ell-1},t_\ell)$, from which
                $(s_{i+1},r_{i+1})\le (u_\ell,t_\ell)$.  Thus, in both
                cases, $(s_{i+1},r_{i+1})\le (u_\ell,t_\ell)$, and so
                the desired conclusion follows.
            \end{proof}

            An immediate corollary of the above proof is that each row
            $i$ ($i\geq m$) in both $R(\sigma)$ and $S(\sigma)$ is
            monotone increasing (from left to right) if and only if
            every iterate $D_{SW}^{(i)}(\sigma)$ ($i\geq m$) is free
            from crossings.
                      
            We conclude this section by noting that
            Theorem~\ref{thm:NoncrossingPilesCondition} only
            characterizes the output of the Extended Patience Sorting
            Algorithm.  At the time of this writing, a full description
            of the permutations themselves remains elusive.  We
            nonetheless provide the following theorem as a first step
            toward characterizing those permutations that result in
            non-crossing lattice paths under Geometric Patience Sorting.

            \begin{Theorem}
            \label{thm:CharacterizedZerothCrossings}
                The set $S_n(3\mathrm{-}\bar{1}\mathrm{-}42,
                31\mathrm{-}\bar{4}\mathrm{-}2)$ consists of all reverse
                patience words having non-intersecting shadow diagrams
                (i.e., no shadowlines cross in the $0^{\rm th}$ iterate
                shadow diagram).  Moreover, given a permutation $\sigma
                \in S_n(3\mathrm{-}\bar{1}\mathrm{-}42,
                31\mathrm{-}\bar{4}\mathrm{-}2)$, the values in the
                bottom rows of $R(\sigma)$ and $S(\sigma)$ increase from
                left to right.
            \end{Theorem}

            \begin{proof}
                From Theorem~\ref{thm:23-1and3-1-42Equivalence},
                $R(S_n(3\mathrm{-}\bar{1}\mathrm{-}42,
                31\mathrm{-}\bar{4}\mathrm{-}2)) =
                R(S_n(23\mathrm{-}1,3\mathrm{-}12))$ consists exactly of
                set partitions of $[n] = \{1, 2, \ldots, n\}$ whose
                components can be ordered so that both the minimal and
                maximal elements of the components simultaneously
                increase.  (These are called \emph{strongly monotone
                partitions} in \cite{refCM2005}).

                Let $\sigma\in
                S_n(3\mathrm{-}\bar{1}\mathrm{-}42,31\mathrm{-}\bar{4}\mathrm{-}2)$.
                Since $\sigma$ avoids $3\mathrm{-}\bar{1}\mathrm{-}42$,
                we must have that $\sigma=RPW(R(\sigma))$ by
                Theorem~\ref{thm:23-1and3-1-42Equivalence}.  Thus, the
                $i^{\rm th}$ shadowline $L^{(0)}_{i}(\sigma)$ is the
                boundary of the union of shadows generated by the
                $i^{\text{th}}$ left-to-right minima subsequence $s_{i}$
                of $\sigma$.  In particular, we can write $s_{i} =
                \varsigma_{i}a_{i}$ where $a_{k} > \cdots > a_{2} >
                a_{1}$ form the right-to-left minima subsequence of
                $\sigma$.  Let $b_i$ be the $i^{\rm th}$ left-to-right
                maximum of $\sigma$.  Then $b_i$ is the left-most (i.e.,
                maximal) entry of $\varsigma_{i}a_{i}$, so
                $\varsigma_{i}a_{i} = b_{i}\varsigma'_{i}a_{i}$ for some
                decreasing subsequence $\varsigma'_{i}$.  Note that
                $\varsigma'_{i}$ may be empty so that $b_{i} = a_{i}$.

                Since $b_{i}$ is the $i^{\rm th}$ left-to-right maximum
                of $\sigma$, it must be at the bottom of the $i^{\rm
                th}$ column of $R(\sigma)$.  (Similarly, $a_i$ is at the
                top of the $i^{\rm th}$ column.)  So the bottom rows of
                both $R(\sigma)$ and $S(\sigma)$ must be in increasing
                order.

                Now consider the $i^{\rm th}$ and $j^{\rm th}$ shadowlines
                $L^{(0)}_{i}(\sigma)$ and $L^{(0)}_{j}(\sigma)$ of
                $\sigma$, respectively, where $i<j$.  We have that
                $b_i<b_j$ from which the initial horizontal segment of the
                $i^{\rm th}$ shadowline is lower than that of the $j^{\rm
                th}$ shadowline.  Moreover, $a_i$ is to the left of $b_j$,
                so the remaining segment of the $i^{\rm th}$ shadowline is
                completely to the left of the remaining segment of the
                $j^{\rm th}$ shadowline.  Thus, $L^{(0)}_{i}(\sigma)$ and
                $L^{(0)}_{j}(\sigma)$ do not intersect.
            \end{proof}
            
    %
    %

\end{document}